  \documentstyle[12pt]{article}
\newcommand{\numero}[1]{
\addtocounter{section}{1}
\begin{center}{\bf \thesection .\
#1\vspace{-.1in}}\end{center}
\setcounter{subsection}{0}
\setcounter{lemma}{0}\indent}

\newcommand{\subnumero}[1]{
\pagebreak[1]\begin{center}{\em #1}\nopagebreak\end{center}
}

\newcommand{\eop}{\hfill $/$\hspace*{-.1cm}$/$\hspace*{-.1cm}$/$\vspace{.1in}}

\newtheorem{lemma}{Lemma}[section]
\newtheorem{theorem}[lemma]{Theorem}

\newtheorem{corollary}[lemma]{Corollary}

\newtheorem{proposition}[lemma]{Proposition}

\newcommand{\Cc}{{\cal C}}

\newcommand{\Ff}{{\cal F}}
\newcommand{\Gg}{{\cal G}}

\newcommand{\Bb}{{\cal B}}

\newcommand{\Uu}{{\cal U}}

\newcommand{\Ii}{{\cal I}}
\newcommand{\Jj}{{\cal J}}

\newcommand{\Xx}{{\cal X}}

\newcommand{\tworightarrows}{\stackrel{\displaystyle \rightarrow}{\rightarrow}}

\begin{document}

\section*{A closed model structure for $n$-categories, internal
$\underline{Hom}$, $n$-stacks and generalized Seifert-Van Kampen}

Carlos Simpson\newline
CNRS, UMR 5580, Universit\'e Paul Sabatier, 31062 Toulouse CEDEX, France.

\bigskip

\numero{Introduction}

The purpose of this paper is to develop some additional techniques for the
weak $n$-categories defined by Tamsamani in \cite{Tamsamani} (which he calls
{\em $n$-nerves}).  The goal is to be able to define the internal $Hom(A,B)$ for
two $n$-nerves $A$ and $B$, which should itself be an $n$-nerve.  This in
turn is for defining the $n+1$-nerve $nCAT$ of all $n$-nerves conjectured in
\cite{Tamsamani}, which we can do quite easily once we have an internal $Hom$.
It is essentially clear {\em a priori} that  we cannot just take an internal
$Hom$ on all of the $n$-nerves of Tamsamani, and in fact some simple examples
support this: any strict $n$-category may be considered in an obvious way as an
$n$-nerve  i.e. a presheaf of sets over $\Delta ^n$ satisfying certain
properties, but the morphisms of the resulting presheaves are the same as the
strict morphisms of the original strict $n$-categories; on the other hand
one can see that these strict morphisms are not enough to reflect all of the
``right'' morphisms.
\footnote{The simplest example which shows that the strict morphisms are not
enough is  where $G$ is a group and $V$ an abelian group and we set $A$
equal to the category with one object and group of automorphisms $G$, and
$B$ equal to the strict $n$-category with only one $i$-morphism for $i<n$ and
group $V$ of $n$-automorphisms of the unique $n-1$-morphism; then for $n=1$
the equivalence classes of strict morphisms from $A$ to $B$ are the elements of
$H^1(G,V)$ so we would expect to get $H^n(G,V)$ in general, but for $n>1$ there
are no nontrivial strict morphisms from $A$ to $B$.
}
Our strategy to get around this problem will be based
on the idea of {\em closed model category} \cite{Quillen}. We will construct a
closed model category containing the  $n$-nerves of Tamsamani. Then we can
simply take  as the ``right''  $n$-nerve of morphisms, the internal
$Hom(A,B)$ whenever $A$ and $B$ are
{\em fibrant} objects in the closed model category (all objects will be
cofibrant in our case). This strategy is standard practice for topologists.

As usual, in order to define a closed model category we first have to
enlarge the class of objects under consideration.  Instead of $n$-nerves as
defined by Tamsamani we look at $n$-pre-nerves (i.e. presheaves of sets over the
cartesian product of $n$ copies of the standard simplicial category) which
satisfy the constancy condition---C1 in Tamsamani's definition of
$n$-nerve---and call these {\em $n$-precats} (this notion being in between the
pre-nerves and  nerves of \cite{Tamsamani}, we take a different notation). An
$n$-precat may be interpreted as a presheaf on a certain quotient $\Theta^n$ of
$\Delta ^n$, in particular we obtain a category $PC_n$ of objects closed under
all limits, with internal $Hom$ etc. We follow the method of constructing a
closed model category developed by Jardine-Joyal \cite{Jardine} \cite{Joyal} in
the case of simplicial presheaves.  The cofibrations are essentially just
monomorphisms (however we cannot---and don't---require injectivity for
top-degree morphisms, just as sets or categories with monomorphisms are not
closed model categories \cite{Quillen}).  The main problem is to define a
notion of weak equivalence. Our key construction is the construction of an
$n$-nerve $Cat(A)$ for any $n$-precat $A$, basically by throwing in freely all
of the elements which are required by the definition of nerve \cite{Tamsamani}
(although to make things simpler we use here a definition of nerve modified
slightly to ``easy nerve'').  Then we say that a morphism of $n$-precats
$A\rightarrow B$ is a {\em weak equivalence} if $Cat(A)\rightarrow Cat(B)$
is an exterior equivalence of $n$-nerves in the sense of \cite{Tamsamani}.
The fibrant morphisms are characterized in terms of cofibrations and weak
equivalences by a lifting property, in the same way as in \cite{Jardine}.

One new thing that we obtain in the process of doing this
is the notion of pushout. The category of $n$-precats is closed under direct
limits and in particular under pushouts.  Applying the operation $Cat$ then
gives an {\em $n$-categorical pushout}: if $A\rightarrow B$ and $A\rightarrow C$
are morphisms of $n$-nerves then the categorical pushout is $Cat(B\cup ^AC)$.

The main lemma which we need to prove (Lemma \ref{pushout} below) is---again
just as in \cite{Jardine}---the statement that a pushout by a trivial
cofibration (i.e. a cofibration which is a weak equivalence) is again a trivial
cofibration.  After that the rest of the arguments needed to obtain the
closed model structure are relatively standard following \cite{Jardine} when
necessary.

Once the closed model structure is established, we can go on to define
internal $Hom $ and construct the $n+1$-nerve $nCAT$.  Using these we can, in
principal, define the notion of $n$-stack. Our discussion of $n$-stacks is
still at a somewhat speculative stage in the present version of the paper,
because there are several slightly different notions of a family of
$n$-categories parametrized by a $1$-category $\Xx$ and ideally we would like
to---but don't yet---know that they are all the same (as happens for
$1$-stacks).

The notion of categorical pushout which we developed as a technical
tool actually has a geometric consequence: we obtain a generalized Seifert-Van
Kampen theorem (Theorem \ref{svk} below) for the Poincar\'e $n$-groupoids
$\Pi _n(X)$ of a space $X$ which were defined by Tamsamani (\cite{Tamsamani}
\S 2.3 ff).  If $X$ is covered by open sets $U$ and $V$ then $\Pi _n(X)$ is
equivalent to the category-theoretic pushout of $\Pi _n(U)$ and $\Pi _n(V)$
along $\Pi _n(U\cap V)$.  We define the {\em nonabelian cohomology} of $X$
with coefficients in a fibrant $n$-precat $A$ as $H(X,A):= Hom (\Pi _n(X), A)$.
The generalized Seifert-Van Kampen theorem implies a Mayer-Vietoris statement
for this nonabelian cohomology.

There are many possible approaches to the notion of $n$-category and, without
pretending to be exhaustive, I would like to point out some of the other
possibilities here for comparison.
\newline
---One of the  pioneering works in the search for an algebraic approach to
homotopy of spaces is the notion of $Cat ^n$-groups of Brown and Loday. This
is what is now known as the ``cubical'' approach where the set
of objects can itself have a structure for example of $n-1$-category, so it
isn't
quite the same as the approach we are looking for (commonly called the
``globular'' case).
\newline
---Gordon, Powers and Street have intensively
investigated the cases $n=3$ and $n=4$ \cite{Gordon-Power-Street}, following the
path set out by Benabou for $2$-categories \cite{Benabou}.
\newline
---In \cite{Grothendieck} A. Grothendieck doesn't
seem to have hit upon any actual definition but gives a lot of nice intuition
about $n$-categories.
\newline
---On p. 41 of \cite{Grothendieck} starts a reproduction of a letter
from Grothendieck to Breen dated July 1975, in which Grothendieck acknowledges
having recieved a proposed definition of non-strict $n$-category from Breen, a
definition which according to {\em loc. cit} ``...has certainly the merit of
existing...''. It is not clear whether this proposed construction was ever
worked out.
\newline
---In \cite{Street}, R. Street proposes a definition of weak $n$-category
as a simplicial set satisfying a certain variant of the Kan condition where one
takes into account the directions of arrows.
\newline
---Kapranov and Voevodsky in \cite{Kapranov-Voevodsky}
construct, for a topological space $X$, a ``Poincar\'e $\infty$-groupoid'' which
is a strictly associative $\infty$-groupoid but where the arrows are
invertible only up to equivalence.  This of course raises the question to
know if
strictly associative $n$-categories would be a sufficient class to yield the
correct $n+1$-category $nCAT$. As pointed out in the footnote above, one wonders
in particular whether there is a closed model structure to go along with these
strict $n$-categories.
\newline
---In his recent preprint \cite{Batanin} M. Batanin develops some ideas
towards a definition of weak $\infty$-category based on operads. In the
introduction he mentions a letter from Baez and Dolan to Street dating to
November 29, 1995 which contains some ideas for a definition of weak
$n$-category; and he states that Makkai, Hermida and Power have worked on the
idea contained in this letter.
\newline
 ---M. Rosellen told me in September 1996 that he was
working on a version using the theory of operads (cf \cite{Adams} for example).
Just as our current effort is based on Segal's delooping machine, there should
probably be an $n$-category machine analogous to any of the other various
delooping machines, and in fact the problems are almost identical:  the basic
problem of doing $n$-categories comes down to doing delooping while keeping
track of the non-connected case and not requiring things to be invertible up to
homotopy (cf the last section of \cite{Tamsamani} for some arguments relating
$n$-categories and delooping machines).
\newline
---J. Baez and J. Dolan have developed their theory originating in the letter
refered to above,
a definition of $n$-categories based on
operads,  in a preprint  \cite{BaezDolanIII} of February 1997.  In this preprint
they discuss operads, give their definition of $n$-category and of certain
morphisms of $n$-categories, and define the homotopy category of $n$-categories
which they conjecture to be equivalent to the homotopy category for other
definitions such as the category $Ho-n-Cat$ mentionned in \cite{Tamsamani}.

The main problem which needs to be accomplished in any of these points of view
is to obtain an $n+1$-category (hopefully within the same point of view) $nCAT$
parametrizing the $n$-categories of that point of view.  This is the main thing
we are doing here for Tamsamani's point of view.  As far as I know, the present
one is the first precise construction of the $n+1$-category $nCAT$.

Once several such points of
view are up and running, the comparison problem will be posed: to find an
appropriate way to compare different points of view on $n$-categories and (one
hopes) to say that the various points of view are equivalent and in particular
that the various $n+1$-categories $nCAT$ are equivalent via these comparisons.
It is not actually clear to me what type of general setup one should use for
such a comparison theory, although the first thing to try would be to explore a
theory of ``internal closed model category'',  a closed model
category with internal $Hom$: any reasonable point of view on $n$-categories
should probably yield an internal closed model category $n--C$ (such as the
$PC_n$ we obtain below) and furthermore $nCAT$ should be an object in
$(n+1)--C$.
Comparison between the theories might then be possible using a version of
Quillen's adjoint functor approach \cite{Quillen}.  We give an indication of
how to start on comparison in \S 11 by sketching how to obtain a functor from
any internal closed model category containing $Cat$, to the our closed model
category of $n$-categories.

Having a good theory of $n$-categories should open up the possibility to pursue
any of the several programs such as that outlined by Grothendieck
\cite{Grothendieck}, the generalization to $n$-stacks and $n$-gerbs of the work
of Breen \cite{Breen}, or the program of Baez and Dolan in
topological quantum field theory \cite{BaezDolan}. Once the theory of
$n$-stacks is off the ground this will give an algebraic approach to the
``geometric $n$-stacks'' considered in \cite{geometricN}.

We clarify the pretentions to rigor of the various sections of this paper.
\S\S 2--7 are supposed to be a first version of something precise and correct
(although at the time of this first version I haven't checked all of the details
in a very thorough way). The same holds for \S 9 on Seifert-Van Kampen.
On the
other hand, the discussion of \S 8 on $n$-stacks is
blatantly speculative; and the
discussion of \S 10 on nonabelian cohomology is very incomplete.

{\em Acknowledgements:}
I would specially like to thank Z. Tamsamani and A. Hirschowitz.
This work follows up the definition and original work on $n$-nerves done by Z.
Tamsamani in his thesis \cite{Tamsamani}. Much of what is done in the
present paper was suggested by discussions  with Tamsamani. More recently in
preparing some joint work with A. Hirschowitz on universal geometric $n$-stacks
related to Brill-Noether, Hirschowitz asked repeatedly for an algebraic approach
to $n$-stacks which would be more natural than the approach passing through
presheaves of topological spaces or simplicial presheaves.  The present work
owes much to A. Hirschowitz's questions and suggestions, as well as to his
perseverance  in asking for an algebraic approach to $n$-stacks.

I would also like to thank R. Brown for pointing out the
importance of the notion of push-out and Seifert-Van Kampen, and G.
Maltsiniotis and A. Brugui\`eres for helpful discussions.

{\em Ce papier est dedi\'e \`a Nicole, Chlo\'e et L\'eo.}

\numero{Preliminaries}

Let $\Delta$ be the standard category of ordered finite sets.  Let $\Theta ^n$
be the quotient of the cartesian product $\Delta ^n$ obtained by identifying
all of the objects $(M, 0, M')$ for fixed $M = (m_1,\ldots , m_k)$ and variable
$M'= (m'_1, \ldots , m'_{n-k-1})$.  The object of $\Theta ^n$ corresponding
to the class of $(M,0,M')$ with all $m_i >0$ will be denoted $M$.  The object
$(1, \ldots , 1)$ ($k$ times) will be denoted $1^k$. We  permit concatenation
in our notation for objects, thus $M, m$ denotes the object
$(m_1, \ldots , m_k , m)$ (when this makes sense, that is when $k<n$).
The class of $(0,\ldots , 0)$ will be denoted by $0$.

We give the explicit construction of $\Theta ^n$.  If $M=(m_1,\ldots , m_k)$ and
$M' = (m'_1,\ldots , m'_l)$ then set $M,0$ equal to the concatenation of
$M$ with
$(0,\ldots , 0)$ in $\Delta ^n$ and similarly for $M',0$.  We define an
equivalence relation on morphisms $\varphi = (\varphi _1, \ldots , \varphi _n)$
from $M, 0$ to $M', 0$ by saying $\varphi \sim \varphi '$ whenever
there exists $j$ such that $\varphi _i = \varphi '_i$ for $i\leq j$ and
$\varphi _j: m_j \rightarrow m'_j$ factors through the object $0\in \Delta$
(which is the one-point set).  This equivalence relation is compatible with
composition so we obtain a category $\Theta ^n$ by taking as morphisms the
quotient of the morphisms in $\Delta ^n$ by this equivalence relation. There
is an obvious projection from $\Delta ^n$ to $\Theta ^n$.

We assume familiarity with \cite{Tamsamani}.
An $n$-precat is a presheaf of sets on $\Theta ^n$.
This corresponds to an $n$-prenerve in Tamsamani's notation (i.e. presheaf of
sets on $\Delta ^n$) which satisfies his axiom C1 in the definition of
$n$-nerve.
Let $PC_n$ denote the category of $n$-precats.
An $n$-precat is an {\em $n$-category} (or {\em $n$-nerve} in the notation of
Tamsamani \cite{Tamsamani} which is the sense which we will always assign to the
terminology ``$n$-category'' below) if it satisfies certain additional
conditions
\cite{Tamsamani}.   We give an
easier version of these conditions which we call an {\em easy $n$-category}. We
start with the notion of easy equivalence between  two easy
$n$-categories---this is not circular because the notion of easy $n$-category
will only use the notion of easy equivalence for morphisms of $n-1$-categories.
If $A$ and $B$ are easy $n$-categories then a morphism $f: A\rightarrow B$ (of
$n$-precats, i.e. of presheaves on $\Theta ^n$) is an {\em easy equivalence} if
for all $v\in B_{1^k}$ (called a $k$-morphism of $B$) and all $a,a' \in
A_{1^{k-1}}$ with $s(v)=f(a)$ and $t(v)=f(a')$ and $s(a)=s(a')$ and
$t(a)=t(a')$ (here $s$ and $t$ denote the morphisms ``source'' and ``target''
from  $T_{1^{k}}$ to $T_{1^{k-1}}$ for any $n$-precat $T$), there exists $u\in
A_{1^k}$ with $s(u)= a$ and $t(u)= a'$ and $f(u)=v$.  A {\em marked easy
equivalence} is the data of a morphism $f$ together with choices $u(a,a', v)$ in
every situation as above.

The reader is cautioned that we will still need Tamsamani's notion of
equivalence (which he calls ``\'equivalence ext\'erieure'' \cite{Tamsamani} \S
1.3) for our closed model category structure below. The notion of easy
equivalence is mainly just used when it is an ingredient in the notion of
$n$-category.

Before giving the definition of easy $n$-category we introduce the following
notation.  If $T$ is an $n$-precat then for any $M= (m_1, \ldots , m_k)$
we denote by $T_{M/}$ the $n-k$-precat obtained by restricting $T$ to
the subcategory of objects of $\Theta^n$ of the form $(M,M')$ for variable $M'$.
This differs from the notation of Tamsamani who called this just $T_M$;
our notation with a slash is necessitated by the notation $M$ for objects of
$\Theta ^n$.  (Sorry about these slight notational changes but is is much
easier for us to use $\Theta ^n$ for what will be done below).

With these notations,
an $n$-precat $A$ is an {\em easy $n$-category} if:
\newline
---for each $m$, $A_{m/}$ is an easy $n-1$-category; and
\newline
---the morphisms
$$
A_{m/} \rightarrow A_{1/} \times _{A_0} \ldots \times _{A_0} A_{1/}
$$
are easy equivalences. Note here that $A_0$ is set which is the fiber over the
object $0\in \Theta ^n$ (which exits slightly from our notational convention; it
is the class of objects $0, M'$ of $\Delta ^n$ but here there is no ``$M$'' to
put into the  notation so we put ``$0$'' instead).

A {\em marked easy $n$-category} is an easy $n$-category provided with the
addional data of markings for the $A_{m/}$ and markings for the easy
equivalences going into the definition.  These two conditions amount
(recursively) to saying that we have markings for all of the morphisms of the
form
$$
A_{M,m/} \rightarrow A_{M,1/} \times _{A_M} \ldots \times _{A_M} A_{M,1/}.
$$

The notion of marking as we have defined above actually makes sense
for any $n$-precat, and an $n$-precat with a marking is automatically an easy
$n$-category.  For this reason, arbitrary inverse limits of marked easy
$n$-categories (indexed by systems of morphisms which preserve the markings)
are again marked easy $n$-categories.

Suppose $A$ is an $n$-precat.  We define the {\em marked easy
$n$-category generated by $A$} denoted $Cat(A)$ by
$$
Cat (A)= \lim _{\leftarrow , \Cc} T
$$
where the limit is taken over the category $\Cc$ whose objects are triples
$(T,\mu ,f)$ with $(T, \mu )$ a marked easy  $n$-category ($\mu$ denotes the
marking) and $f: A\rightarrow T$ is a morphism of $n$-precats.  The morphisms
of $\Cc$ are morphisms of $n$-precats (i.e. morphisms of presheaves on
$\Theta ^n$) required to preserve $f$ and the marking $\mu$.  By the principle
given in the previous paragraph, this inverse limit is again a marked easy
$n$-category.

The construction $Cat(A)$ is the key to the rest of what we are going to say.
The description of $Cat(A)$ given above is one of cutting it down to size.
There is also a creative description.  In order to explain this
we first discuss certain push-outs of $n$-precats.  An object of $\Theta ^n$
represents a presheaf (i.e. $n$-precat). If $M$ is an object we denote the
$n$-precat represented by $M$ as $h(M)$.  A morphism of $n$-precats
$h(M)\rightarrow A$ is the same thing as an element of $A_M$.  Note that
direct limits exist in the category of $n$-precats (as in any category of
presheaves). In particular push-outs exist.

We construct the following
standard $n$-precats. Let $M= (m_1, \ldots , m_l)$ with $l\leq n-1$, and
let $m\geq 1$ (although by the remark below we could also restrict to $m\geq
2$).  Let $-1 \leq k \leq n-l-1$. We will state the constructions by universal
properties (although we give an explicit construction later).  Note that these
universal properties admit solutions because we work in the category of
presheaves over a given category $\Theta ^n$ so the necessary limits exist.

Define $\Sigma = \Sigma (M, [m], \langle k,k+1\rangle )$ to be the universal
$n$-precat with $k$-morphisms $a,b$ in $\Sigma _{M,m/}$ (i.e.
$a,b\in \Sigma _{M,m, 1^k}$) and a $k+1$-morphism
$$
v= (v_1, \ldots , v_m) \in (\Sigma _{M,1/} \times _{\Sigma _{M,0}}
\ldots  \times _{\Sigma _{M,0}}\Sigma _{M,1/} )_{1^{k+1}}
$$
such that $s(a)=s(b)$, $t(a)=t(b)$, and such that the images of $a$ and $b$
by the usual map to the product of $\Sigma _{M,1/}$ are $s(v)$ and $t(v)$
respectively.  Note that $h=h(M,m, 1^{k+1})$ is the universal $n$-precat
with a $k+1$-morphism $u$ in $\Sigma _{M,m/}$ (i.e. $u\in \Sigma
_{M,m,1^{k+1}}$).  Note that setting $a$ to $s(u)$, $b$ to $t(u)$ and
$v$ to the image of $u$ by the usual map to the product, we obtain (by the
universal property of $\Sigma$) a morphism
$$
\varphi = \varphi (M, [m], \langle k,k+1\rangle ) :
\Sigma (M, [m], \langle k,k+1\rangle )
\rightarrow h(M,m,1^{k+1})
$$
We will show below that $\varphi$ is a cofibration,
\footnote{
The definition, from \S 3 below, is that a morphism
$A\rightarrow B$ of $n$-precats is a {\em cofibration} if for every $M= (m_1,
\ldots , m_k)$ with $k <n$, the morphism $A_M \rightarrow B_M$ is injective.
}
at the same time giving an
explicit construction of $\Sigma$ as a pushout of representable presheaves.
Before doing that, we mention the modifications to the above definition
necessary for the boundary cases $k=-1$ and $k=n-l-1$.

For $k=-1$, $\Sigma (M, [m], \langle -1,0\rangle )$ is the universal $n$-precat
with an object
$$
v= (v_1, \ldots , v_m) \in \Sigma _{M,1} \times _{\Sigma _{M,0}}
\ldots  \times _{\Sigma _{M,0}}\Sigma _{M,1} .
$$
As $h(M,m,0)$ is the universal $n$-precat with an object $u\in h_{M,m}$ we have
an object $v$ as above for $h$ (the image of $u$ by the usual map)
so we obtain $\Sigma \rightarrow h$.

For $k=n-l-1$,
$\Sigma (M, [m], \langle n-l-1,n-l\rangle )$ is the universal $n$-precat
with $a,b\in \Sigma _{M,m, 1^{n-l-1}}$ such that $s(a)=s(b)$ and $t(a)=t(b)$
and such that $a$ and $b$ map to the same elements of
$$
(\Sigma _{M,1/} \times _{\Sigma _{M,0}}
\ldots  \times _{\Sigma _{M,0}}\Sigma _{M,1/} )_{1^{n-l-1}} .
$$
Note that $h=h(M, m, 1^{n-l})$ is normally speaking not defined because
the length of the multiindex $(M,m,1^{n-l})$ is $n+1$.  Thus we formally
define this $h$ to be equal to  $h(M,m, 1^{n-l-1})$ and take the elements $a=b$
equal to the canonical $n-l-1$-morphism in $h_{M,m/}$.  This gives a morphism
$\Sigma \rightarrow h$.

We will now give an explicit construction of $\Sigma$ and use this to show that
$\Sigma \rightarrow h$ is a cofibration. (The boundary cases will be left to
the reader). In general the universal $n$-precat with a collection of elements
with certain equalities required, is a quotient of the disjoint union of the
representable $n$-precats corresponding to the elements we want, by identifying
pairs of morphisms from the representable $n$-precats corresponding to the
elements which need to be equal.
We do this in several steps. First, the universal $n$-precat
$\Upsilon = \Upsilon (M, [m], 1^k)$ with
element
$$
v= (v_1, \ldots , v_m) \in (\Upsilon _{M,1/} \times _{\Sigma _{M,0}}
\ldots  \times _{\Upsilon _{M,0}}\Upsilon _{M,1/} )_{1^{k}}
$$
is constructed as the quotient of the disjoint union of $m$ copies of
$h(M,1,1^k)$ making $m-1$ identifications over pairs of maps
$$
h(M,1,1^k) \leftarrow h(M) \rightarrow h(M,1,1^k).
$$
This is the same as taking the pushout of the diagram
$$
h(M,1,1^k) \leftarrow h(M) \rightarrow h(M,1,1^k) \leftarrow \ldots
\leftarrow h(M) \rightarrow h(M,1, 1^k).
$$
Now $\Sigma (M, [m], \langle k,k+1\rangle )$ is the quotient of the disjoint
union
$$
h(M, m, 1^k)^a\sqcup  h(M, m, 1^k)^b \sqcup  \Upsilon (M, [m], 1^{k+1})
$$
by the following identifications (the superscripts $a$ and $b$ in the above
notation are there to distinguish between the two components, which correspond
respectively to choosing $a$ and $b$). There are two maps (dual to $s$ and $t$)
$s^{\ast}, t^{\ast}: h(M,m,1^{k-1})\rightarrow h(M,m, 1^k)$, and
we identify over the pairs of morphisms
$$
h(M, m, 1^k)^a \stackrel{s^{\ast}}{\leftarrow}
h(M,m, 1^{k-1})\stackrel{s^{\ast}}{\rightarrow} h(M, m, 1^k)^b
$$
and
$$
h(M, m, 1^k)^a \stackrel{t^{\ast}}{\leftarrow}
h(M,m, 1^{k-1})\stackrel{t^{\ast}}{\rightarrow} h(M, m, 1^k)^b.
$$
Then we also identify over the pairs of maps
$$
h(M,m,1^k)^a\leftarrow \Upsilon (M, [m], 1^{k})
\stackrel{s^{\ast}}{\rightarrow}\Upsilon (M, [m], 1^{k+1})
$$
and
$$
h(M,m,1^k)^b\leftarrow \Upsilon (M, [m], 1^{k})
\stackrel{t^{\ast}}{\rightarrow}\Upsilon (M, [m], 1^{k+1})
$$
where the left maps are induced by the collection of principal morphisms
$1\rightarrow m$. The result of all these identifications is
$\Sigma (M, [m], \langle k,k+1\rangle )$.  This gives an explicit construction
for those wary of just defining things by the universal property.

We now show that the morphism $\varphi : \Sigma \rightarrow h$ is
a cofibration in the sense of \S 3 below, i.e. injective on all levels
except the
top one.

We will say that a diagram of $n$-precats of the form
$$
A\tworightarrows B \rightarrow C
$$
(where the two compositions are the same)
is {\em semiexact} if the morphism from the coequalizer of the two arrows
to $C$ is a cofibration in the sense of \S 3 below.
Our above construction gives $\Sigma$ as the coequalizer of
$$
\Upsilon ' \sqcup  \Upsilon ' \sqcup  h' \sqcup  h' \tworightarrows
\Upsilon \sqcup  h^a \sqcup  h^b
$$
where
$$
\Upsilon =  \Upsilon (M, [m], 1^{k+1}), \;\;\;
\Upsilon ' = \Upsilon (M, [m], 1^{k}),
$$
$$
h^a (\mbox{resp.}\,\,  h^b) := h(M,m,1^k),\;\;\; h' := h(M,m, 1^{k-1}).
$$
We would like to prove that
$$
\Upsilon ' \sqcup  \Upsilon ' \sqcup  h' \sqcup  h' \tworightarrows
\Upsilon \sqcup  h^a \sqcup  h^b \rightarrow h(M, m, 1^{k+1})
$$
is semiexact. To prove this it suffices (by a simple set-theoretic
consideration)  to show that
$$
h' \sqcup  h' \tworightarrows h^a \sqcup  h^b \rightarrow h(M, m, 1^{k+1})
$$
$$
\Upsilon ' \tworightarrows \Upsilon \sqcup  h^a \rightarrow h(M, m, 1^{k+1})
$$
and
$$
\Upsilon ' \tworightarrows \Upsilon \sqcup  h^b \rightarrow h(M, m, 1^{k+1})
$$
are semiexact.

The first statement follows from the claim that for any $M= (m_1, \ldots ,
m_l)$,
$$
h(M) \sqcup  h(M)\tworightarrows h(M, 1) \sqcup  h(M, 1) \rightarrow h(M,1,1)
$$
is semi-exact. Let $P= (p_1, \ldots , p_n)$ be
an element of $\Theta ^n$ (some of the $p_j$ may be zero). A morphism
$P\rightarrow (M, 1)$ corresponds to a collection of morphisms $f_i:p_i
\rightarrow m_i$ for $i\leq l$ and $f_{l+1}:p_{l+1} \rightarrow 1$, up to
equivalence. The equivalence relation is obtained by saying that if one of the
morphisms factors through $0$ then the subsequent ones don't matter.
The first thing to note is that the two morphisms $(I, s), (I,t): h(M,
1)\rightarrow h(M, 1, 1)$ are injective, as follows directly from the above
description.  Now suppose  that $f, g: P\rightarrow (M, 1)$ are two  morphisms
such that $(I,s)\circ f = (I, t)\circ g$.  Since $s: 0\rightarrow 1$ composed
with anything is different from $t: 0\rightarrow 1$ composed with the same, this
means that one of the $f_i$ must factor through $0$ for $i\leq l+1$, and that
$f_j=g_j$ for $j\leq i$. If it is the case for $i\leq l$ then $f$ and $g$ both
come from the morphism $(f_1, \ldots , f_l): P \rightarrow M$, which is
equivalent to $(g_1, \ldots , g_l)$, via either one of the morphisms
$M\rightarrow (M, 1)$. If $i=l+1$ then $f_{l+1}=g_{l+1}$ factors through
one of the two morphisms $0\rightarrow 1$, so $f$ and $g$ both come from the
morphism   $(f_1, \ldots , f_l)=(g_1, \ldots , g_l)$ via the morphism
$M\rightarrow (M, 1)$ corresponding to the morphism $0\rightarrow 1$ occuring
above. Thus $f$ and $g$ are equivalent in the coequalizer, giving the claim for
this paragraph and thus the first of our semiexactness statements.

For the next semiexactness statement, we first note that $\Upsilon \rightarrow
h(M, m, 1^{k+1})$ is cofibrant. In fact we can describe $\Upsilon$ as a
subsheaf of $h(M,m, 1^{k+1})$ as follows. For $P= (p_1 , \ldots , p_n)$
the morphisms from $P$ to $(M,m, 1^{k+1})$ are the sequences of morphisms
$f= (f_1, \ldots , f_{l+k +2})$ with $f_i : p_i \rightarrow m_i$ (or taking
values in $m$ or $1$ as appropriate depending on $i$).  Such a morphism is
contained in $\Upsilon _P$ if and only if the morphism
$f_{l+1}: p_{l+1}\rightarrow m$ factors through one of the principal morphisms
$1\rightarrow m$ (we leave to the reader to verify that $\Upsilon$ is equal to
this subsheaf). Suppose $f\in \Upsilon _P$ and $g= (g_1, \ldots , g_{l+k+1})\in
h^a_P$, projecting to the same element of $h(M,m,1^{k+1})_P$. Note that $g$
projects to the element $(g_1, \ldots , g_{l+k+1}, s)$ where $s: 0\rightarrow
1$ denotes the source map (or really its dual but for purposes of the present
argument we omit the dual notation).  In particular $f$ is equivalent to
$(g_1, \ldots , g_{l+k+1}, s)$, which implies that (up to changing $f$ and $g$
in their equivalence classes) $g_{l+1}$ factors through one of the
principal maps
$1\rightarrow m$ and $f_{l+k+2}=s$. This exactly means that $f$ comes from
$\Upsilon '_P \rightarrow \Upsilon _P$ and $g$ from
$\Upsilon ' _P \rightarrow h^a_P$. Thus $f$ and $g$ are equivalent in the
coequalizer, giving the second of our semiexactness statements.

The proof of the third semiexactness statement is
the same as that of the second (although $s$ above would be replaced by $t$).
This completes the proof that the standard morphisms $\Sigma \rightarrow h$ are
cofibrations (modulo the boundary cases which we have left to the reader).

An $n$-precat $A$ is an easy $n$-category if and only if every morphism
$\Sigma (M,[m], \langle k,k+1\rangle )\rightarrow A$ extends to a morphism
$h(M, m, 1^{k+1})\rightarrow A$.   A marked easy $n$-category is an $n$-precat
$A$ together with choice of extension for every morphism
$\Sigma (M,[m], \langle k,k+1\rangle )\rightarrow A$. Finally, we say that a
{\em partially marked $n$-precat} is an $n$-precategory provided with a
distinguished subset $\mu$ of the set of all morphisms of the form
$f:\Sigma (M,[m], \langle k,k+1\rangle )\rightarrow A$, and for each such
morphism, a chosen extension $f^{\mu}$ to $h(M, m, 1^{k+1})$.

If $(A, \mu )$ is a partially marked $n$-precat, then we define a new partially
marked $n$-precat $Raj(A, \mu )$ by taking the pushout via
$\varphi (M, [m], \langle k,k+1\rangle )$ for all morphisms
$$
\Sigma (M,[m], \langle k,k+1\rangle )\rightarrow A
$$
which are not in the
subset $\mu$ of marked ones.

{\em Remark:}
In the above notations, if $m=1$ then
$$
\Sigma (M, [1], \langle k,k+1\rangle ) =
\Upsilon (M, [1], 1^{k+1})  = h(M, 1 , 1^{k+1})
$$
so the pushout by $\varphi (M, [1], \langle k,k+1\rangle )$ is trivial and we
can ignore these cases if we like in the previous notation (and also in the
notion of marking).

\begin{lemma}
If $A$ is an $n$-precat then the marked easy $n$-category $Cat(A)$ is
obtained by iterating infinitely many times (i.e. over the first countable
ordinal) the operation $(A', \mu ') \mapsto Raj(A', \mu ')$, starting with $(A,
\emptyset )$.
\end{lemma}
{\em Proof:}
If $(B, \nu )$ is a marked easy $n$-category and $(A', \mu ') \rightarrow
(B, \nu
)$ is a morphism compatible with the partial marking of $A'$, then
there is a unique extension to a morphism $Raj( A', \mu ' )\rightarrow B$
compatible with the partial marking of $Raj(A', \mu ')$.
It follows that if we set $Cat '(A)$ equal to the result of the iteration
described in the lemma, then there is a unique morphism $Cat '(A)\rightarrow B$
compatible with the partial marking of $Cat '(A)$ and extending the given
morphism $A\rightarrow B$.  But $Cat '(A)$ is fully marked. By the universal
property of $Cat (A)$ this implies that $Cat (A)=Cat '(A)$.
\eop

We will often have a need for the following construction. If $A$ is an
$n$-precat then iterate (over the first countable ordinal) the operation
$(A', \mu ') \mapsto Raj (A', \emptyset )$. Call this $BigCat(A)$.  Another way
to describe this consruction is that we throw in an infinite number of times
the pushouts of all of the required diagrams (which is in some sense a more
obvious way to obtain an $n$-category).  There is an obvious morphism
$Cat (A) \rightarrow BigCat (A)$.  One of the advantages of the $BigCat$
construction is that $BigCat(A) \cong  BigCat (BigCat(A))$ (although the natural
maps are not this isomorphism).  More generally, we will use the terminology
``reordering'' below to indicate that a sequence of pushouts can be done in any
order (subject to the obvious condition that the things over which the pushouts
are being done exist at the time they are done!), which yields isomorphisms
such as  $BigCat(A) \cong  BigCat (BigCat(A))$.

If $B\leftarrow A \rightarrow C$ is  a diagram of $n$-categories, then
we define the {\em category-theoretic pushout} to be $Cat (B\cup ^AC)$.
It is again an $n$-category.  We will also often use just the pushout of
$n$-precats, i.e. the pushout of presheaves over $\Theta ^n$.

\numero{The closed model category structure}

We now come to the first main definition. A morphism $A\rightarrow B$ of
$n$-precats (that is, a morphism of presheaves on $\Theta ^n$) is a {\em weak
equivalence} if the induced morphism $Cat(A)\rightarrow Cat(B)$ is an
exterior equivalence of $n$-categories in the sense of Tamsamani
(\cite{Tamsamani} \S 1.3).  Note in particular that we don't require it to be an
easy equivalence---which would be too strong a condition.

The second main definition is relatively easy: we would like to say that a
morphism $A\rightarrow B$ of $n$-precats is a  cofibration if it is a
monomorphism of presheaves on $\Theta ^n$.  However, this doesn't work out well
at the top degree (for example, the category of sets with isomorphisms as weak
equivalences and injections as cofibrations, is not a closed model category
\cite{Quillen}).  Thus we leave the top level alone and say that a morphism
$A\rightarrow B$ of $n$-precats is a {\em cofibration} if for every $M= (m_1,
\ldots , m_k)$ with $k <n$, the morphism $A_M \rightarrow B_M$ is injective. A
cofibration which is a weak equivalence is called a {\em trivial cofibration}.

The third definition which goes along automatically with these two is that a
morphism $A\rightarrow B$ of $n$-precats is a {\em fibration} if it satisfies
the lifting property for trivial cofibrations, that is if every time
$U\hookrightarrow V$ is a trivial cofibration and $U\rightarrow A$ and
$V\rightarrow B$ are morphisms inducing the same $U\rightarrow B$ then there
exists a lifting to $V\rightarrow A$ compatible with the first two morphisms.

We recall from \cite{Quillen} the definition of {\em closed model category},
as well as from \cite{QuillenAnnals} an equivalent set of
axioms.

\begin{theorem}
\label{cmc}
The category of $n$-precats with the weak equivalences, cofibrations and
fibrations defined above, is a closed model category.
\end{theorem}

\subnumero{Some lemmas}
The proof of Theorem \ref{cmc} is by induction on $n$. Thus we may assume
that the theorem and all of the lemmas contained in the present section and \S\S
4-6 are true for $n'$-precats for all $n'<n$.  In view of this, we state all (or
most) of the lemmas before getting to the proofs.

Our proof will be modelled on the proof of Jardine that simplicial presheaves
on a site form a closed model category \cite{Jardine}.
The main lemma that we
need (which corresponds to the main lemma in \cite{Jardine}) is

\begin{lemma}
\label{pushout}
Suppose $A\rightarrow B$ is a trivial cofibration and $A\rightarrow C$
is any morphism. Let  $D = B\cup ^AC$ be the push-out of these two
morphisms (the
push-out of $n$-precats). Then the morphism $C\rightarrow D$ is a
weak equivalence.
\end{lemma}

This lemma speaks of push-out of $n$-precats. Applying the construction $Cat$
we obtain a notion of push-out of $n$-categories: if $A\rightarrow B$ and
$A\rightarrow C$ are morphisms of $n$-categories (i.e. morphisms of the
corresponding $n$-precats) then define the {\em push-out $n$-category} to be
$$
Cat ( B\cup ^AC ).
$$
\begin{corollary}
If $A\rightarrow B$ is an equivalence of $n$-categories then
$$
C\rightarrow Cat ( B\cup ^AC )
$$
is an equivalence of $n$-categories.
\end{corollary}

We will come back to push-out below in the section on Siefert-Van Kampen.

Going along with the previous lemma is something that we would like to know:

\begin{lemma}
\label{equiv}
If an $n$-precat $A$ is an $n$-category in the sense of
\cite{Tamsamani} then the morphism $A\rightarrow Cat(A)$ (resp. the morphism
$A\rightarrow BigCat(A)$) is an equivalence of $n$-categories in the sense of
\cite{Tamsamani}. \end{lemma}

Another lemma which is an important technical point in the proof of everything
is the following.  An $n$-precat $A$ can be considered as a collection $\{
A_{m/}\}$ of $n-1$-precats (functor of $\Delta$ and the first element is a
set).  We obtain the collection $\{
Cat(A_{m/})\}$ which is a functor from $\Delta$ to the category of
$n-1$-precats.  Divide by
 the equivalence relation setting the $0$-th element to a constant
$n-1$-precat, in this way we obtain a new $n$-precat denoted $Cat _{\geq 1}
(A)$.

\begin{lemma}
\label{partialCat1}
Suppose that $A$ and $B$ are $n$-precats and  $f: A\rightarrow B$ is a morphism
which induces an equivalence on the $n-1$-categories $Cat (A_{m/})\rightarrow
Cat (B_{m/})$. Then $Cat(A)\rightarrow Cat(B)$ (resp. $BigCat(A)\rightarrow
BigCat(B)$) is an equivalence.
\end{lemma}

\begin{corollary}
\label{partialCat}
The morphism
$$
Cat (A) \rightarrow Cat (Cat _{\geq 1}(A))
$$
is an equivalence of $n$-categories.
\end{corollary}
{\em  Proof:} The morphism $A \rightarrow Cat _{\geq 1}(A)$ satisfies the
hypotheses of the previous lemma so the corollary follows from the lemma.
\eop

\begin{lemma}
\label{coherence}
For any $n$-precat $A$, the morphism $A\rightarrow Cat(A)$ (resp. $A\rightarrow
BigCat(A)$) is a weak equivalence.
\end{lemma}

The closed model structure that we already have by induction for $n-1$-precats
allows us to deduce some things about $n$-categories.  Let $HC_{n-1}$
denote the
localization of $PC_{n-1}$ by inverting the set weak equivalences, which is also
(see \ref{htytype} below) the localization of $n-1$-categories by
inverting the set of
equivalences.   We know from the closed model structure \cite{Quillen} that this
is equivalent to the category of fibrant (and automatically cofibrant) objects
where we take as morphisms, the homotopy classes of morphisms.
We also know that a morphism in $PC_{n-1}$ is a weak equivalence if and only if
it projects to an isomorphism in $HC_{n-1}$ (\cite{Quillen},
Proposition 1, p. 5.5).  In particular by \ref{equiv} in degree $n-1$ we know
that a morphism of $n-1$-categories is an equivalence if and only if it projects
to an isomorphism in $HC_{n-1}$.

Suppose $A$ is an $n$-category. For $x,y\in A_0$ we have an
$n-1$-category $A_1(x,y)$ which we could denote by $Hom_A(x,y)$. Let
$LHom_A(x,y)$ denote the image of this object in the localization $HC_{n-1}$. On
the other hand, the truncation $T^{n-1}A$ is a $1$-category. We claim that for
$x$ fixed, the mapping $y\mapsto LHom_A(x,y)$ is a functor from $T^{n-1}A$ to
$HC_{n-1}$.  Similarly we claim that for $y$ fixed the mapping $x\mapsto
LHom_A(x,y)$ is a contravariant functor from $T^{n-1}A$ to $HC_{n-1}$.  These
claims give some meaning at least in a homotopic sense to the notion of
``composition with $f: y\rightarrow z$'' as a map  $LHom _A(x,y)\rightarrow
LHom_A(x,z)$.

We prove the first of the two claims, the proof for the second one being
identical. Note that these arguments are generalizations of what is mentionned
in \cite{Tamsamani} Proposition 2.2.8 and the following remark. Suppose $f\in
A_1(y,z)$. Then let $A_2(x,y,f)$ be the homotopy fiber of $A_2(x,y,z)\rightarrow
A_1(y,z)$ over the object $f$ (this is calculated by  replacing the above map by
a fibrant map and taking the fiber). The condition that $A_2$ be equivalent to
$A_1\times _{A_0}A_1$ implies that this homotopy fiber maps by an equivalence to
$A_1(x,y)$. On the other  hand it maps to $A_1(x,z)$ and this diagram gives a
morphism $LHom _A(x,y)\rightarrow  LHom _A(x,z)$ in the localized category
$H_{n-1}$.  We just have to check that this morphism is independent of the
choice of $f$ in its equivalence class.
For this we use Proposition \ref{intervalK} below (there is no circularity
because we are discussing $n-1$-categories here).  If $f$ is equivalent to
$g$ as
elements of the $n-1$-category $A_1(y,z)$
then let $K$ denote the $n-1$ category given
by \ref{intervalK}; there is a morphism $K\rightarrow A_1(y,z)$ sending $0$
to $f$ and $1$ to $g$, and since $K$ is a contractible object (weakly equivalent
to $\ast$) this proves that the homotopy fibers over $f$ or $g$ are
equivalent to
the homotopy fiber product with $K$; we have a single map from here to
$A_1(x,z)$ so our two maps induced by $f$ and $g$ are homotopic.

Associativity is given by a
similar argument using $A_3$ which we omit.

Once we have our functors $T^{n-1}A\rightarrow HC_{n-1}$ we obtain the
following type of statement: suppose $f$ is an equivalence between $u$ and $x$,
then composition with $f$ induces an equivalence $LHom_A(x,y)\cong LHom
_A(u,y)$ (and similarly for composition in the second variable).

\begin{lemma}
\label{remark}
If
$$
A\stackrel{f}{\rightarrow }B \stackrel{g}{\rightarrow } C
$$
is a pair of morphisms of $n$-categories such that any two of $f$, $g$ or
$g\circ f$ are equivalences in the sense of \cite{Tamsamani} then the third is
also an equivalence.

If $f: A\rightarrow B$ and $g:B\rightarrow A$ are two morphisms of
$n$-categories such that  $fg$ is an  equivalence and $gf$ is the identity then
$f$ and $g$ are equivalences.
\end{lemma}
{\em Proof:}
The fact that composition of equivalences is an equivalence is \cite{Tamsamani}
Lemme 1.3.5.  The statement concluding that $f$ is an equivalence if $g$ and
$g\circ f$ are, is a direct consequence of  Tamsamani's interpretation
of equivalence in terms of truncation operations (\cite{Tamsamani}
Proposition 1.3.1).

For the conclusion for $g$, note first of all that on the level of truncations
$T^nA\rightarrow T^nB \rightarrow T^nC$ the fact that
$f$ and $gf$ are isomorphisms of sets implies that $g$ is an isomorphism of
sets. This gives essential surjectivity.  Now suppose $x,y$ are objects of $B$.
Choose objects $u,v$ of $A$ and equivalences $f(u)\sim x$ and $f(v)\sim y$.
Then composition with these equivalences induces an isomorphism in
the localized category $HC_{n-1}$ between $LHom _B(x,y)$ and $LHom_B(f(u),
f(v))$ (see the discussion preceeding this lemma).  The image under $g$ of this
isomorphism is the same as composition with the images of the equivalences, so
we have a diagram
$$
\begin{array}{ccc}
LHom _B(x,y) & \rightarrow & LHom _C(g(x), g(y)) \\
\downarrow &&\downarrow \\
LHom _B(f(u), f(v)) & \rightarrow & LHom _C(gf(u), gf(v))
\end{array}
$$
in the category  $HC_{n-1}$. The horizontal arrows are the localizations of the
arrows $g:B_1(x,y)\rightarrow C_1(g(x), g(y))$ etc., and the vertical arrows
are composition with our chosen equivalences, isomorphisms in $HC_{n-1}$.
On the other hand, the bottom arrow fits into a diagram
$$
LHom _A(u,v)
LHom _B(f(u), f(v))  \rightarrow  LHom _C(gf(u), gf(v))
$$
where the first arrow induced by $f$ is an isomorphism, and the composed arrow
induced by $gf$ is an isomorphism; thus the bottom arrow of the previous
diagram is an isomorphism therefore the top arrow is an isomorphism in
$H_{n-1}$. This implies that the morphism $B_1(x,y)\rightarrow C_1(g(x), g(y))$
is an equivalence of $n-1$-categories.  This is what we needed to prove to
complete the proof that $g$ is an equivalence.

We turn now to the second paragraph of the lemma: suppose $f$ and $g$ are
morphisms of $n$-categories such that  $fg$ is an equivalence and $gf$ is the
identity. The corresponding fact for sets shows that $T^nf$ and $T^ng$ are
isomorphisms between the sets of equivalence classes of objects $T^nA$ and
$T^nB$.  Suppose $x,y\in A_0$. Note that $gf(x)=x$ and $gf(y)=y$.
We obtain morphisms
$$
A_1(x,y) \stackrel{f}{\rightarrow} B_1(fx,fy)
 \stackrel{g}{\rightarrow} A_1(gfx,gfy)=A_1(x,y)
 \stackrel{f}{\rightarrow} B_1(fgfx, fgfy) = B_1(fx,fy).
$$
We have again that $gf$ is the identity on $A_1(x,y)$ and $fg$ is an
equivalence on $B_1(fx,fy)$.  Inductively by our statement for
$n-1$-categories, the morphism $f: A_1(x,y)\rightarrow B_1(fx,fy)$ is an
equivalence. This implies that $f: A\rightarrow B$ is an equivalence and hence
that $g$ is an equivalence (by the first paragraph of the lemma).
\eop

\begin{corollary}
\label{htytype}
The localized category of $n$-precats modulo weak equivalence is equivalent to
the category $Ho-n-Cat$ of $n$-categories localized by equivalence defined in
\cite{Tamsamani}.
\end{corollary}
{\em Proof:}
The functor $Cat$ sends weak equivalences to equivalences (by Lemma \ref{equiv}
together with Lemma \ref{remark}).  Thus it induces a functor $c$ on
localizations. Let $i$ be the functor induced on localizations by the inclusion
of $n$-categories in $n$-precats. The natural transformation $A\rightarrow
Cat(A)$ gives a natural isomorphism $1 \cong c\circ i$ of functors on the
localization of $n$-categories.  On the other hand, Lemma \ref{coherence} says
that the same natural transformation induces a natural isomorphism $1\cong
i\circ c$ of functors on the localization of $n$-precats.
\eop

We will prove lemmas \ref{pushout}, \ref{equiv} and \ref{partialCat1} all at
once in one big induction on $n$. Thus we may assume that they hold for $n' <
n$.  All lemmas from here until the end of the big induction presuppose that we
know the inductive statement for $n'<n$.

{\em Remark on the passage between $Cat$ and $BigCat$ in \ref{equiv} and
\ref{partialCat1}:}
The statements for $Cat$ and $BigCat$ are equivalent.  Take Lemma \ref{equiv}
for example.  If $A\rightarrow Cat(A)$ is an equivalence for any $n$-category
$A$ then $BigCat(A)$ can be constructed as the iteration over the first
countable ordinal of the operation $A' \mapsto Cat (A')$ (and starting at $A$).
The morphisms at each stage in the iteration are equivalences, so it follows
that the morphism $A\rightarrow BigCat(A)$ is an equivalence. On the other
hand, suppose we know that $A\rightarrow BigCat(A)$ is an equivalence for any
$n$-category $A$.  Then $Cat(A)\rightarrow BigCat(Cat(A))$ is an equivalence,
but by reordering $BigCat(Cat(A))=BigCat(A)$. Thus the hypothesis also gives
that $A\rightarrow BigCat(Cat(A))$ is an equivalence. Lemma \ref{remark} then
implies that $A\rightarrow Cat(A)$ is an equivalence. We obtain the required
statement concerning \ref{partialCat1} by using the fact that
$Cat(A)\rightarrow BigCat(A)$ (resp. $Cat(B)\rightarrow BigCat(B)$) is an
equivalence---note that our proof of \ref{partialCat1} comes after our proof
of \ref{equiv} below---and applying \ref{remark}.

\subnumero{A simplified point of view}
We started to see, in the proof of Lemma \ref{remark}, a simplified or
``derived'' point of view on $n$-categories. We will expand on that a bit more
here. When we use the statements of the above lemmas for $n-1$-categories, they
may be considered as proved in view of our global induction. The
homotopy or localized category $HC_{n-1}$ of $n-1$-precats modulo weak
equivalence, also equal to the localization of Tamsamani's $(n-1)-Cat$ by
equivalences, admits direct products.  There is a functor $T^{n-1}: HC_{n-1}
\rightarrow Sets$ related to the inclusion $i: Sets \subset HC_{n-1}$ by
morphisms $iT^{n-1}(X)\rightarrow X$ and $T^{n-1}iS \cong S$ (the first is only
well defined in the localized category). Thus if $A\times B\rightarrow C$ is a
morphism in $HC_{n-1}$ then we obtain the map $A\times iT^nB \rightarrow C$.

Note that $HC_{n-1}$ admits fibered products over objects of the form $i(S)$
for $S$ a set, since these are essentially just direct products. (However the
homotopy category does not admit general fibered products nor, dually, does it
admit pushouts.)

We can define the notion of $HC_{n-1}$-category, as
simply being a category in the category $HC_{n-1}$ such that the object object
is a set.  Applying the functor $T^{n-1}$ yields a category, and this category
acts on the morphism objects of the previous one, using the above remark.

If $A$ is an $n$-category then taking $A_0$ as set of objects and using the
object $LHom_A(x,y)$ as morphism object in $HC_{n-1}$ we obtain an
$HC_{n-1}$-category which we denote $HC_{n-1}(A)$. We can write
$$
Hom _{HC_{n-1}(A)}(x,y) := LHom_A(x,y)
$$
which in turn is, we recall, the image of $A_1(x,y)$ in the localization of
the category of $n-1$-precats.  This is what we used in the proof of Lemma
\ref{remark} above. The truncation operation $T^{n-1}$ applied to $HC_{n-1}(A)$
gives the $1$-category $T^{n-1}A$. We obtain again the action of this category
on the  morphism objects in $HC_{n-1}(A)$.

In the next section we will be interested in the notion of
{\em $HC_{n-1}$-precategory}, a functor $F:\Delta \rightarrow HC_{n-1}$
sending $0$ to a set.   An $HC_{n-1}$-precategory $F$ is an $HC_{n-1}$-category
if and only if the usual morphisms
$$
F_p\rightarrow F_1\times _{F_0} \ldots \times _{F_0}F_1
$$
are isomorphisms.
If $A$ is an $n$-precat then let $HC_{n-1}(A)$
denote the $HC_{n-1}$-precategory which to $p\in \Delta$ associates the image of
$A_{p/}$ in the localized category $HC_{n-1}$.

Here is a small remark which is sometimes useful.

\begin{lemma}
\label{HCequivCat}
Suppose $f:A\rightarrow B$ is a morphism of $n$-categories
and suppose that $HC_{n-1}(A)\rightarrow HC_{n-1}(B)$ is an equivalence in
$HC_{n-1}Cat$.  Then $f$ is an
equivalence.
\end{lemma}
\eop

We can also make a similar statement for $HC_{n-1}$-precats under the condition
of requiring an isomorphism on the set of objects.
\begin{lemma}
\label{HCequivPreCat}
Suppose $f:A\rightarrow B$ is a morphism of $n$-precats
and suppose $HC_{n-1}(A)\rightarrow HC_{n-1}(B)$ is an isomorphism of functors
$\Delta \rightarrow HC_{n-1}Cat$.  Then $f$ is a weak
equivalence.
\end{lemma}
{\em  Proof:}
This is just  a restatement of Lemma \ref{partialCat1} (in particular it is not
available for use in degree $n$ until we have proved \ref{partialCat1} below).
\eop

It would have been nice to be able to
have an operation on $HC_{n-1}$-precategories which, when applied to
$HC_{n-1}(A)$ yields $HC_{n-1}(Cat(A))$.  This doesn't seem to be possible
(although I don't have a counterexample) because the construction we discuss in
the next section relies heavily on pushouts but these don't exist in $HC_{n-1}$.
If this had been possible we would have been able to formulate a notion of weak
equivalence for $HC_{n-1}$-precats and in particular we would have been able to
give a stronger formulation in the previous lemma.

We end this discussion by pointing out that information is lost in passing from
$A$ to $HC_{n-1}(A)$. (See the next paragraph for some counterexamples but I
don't have counterexamples for all of the nonexistence statements which are
made.)
 Let $HC_{n-1}Cat$ (resp. $HC_{n-1}PreCat$) denote the
categories of $HC_{n-1}$-categories (resp. $HC_{n-1}$-precategories).
The functors $n-Cat\rightarrow HC_{n-1}Cat$ and $n-Cat \rightarrow HC_n$ do
not enter into a commutative triangle with a morphism between  $HC_n$ and
$HC_{n-1}Cat$ in either direction. The only thing we can say is that there is
an obvious notion of equivalence between two $HC_{n-1}$-categories, and if
we let $Ho-HC_{n-1}-Cat$ denote the category of $HC_{n-1}$-categories
localized by inverting these equivalences, then there is a factorization
$$
n-CAT \rightarrow HC_n \rightarrow Ho-HC_{n-1}-Cat
$$
but the second arrow in the factorization is not an isomorphism.  In
particular, when we pass from $A$ to $HC_{n-1}(A)$ we lose information.
Nonetheless, it may be helpful especially from an intuitive point of view to
think of an $n$-category in terms of its associated object $HC_{n-1}(A)$ which
is a category in the homotopy category of $n-1$-categories.

The topological analogy of the above situation (which can be made precise using
the Poincar\'e groupoid and realization constructions \cite{Tamsamani}---thus
providing some counterexamples to support some of the the nonexistence
statements
made in the previous paragraph) is the following: if $X$ is a space then
for each
$x,y\in X$ we can take  as $h(x,y)$ the space of paths from $x$ to $y$ viewed as
an object in the homotopy category $Ho(Top)$.  We obtain a category in
$Ho(Top)$.  If $X$ is connected it is a groupoid with one isomorphism class,
thus essentially a group in $Ho(Top)$.  This group is just the loop space based
at any choice of point,  viewed as a group in $Ho(Top)$. It is well known
(\cite{Adams} \cite{Tanre}) that this object does not suffice to reconstitute
the homotopy type of $X$,  thus our functor from $Top$ to the category of
groupoids in $Ho(Top)$ does not yield a factorization  of the localization
functor $Top\rightarrow Ho(Top)$. On the other hand, since there is no way to
canonically choose a collection of basepoints for an object in $Ho(Top)$, there
probably is not a factorization in the other direction either.

\subnumero{Another simplified point of view}
  We now give another set of remarks relating the present approach to
$n$-categories with the usual standard ideas.  This is based on the following
observation. The proof of the lemma is based on some ideas from the next
section so the reader should look there before trying to follow the proof. We
have put the lemma here for expository reasons.

\begin{lemma}
\label{fibrantpieces}
If $A$ is a fibrant $n$-precat then the $A_{p/}$ are fibrant $n-1$-precats.
\end{lemma}
{\em Proof:}
Fix objects $x_0, \ldots , x_p\in A_0$.  We show that $A_{p/}(x_0, \ldots ,
x_p)$
is fibrant. Suppose $U\rightarrow V$ is a trivial cofibration of $n-1$-precats.
Let $B$ (resp. $C$) be the $n$-precat with objects $0, \ldots , p$
and such that $B_{q/}(i_0, \ldots , i_q)$ (resp. $C_{q/}(i_0, \ldots , i_q)$) is
the disjoint union of $U$ (resp. $V$)
over all morphisms $f:q\rightarrow p$ such that $f(q)=i_q$.  Then (as can be
seen by the discussion of the next section) $B\rightarrow C$ is a trivial
cofibration.  A morphism $B\rightarrow A$ (resp. $C\rightarrow A$) is the same
thing as a morphism $U\rightarrow A_{p/}(x_0,\ldots , x_p)$ (resp. $V\rightarrow
A_{p/} (x_0,\ldots , x_p)$).
It follows immediately that if $A$ is fibrant then $A_{p/}(x_0,\ldots , x_p)$
has the required
lifting property to be fibrant. \eop

Now we can use the closed model category structure on $PC_{n-1}$ to analyze the
collection of $A_{p/}$ when $A$ is fibrant.  Recall that morphisms in the
localized category between fibrant and cofibrant objects are represented by
actual morphisms \cite{Quillen}.  Thus the morphism
$$
A_{2/} \rightarrow A_{1/} \times _{A_0} A_{1/}
$$
which is an equivalence, can be inverted and then followed by the projection
to the third edge of the triangle to give
$$
A_{1/} \times _{A_0} A_{1/}\rightarrow A_{2/} \rightarrow A_{1/}.
$$
We get a morphism ``composition''
$$
m:A_1(x,y)\times A_1(y,z) \rightarrow A_1(x,z)
$$
which  represents the composition
$$
LHom_A(x,y)\times LHom_A(y,z)\rightarrow
LHom_A(x,z)
$$
of the previous ``simplified point of view''.
Of course our composition morphism $m$ is not uniquely determined but depends
on the choice of inversion of the original equivalence.  In particular $m$ will
not in general be associative. However $A_{3/}$ gives a homotopy in the
sense of Quillen between  $m(m(f,g),h)$ and $m(f, m(g,h))$.  This can be turned
into a homotopy in the sense of the $n-1$-categories of morphisms (an exercise
left to the reader).

\numero{Calculus of ``generators and relations''}
For the proofs of \ref{equiv} and \ref{partialCat1} we need a close analysis of
an operation which when iterated yields $BigCat$.
This analysis will lead us to a point of view which generalizes the idea of
generators and relations for an associative monoid. At the end we draw as a
consequence one of the main special lemmas needed to treat the special case
\ref{specialcase} in the proof of \ref{pushout}.

The overall goal of this section is to investigate the operation $A\mapsto
Cat(A)$ in the spirit of looking at the simplicial collection of
$n-1$-precats $A_{p/}$ as a functor from $\Delta$ to our closed model
category in degree $n-1$. We would like to understand the transformation which
this functor undergoes when we apply the operation $Cat$ to $A$.

We first describe a general type of operation which we often encounter. Suppose
$A$ is an $n$-precat and suppose $A_{m/}\rightarrow B$ is a cofibration of
$n-1$-precats provided with a morphism $\pi :B\rightarrow A_0 \times \ldots
\times A_0$ making the composition
$$
A_{m/}\rightarrow
B\rightarrow A_0 \times \ldots \times
A_0
$$
equal to the usual morphism (there are $m+1$ factors $A_0$ in the product).  We
can alternatively think of this as a collection of cofibrations
$$
A_{m/} (x_0, \ldots , x_m) \rightarrow B(x_0,\ldots , x_m)
$$
for all sequences of objects $x_i \in A_0$.  Then we define
the cofibration of $n$-precats
$$
A\rightarrow \Ii (A; A_{m/}\rightarrow B)
$$
as follows (the projection $\pi$ is  part of the data even though
it is not contained in the notation). For any $p$, $\Ii (A; A_{m/}\rightarrow
B)_{p/}$ is the multiple pushout of $A_{p/}$ and $A_{m/}\rightarrow B$ over all
morphisms  $A_{m/} \rightarrow A_{p/}$ coming from morphisms $p\rightarrow m$
which do not factor through $0$.  Functoriality is defined as follows: if
$q\rightarrow p$ is a morphism then for any $f:p\rightarrow m$ such that the
composition $q\rightarrow m$ doesn't factor through $0$, we define the morphism
of functoriality on the part of the pushout corresponding to $f$ as the identity
in the obvious way; on the other hand, if $f: p\rightarrow m$ is a morphism such
that the composition $q\rightarrow m$ factors through $0\rightarrow m$ then we
obtain (from the projection $\pi$) a morphism $B\rightarrow A_0$ extending the
morphism $A_{m/}\rightarrow A_0$ and so that part of the pushout is sent into
the image of $A_0$ in $A_{q/}$.

We call $A\rightarrow \Ii (A; A_{m/}\rightarrow B)$ the {\em pushout of $A$
induced by $A_{m/}\rightarrow B$}.
Using Lemma \ref{pushout}  in degree $n-1$ we find that if $A_{m/}\rightarrow
B$ is a trivial cofibration then the morphisms
$$
A_{p/} \rightarrow \Ii (A; A_{m/}\rightarrow B)_{p/}
$$
are trivial cofibrations.

This operation occurs notably in the process of doing $Cat$ or $BigCat$ to $A$.
Fix $m\geq 1$, $M=(m_1, \ldots , m_l)$, $m'$ and $k$.  Let
$$
\Sigma := \Sigma (m, M, [m'], \langle
k,k+1 \rangle ),
$$
and
$$
\varphi := \varphi (m, M, [m'], \langle k,k+1 \rangle ):\Sigma
\rightarrow h(m,M,m', 1^{k+1}).
$$
Suppose again $a: \Sigma\rightarrow A$ is a morphism and let $C$ be the pushout
$n$-precat
of $A$ and $\varphi$ over $a$. In this case note that $\Sigma$ and $h(m,M,m',
1^{k+1})$ are pushouts of diagrams of objects entirely within the category
$(m, \Theta ^{n-1})$ of objects of the form $(m,M')$.  The restriction of $A$
to this category is just $A_{m/}$. Let $\psi : A_{m/}\rightarrow F$ be the
pushout   $n-1$-precat of  $\varphi$ over $a$ considered in this way.
Then $C = \Ii (A; A_{m/}\rightarrow F)$ is the pushout of $A$ induced by
$\psi$ (note that $\psi$ admits a projection $\pi$ in an obvious way).
The proof
is that $h(m,M,m', 1^{k+1})$ has exactly the same description as a pushout of
$\Sigma$.

Suppose $A$ is an $n$-precat.
Define a new $n$-precat $Fix(A)$ by iterating the above operation of pushout
by all standard cofibrations $\varphi (m, M, [m'], \langle k,k+1 \rangle )$,
over all possible values of $m$, $M$, $m'$ and $k$, and repeating this
operation a countable number of times. By reordering, $Fix(A)$ may be seen as
obtained from $A$ by a sequence of standard pushouts of the form
$$
A' \rightarrow \Ii (A', A'_{m/} \rightarrow BigCat(A'_{m/})).
$$
In particular it is clear that each $A_{p/}\rightarrow Fix(A)_{p/}$ is a trivial
cofibration. On the other hand it is also clear that the $Fix(A)_{p/}$ are
$n-1$-categories (they are obtained by iterating operations of the form,
taking $BigCat$ then taking a bunch of pushouts
then taking $BigCat$ and so on an infinite number of times---and such an
iteration is automatically an easy $n-1$-category).

In order to get to $Cat(A)$ or $BigCat(A)$ we need another type of operation
which relates the different $A_{m/}$.
Suppose $A$ is an $n$-precat, fix $m\geq 2$ and suppose  that we have a diagram
$$
A_{m/}
\stackrel{f}{\rightarrow} Q \stackrel{g}{\rightarrow} A_{1/} \times _{A_0}
\ldots \times _{A_0}A_{1/}
$$
with the first arrow cofibrant. Then we define the pushout
$A\rightarrow \Jj (A; f,g)$ as follows.
$\Jj (A; f,g)_{p/}$ is the multiple pushout of $A_{m/}\rightarrow A_{p/}$
and $A_{m/} \rightarrow Q$ over all maps $p\rightarrow m$ not factoring through
any of the principal maps $1\rightarrow m$.  The morphisms of functoriality are
defined in the same way as for the construction $\Ii$ using the map $g$ here.

{\em Remark:}  This pushout changes the object over $1\in \Delta$ because
there are morphisms $1\rightarrow m$ (the faces other than the principal ones)
which don't factor through the principal face maps.

The remaining of our standard pushouts which are not covered by the operation
$\Ii$ are covered by this operation $\Jj$.
Fix some $m\geq 2$ and $k$.  Write
$$
\Sigma\;\;\; \mbox{for} \;\;\; \Sigma ([m], \langle
k,k+1 \rangle ),
$$
and
$$
\varphi := \varphi ([m], \langle k,k+1 \rangle ):\Sigma
\rightarrow h(m, 1^{k+1}).
$$
We have a diagram of $n-1$-precats
$$
\Sigma _{m/}
\stackrel{f(m,k)}{\rightarrow} h(1^{k+1}) \stackrel{g(m,k)}{\rightarrow}
\Sigma _{1/} \times _{\Sigma
_0} \ldots \times _{\Sigma _0} \Sigma _{1/},
$$
and via this diagram
$$
h(m, 1^{k+1}) = \Jj (\Sigma ; f(m,k),g(m,k)).
$$
It follows that if $A$ is an $n$-precat and $\Sigma \rightarrow A$ is a
morphism then the standard pushout $B$ of $A$ along $\varphi$ is of the form
$B = \Jj (A; f,g)$ for appropriate maps $f$ and $g$  induced by the above ones.

We need to have some information about decomposing and commuting the operations
$\Ii$ and $\Jj$. Suppose
$$
A_{m/} \stackrel{\varphi}{\rightarrow}
P \stackrel{\pi}{\rightarrow} A_0 \times \ldots \times
A_0
$$
is a morphism. Let
$$
\eta :A_{1/}\rightarrow B
$$
denote the multiple pushout of $A_{1/}$ by $\varphi$ over all of the principal
morphisms $1\rightarrow m$ (with projection $\nu : B\rightarrow A_0\times
A_0$).  We obtain a factorization
$$
A_{m/} \stackrel{\varphi}{\rightarrow}
P \stackrel{\psi}{\rightarrow} B\times _{A_0} \ldots \times _{A_0} B
$$
and we have
$$
\Ii (A; \varphi , \pi )= \Jj (\Ii (A; \eta ,\nu ); \varphi , \psi ).
$$
In this way we turn an operation of the form $\Ii$ for $m$ into an operation of
the form $\Ii$ for $1$ followed by an operation of the form $\Jj$ for $m$.

We define a type of operation combining operations of the form $\Jj$ for $m$
with operations of the form $\Ii$ for $1$.
However, we would like to keep track of certain sub-$n-1$-precats of $A_{m/}$
and $A_{1/}$. So we say that an {\em $(m,1)$-painted $n$-precat} (or just {\em
painted $n$-precat} in the current context where $m$ is fixed) is an $n$-precat
$A$ together with cofibrations of $n-1$-precats  $A^{\ast}_{m/} \rightarrow
A_{m/}$ and $A^{\ast}_{1/} \rightarrow A_{1/}$.  We require a lifting of the
standard morphism to
$$
A^{\ast}_{m/} \rightarrow A^{\ast}_{1/} \times _{A_0} \ldots \times _{A_0}
A^{\ast}_{1/}.
$$

Suppose $(A,A^{\ast}_{m/},A^{\ast}_{1/})$ is a painted $n$-precat, and
suppose that we have morphisms
$$
A^{\ast}_{1/}\stackrel{\eta}{\rightarrow} B\stackrel{\nu}{\rightarrow}A_0\times
A_0
$$
and
$$
A^{\ast}_{m/} \stackrel{\varphi}{\rightarrow}
P \stackrel{\psi}{\rightarrow} B\times _{A_0} \ldots \times _{A_0} B
$$
compatible with the previous lifting of the standard morphism to the painted
parts. Let $\eta '$, $\nu '$, $\varphi '$ and $\psi '$ be obtained by taking the
pushouts of the above with $A_{1/}$ or $A_{m/}$. Then we define a new  painted
$n$-precat
$$
\Jj '(A; \eta , \nu ; \varphi , \psi )
:= \Jj (\Ii (A; \eta ',\nu '); \varphi ', \psi '),
$$
with painted parts $(P, B)$ replacing $(A_{m/}^{\ast}, A_{1/}^{\ast})$. This
operation now behaves well under iteration: the composition of two such
operations is again an operation of the same form.  Furthermore, our operations
$\Ii$ and $\Jj$ coming from standard trivial cofibrations can be interpreted as
operations of the above type if the original $\Sigma \rightarrow A$ sends the
arrows $(a,b, v_i)$ into the painted parts $A_{m/}^{\ast}, A_{1/}^{\ast}$.
These operations are exactly designed to do two things: replacing the painted
parts by their associated $n$-categories; and getting the standard map to being
an equivalence. In particular,
starting with
$A^{\ast}_{1/} = A_{1/} $ and $A^{\ast}_{m/} = A_{m/}$,
there is a sequence of
operations coming from standard trivial cofibrations (concerning only $m$ and
$1$) such that, when interpreted as operations on painted $n$-precats, combine
into one big operation of the form $\Jj'$ where $\eta : A^{\ast}_{1/}
\rightarrow B$ is a trivial cofibration to an $n-1$-category, and where the
morphism
$$
P \stackrel{\psi}{\rightarrow} B\times _{A_0} \ldots \times _{A_0} B
$$
is an equivalence of $n$-categories.

Going back to the original definition of the operation $\Jj '$ in terms of
$\Jj$ and $\Ii$ we find that an appropriate sequence of trivial cofibrations
can be reordered into an operation  of the form $A\mapsto A' = \Ii (A, \eta ,
\nu )$ followed by  $\Jj (A'; f,g)$ for
$$
A'_{m/}\stackrel{f}{\leftarrow} \Gg [m](A)
\stackrel{g}{\rightarrow} A'_{1/}\times _{A_0} \times \ldots \times _{A_0}
A'_{1/}
$$
where $g$ is a weak equivalence.
Recall that the morphism $A\rightarrow A'$ coming before the operation $\Jj$ has
the property that the $A_{p/}\rightarrow A'_{p/}$ are weak equivalences. Note
also that we can assume that $A'_{1/}$ is an $n-1$-category, because it is equal
to $B$---there are no morphisms $1\rightarrow 1$ other than the identity and
those which factor through $0$---and $B$ can be chosen to be an $n-1$-category.
(This paragraph is the conclusion we want; the discussion of painted $n$-precats
was just a means to arrive here and will not be used any further below.)

Let $Gen [m](A)$ denote the result of the previous operation, which we can thus
write as
$$
Gen [m](A)= \Jj (A'; A'_{m/}\stackrel{f}{\rightarrow} \Gg [m](A)
\stackrel{g}{\rightarrow} A'_{1/}\times _{A_0} \times \ldots \times _{A_0}
A'_{1/}).
$$
The pushouts chosen as above may be assumed to contain, in particular,
all of the standard pushouts of the second type for $m$.

Put
$Gen _1(A)=Fix(A)$ and $Gen _i (A):= Fix(Gen[i](Gen _{i-1}(A)))$ for
$i\geq 2$. Let $Gen(A)$ be the inductive limit of the $Gen _i(A)$.  Finally,
iterate the operation $A'\mapsto Gen(A')$ a countable number of times.
It is clear that, by reordering, this yields $BigCat(A)$, since on the one
hand all of the necessary pushouts occur, whereas on the other hand only the
standard pushouts are used.

\subnumero{Proofs of \ref{equiv} and \ref{partialCat1}}

The above description yields immediately the proofs of these two lemmas.

{\em Proof of  \ref{equiv}:}

Suppose
$A$ is an $n$-precat such that $A_{1/}$ is an $n-1$-category and such that
$(\ast )$ for all $m$ the morphisms
$$
A_{m/} \rightarrow A_{1/}\times _{A_0} \ldots \times _{A_0} A_{1/}
$$
are weak equivalences of $n-1$-precats.  Fix $m$ and  apply the operation
$Gen [m](A)$. Let $A'$ be the intermediate result of doing the preliminary
operations $\Ii$. The morphism
$$
f:A'_{m/} \rightarrow \Gg [m](A)
$$
is a trivial cofibration of $n$-precats,
using:
\newline
(1)\,\,  the fact that $\Gg [m](A)\rightarrow
A'_{1/}\times _{A_0} \ldots \times _{A_0} A'_{1/}$  is a weak equivalence;
\newline
(2)\,\, the fact that $A'_{1/}$ are $n$-categories equivalent to $A_{1/}$,
and noting that direct products of $n$-categories (or fibered products over
sets) preserve equivalences; and
\newline
(3)\,\, the hypothesis that $A_{m/}$ is weakly equivalent to the product of the
$A_{1/}$, again coupled with the fact that $A'_{m/}$ is weakly equivalent to
$A_{m/}$ because the operations $\Ii$ preserve the weak equivalence type of the
$A_{p/}$.

It now follows from the definition of the operation $\Jj (A' ; f,g)$ that
the morphisms
$$
A_{p/} \rightarrow A'_{p/} \rightarrow
\Jj (A' ; f,g)_{p/}
$$
are weak equivalences. Thus (under the hypothesis $(\ast )$ above) the
morphism $A\rightarrow Gen [m](A)$ induces weak equivalences
$A_{p/}\rightarrow  Gen [m](A)_{p/}$.  The same holds always for the operation
$Fix$, and iterating these we obtain the conclusion (under hypothesis $(\ast
)$) that
$$
A_{p/}\rightarrow BigCat(A)_{p/}
$$
are weak equivalences.

In the hypotheses of \ref{equiv}, $A$ is an $n$-category, so $A_{1/}$ is an
$n-1$-category and hypothesis $(\ast )$ is satisfied.  It follows from above
that the morphisms
$$
A_{p/} \rightarrow BigCat(A)_{p/}
$$
are weak equivalences of $n-1$-precats, but since both sides are
$n-1$-categories this implies that they are equivalences of $n-1$-categories
(using Lemma \ref{equiv} in degree $n-1$).  Therefore $A\rightarrow BigCat(A)$
is an equivalence of $n$-categories,
completing the proof of \ref{equiv}.
\eop

{\em Proof of \ref{partialCat1}:}

Suppose $A\rightarrow B$ is a morphism of $n$-precats which induces weak
equivalences of $n-1$-precats $A_{m/} \rightarrow B_{m/}$ for all $m$.
Replacing $A$ by $Fix(A)$ and $B$ by $Fix(B)$ conserves the hypothesis. We show
that replacement of $A$ by $Gen [m](A)$ and $B$ by $Gen[m](B)$ conserves the
hypothesis. Let $A'$ (resp. $B'$) denote the intermediate $n$-precats used in
the definition of $Gen[m](A)$ (resp. $Gen[m](B)$). These are the results of
applying operations $\Ii$ to $A$ and $B$, and we may assume as in the previous
proof that $A'_{1/}$ and $B'_{1/}$ are $n-1$-categories. We do these
operations in a canonical way so as to preserve morphisms $A'\rightarrow
B'$ and $Gen[m](A)\rightarrow Gen[m](B)$. Note that  $A_{1/}\rightarrow A'_{1/}$
is a weak equivalence and the same for $B$. Our hypothesis now implies that the
morphism $A'_{1/}\rightarrow B'_{1/}$ is an equivalence of $n-1$-categories.
Therefore
$$
A'_{1/} \times _{A'_0} \ldots \times
_{A'_0}A'_{1/} \rightarrow
B'_{1/} \times _{B'_0} \ldots \times _{B'_0}B'_{1/}
$$
is a weak equivalence.

As before
$$
Gen [m](A)= \Jj (A'; A'_{m/}\stackrel{f_A}{\rightarrow} \Gg [m](A)
\stackrel{g_A}{\rightarrow} A'_{1/}\times _{A_0} \times \ldots \times _{A_0}
A'_{1/}),
$$
with $g_A$ being a weak equivalence.  Similarly
$$
Gen [m](B)= \Jj (B'; B'_{m/}\stackrel{f_B}{\rightarrow} \Gg [m](B)
\stackrel{g_B}{\rightarrow} B'_{1/}\times _{B_0} \times \ldots \times _{B_0}
B'_{1/}),
$$
with $g_B$ a weak equivalence.  It follows immediately that
$$
\Gg [m] (A) \rightarrow \Gg [m](B)
$$
is a weak equivalence and hence (using the descripition of $\Jj$ as well as
the fact that weak equivalences on the components induce weak equivalences of
pushouts, which we know by induction for $n-1$-precats) the morphism
$Gen [m](A)\rightarrow Gen [m](B)$ induces weak equivalences
$$
Gen [m](A)_{p/}\rightarrow Gen [m](B)_{p/}.
$$
This shows that the operation $Gen[m]$ preserves the hypothesis of
\ref{partialCat1}. Taking limits we get that $Gen$ and finally that $BigCat$
preserve the hypothesis:
we get that for all $p$,
$$
BigCat (A)_{p/} \rightarrow BigCat (B)_{p/}
$$
is a weak equivalence. This implies that $BigCat(A)\rightarrow BigCat(B)$ is an
equivalence of $n$-categories, finishing the proof of \ref{partialCat1}.
\eop

\subnumero{$1$-free ordered precats}
Suppose $A$ is an $n$-precat. We say that $A$ is {\em $1$-free ordered}
if there is a total order on the set $A_0$ of objects (which we suppose for
simplicity to be finite) such that the following properties are satisfied:
\newline
(FO1)---for any sequence $x_0,\ldots , x_n$ which is out of order (i.e. some
$x_i$ is strictly bigger than $x_{i+1}$), $A_{m/}(x_0,\ldots , x_m)=\emptyset$;
\newline
(FO2)---for any  sequence $x_0, \ldots , x_n$ with $x_{i-1}\leq x_i$ the
morphism
$$
A_{m/}(x_0,\ldots , x_m)\rightarrow A_{1/}(x_0,x_m)
$$
is a weak equivalence; and
\newline
(FO3)---for any stationary sequence $A_{m/}(x,\ldots , x)$ is weakly equivalent
to $\ast$.

Properties (FO1) and (FO2) properties are preserved under the operation
$BigCat$.
Indeed, the standard cofibrations $\Sigma \rightarrow h$ go between $n$-precats
which satisfy these conditions, and these conditions are preserved by pushouts
over diagrams of
morphisms which are order-respecting (i.e. morphisms respecting $\leq$) between
$1$-free ordered $n$-precats. Note that if $A$ satisfies (FO1) then any morphism
$\Sigma \rightarrow A$ must respect the order.

The condition (FO3) is preserved by pushouts of morphisms which are strictly
order-preserving, and also by the operation $Fix$ (which is essentially
the same as $Cat _{\leq 1}$). If $A$ is $1$-free ordered, using (FO1) and (FO3)
we get that in order to obtain $Cat(A)$ it suffices to use the operation $Fix$
and pushouts for trivial cofibrations $\Sigma \rightarrow h(m, 1^{k+1})$ via
morphisms $\Sigma \rightarrow A$ which are strictly order preserving.  Thus
$BigCat(A)$ again satisfies (FO3). Furthermore, note that these trivial
cofibrations do not change the homotopy type for adjacent objects, so  if
$x,y$ are adjacent in the ordering then
$$
A_{1/}(x,y)\rightarrow BigCat(A)_{1/}(x,y)
$$
is a weak equivalence.  We obtain the following conclusion:

\begin{lemma}
\label{freeness}
Suppose $A$ is a $1$-free ordered $n$-precat with finite object set.  For two
objects $x,y\in A_0$ let $x=x_0, \ldots , x_m = y$ be the maximal strictly
increasing ordered sequence  going from $x$ to $y$. Let $A':= BicCat(A)$. Then
the morphisms  $$
A'_{1/}(x,y)\leftarrow A'_{m/}(x_0,\ldots , x_m) \rightarrow
A_{1/}(x_0,x_1)\times _{A_0} \ldots \times _{A_0} A_{1/}(x_{m-1}, x_m)
$$
are weak equivalences.  In particular if $A\rightarrow B$ is a morphism of
$1$-free ordered $n$-precats (with finite object sets) preserving the ordering,
inducing an isomorphism on object sets and inducing equivalences $A_{1/}
(x,y)\rightarrow B_{1/}(x,y)$ for all pairs of adjacent objects $(x,y)$ then
$A\rightarrow B$ is a weak equivalence.
\end{lemma}
{\em  Proof:}
In the diagram
$$
A'_{1/}(x,y)\leftarrow A'_{m/}(x_0,\ldots , x_m) \rightarrow
A'_{1/}(x_0,x_1)\times _{A_0} \ldots \times _{A_0} A'_{1/}(x_{m-1}, x_m)
$$
the left arrow is a weak equivalence by property (FO2) for $A'$ (we have shown
that that property is preserved under passage from $A$ to $A'$); and the right
arrow is a weak equivalence because $A'$ is an $n$-category.  The product on
the right may be replaced by that appearing in the statement of the lemma since
$A_{1/}(x_{i-1},x_{i})\rightarrow A'_{1/}(x_{i-1},x_{i})$ is an equivalence
because $x_{i-1}$ and $x_i$ are adjacent. This gives the first statement of the
lemma. For the second statement, note that (using the same notation for objects
of $A$ and $B$, and using the notations $A'$ and $B'$ for associated
$n$-categories) the first statement implies that $A'_{1/} (x,y)\rightarrow
B'_{1/}(x,y)$ is a weak equivalence for any pair of objects $x,y$.  This
implies that $A'\rightarrow B'$ is an equivalence of $n$-categories.
\eop

\subnumero{Characterization of weak equivalence}
We close this section by mentioning a proposition which gives a sort of
uniqueness for the notion of weak equivalence.

\begin{proposition}
\label{general}
Suppose $F: PC_n \rightarrow PC_n$ is a functor with natural transformation
$i_A: A\rightarrow F(A)$ such
that: \newline
(a)---for all $A$, $F(A)$ is an $n$-category;
\newline
(b)---if $A$ is an $n$-category then $i_A$ is an iso-equivalence of
$n$-categories (recall that this means an equivalence inducing an isomorphism
on sets of objects); and
\newline
(c)---for any $n$-precat $A$ the morphism $F(i_A): F(A) \rightarrow F(F(A))$
is an equivalence of $n$-categories.
\newline
Then for any $n$-precat $A$ the morphism $A\rightarrow F(A)$ is a weak
equivalence.
\end{proposition}
{\em Proof:}
Put $F'(A):= Cat (F(A))$. It is a marked easy $n$-category.  We have a morphism
$k_A := Cat (i_A):Cat(A)\rightarrow F'(A)$. Letting $j_A: A\rightarrow Cat(A)$
denote the inclusion and $i'_A$ the map $A\rightarrow F'(A)$, note that
$k_Aj_A = i'_A$.

The functor  $F'$  again satisfies the properties (a), (b) and (c) above. For
property (c) note that the map $F'(i'_A)$ is  obtained by applying $Cat$ to
the composed map
$$
F(A) \stackrel{F(i_A)}{\rightarrow} F(F(A))\stackrel{F(j_{F(A)})}{\rightarrow}
F(Cat(F(A))),
$$
but the first map is an equivalence by hypothesis and the second is an
equivalence because of the diagram
$$
\begin{array}{ccc}
F(A) & \rightarrow & F(F(A)) \\
\downarrow && \downarrow \\
Cat (F(A)) & \rightarrow & F(Cat (F(A)))
\end{array}
$$
where the top arrow is $i_{F(A)}$ which is an equivalence by (b) (but it
is different from $F(i_A)$!), the left vertical arrow is the equivalence
$i_{F(A)}$ (by (a)) and the bottom arrow is the equivalence $i_{Cat(F(A))}$
again an equivalence by (b).

We have the
following diagram:
$$
\begin{array}{ccccc}
Cat(A) & \rightarrow & Cat(Cat(A)) &&\\
\downarrow & & \downarrow &&\\
F' (A) & \rightarrow & F'(Cat(A)) & \rightarrow & F'(F'(A))
\end{array} .
$$
The morphism on the top is $Cat(j_A)$, and the vertical
morphisms are $k_A$ and $k_{Cat(A)}$ respectively. The
morphisms on the bottom are $F'(j_A)$ and $F'(k_A)$ respectively. The diagram
comes from naturality of $k$.
The top arrow
$Cat(j_A)$ and  the middle vertical  arrow $k_{Cat(A)}$ are equivalences, as is
the composition along the bottom
$F'(i'_A)$.  On the other hand, all of the arrows are identities on the sets of
objects.  Thus, when morphisms are equivalences they are in fact
iso-equivalences, and in particular equivalences on the level of the
$n-1$-categories $(\cdot )_{p/}$. The closed model structure for $n-1$ implies
that a morphism of $n-1$-categories is a weak equivalence (hence an equivalence)
if and only if it projects to an isomorphism in the localized category.  Look
at the images of the above diagram in the localized category of $n-1$-precats
after applying the operation $(\cdot )_{p/}$.  The equivalences that we know
show that the bottom left arrow goes to an arrow which has a left and right
inverse. It follows that it goes to an invertible arrow, i.e. an isomorphism, in
the localized category. Its left inverse, the left vertical arrow, must also
go to an isomorphism. This implies that for each $p$, the map
$Cat(A)_{p/} \rightarrow F'(A)_{p/}$ is an equivalence. By definition then
$A\rightarrow F(A)$  is a weak equivalence.
\eop

\numero{Compatibility with products}

The goal of the present section is to prove the following theorem.

\begin{theorem}
\label{ce}
Suppose $A$ and $B$ are $n$-precats. Then the morphism
$$
A\times B \rightarrow Cat(A)\times Cat(B)
$$
is a weak equivalence.
\end{theorem}

Before getting to the proof, we give some corollaries.

\begin{corollary}
\label{ProdInterval}
Suppose $B$ is an $n$-precat. Let $\overline{I}$ be the $1$-category with two
isomorphic objects denoted $0$ and $1$, considered as an $n$-precat. Then the
morphisms
$$
Cat(B)\stackrel{i_0, i_1}{\rightarrow } Cat (B\times \overline{I})
\stackrel{p}{\rightarrow} Cat(B)
$$
are equivalences of $n$-categories, where $i_0$ and $i_1$ come from the
inclusions $0\rightarrow \overline{I}$ and $1\rightarrow \overline{I}$ and $p$
comes from the projection on the first factor.
\end{corollary}
{\em Proof:}
Note that the morphism $\overline{I} \rightarrow Cat(\overline{I})$
is a weak equivalence by Lemma \ref{equiv}.  Thus Theorem \ref{ce}
says that $B\times \overline{I} \rightarrow Cat (B)\times \overline{I}$
is a weak equivalence. On the other hand, the morphism
$Cat(B)\times \overline{I} \rightarrow Cat(B)$ is a weak equivalence,
so $B\times \overline{I} \rightarrow Cat(B)$ is a weak equivalence.  The
morphism $B\rightarrow Cat(B)$ is of course a weak equivalence, so
Lemma \ref{remark} implies that the two morphisms
$$
B\stackrel{i_0, i_1}{\rightarrow } B\times \overline{I}
\stackrel{p}{\rightarrow} B
$$
are weak equivalences, which is the same statement as the  corollary.
\eop

\begin{corollary}
\label{forInternalHom}
Suppose $A\rightarrow A'$ is a weak equivalence. Then for any $B$,
$A\times B\rightarrow A'\times B$ is a weak equivalence.
\end{corollary}
{\em Proof:}
By Theorem \ref{ce} we have that
$$
A\times B \rightarrow Cat(A)\times Cat(B), \;\;\;
A'\times B \rightarrow Cat(A')\times Cat(B)
$$
are weak equivalences. By hypothesis, the map $Cat(A)\rightarrow Cat(A')$ is an
equivalence of $n$-categories.  It follows
(from any of several characterizations of equivalences of $n$-categories,
see for
example \cite{Tamsamani} Proposition 1.3.1) that
$Cat(A)\times Cat(B)\rightarrow
Cat(A')\times Cat(B)$  is an equivalence, which gives the corollary.
\eop

\subnumero{Proof of Theorem \ref{ce}}
 We start by making some preliminary
reductions.  First we claim that it suffices to prove that if $B$ is any
$n$-precat and $\Sigma = \Sigma  (M,[m], \langle k , k+1 \rangle )$ and
$h= h(M,m, 1^{k+1})$ then  the morphism
$\Sigma \times B \rightarrow h\times B$ is a weak equivalence.  Suppose that we
know this statement. Then, noting that the proof of \ref{pushout} below doesn't
use Theorem \ref{ce} in degree $n$ in the case of a pushout by a trivial
cofibration which is an isomorphism on objects (which is the case for
$\Sigma \times B \rightarrow h\times B$) we obtain that for any $A$ and any
morphism $\Sigma \rightarrow A$ the morphism
$$
A\times B \rightarrow A\times B \cup ^{\Sigma \times B} h\times B
$$
is a weak equivalence.  The morphism $A\times B \rightarrow Cat(A)\times Cat(B)$
is obtained by iterating operations of this form (either on the variable $A$ or
on the variable $B$ which works the same way).  Therefore it would follow from
the hypothesis of our claim that the morphism of \ref{ce} is a weak equivalence.

We are now reduced to proving
that $\Sigma \times B \rightarrow h\times B$ is a weak equivalence.  In the
previous notations if $M$ has length strictly greater than $0$ then the
$\Sigma _{p/} \rightarrow h_{p/}$ are weak equivalences. Thus
$$
(\Sigma \times B)_{p/} = \Sigma _{p/} \times B_{p/}
\rightarrow h_{p/} \times B_{p/} = (h\times B)_{p/}
$$
is a weak equivalence, and by Lemma \ref{partialCat1} it follows that $\Sigma
\times B \rightarrow h\times B$ is a weak equivalence. Thus we are reduced to
treating the case where $M$ has length zero, that is
$$
\Sigma = \Sigma  ([m], \langle k , k+1 \rangle )\;\;\;
\mbox{and} \;\;\; h= h(m, 1^{k+1}).
$$
Let $\Sigma ^{nu}=\Sigma ^{nu}([m], 1^{k+1}) $ denote the pushout of $m$ copies
of $h(1, 1^{k+1})$ over the standard $h(0)$.  We claim that it suffices to prove
that  $$
\Sigma ^{nu} \times B \rightarrow h\times B
$$
is a weak equivalence.  This claim is proved by induction on $k$.  Note that
$\Sigma$ is obtained from $\Sigma ^{nu}$ by a sequence of standard cofibrations
over $\Sigma ' \rightarrow h'$ which are for smaller values of $k$. Assuming
that we have treated all of the cases $\Sigma ^{nu}\times B\rightarrow h\times
B$ and assuming our present claim for smaller values of $k$, we obtain that
$\Sigma ^{nu}\times B \rightarrow \Sigma \times B$ is a trivial cofibration.
It follows from Lemma \ref{remark} that $\Sigma \times B \rightarrow h\times B$
will be a trivial cofibration.

We are now reduced to proving that the morphism
$$
\Sigma ^{nu}([m], 1^k) \rightarrow h(m, 1^k)
$$
induces a weak equivalence $\Sigma ^{nu}\times B \rightarrow h\times B$ (we have
changed the indexing $k$ here from the previous paragraphs).  Our next
reduction is based on the following observation. Suppose we know this
statement for  $n$-precats $A$, $B$ and $C$ with cofibrations $A\rightarrow B$
and $A\rightarrow C$. Let $P:= B\cup ^AC$.  The morphisms
$\Sigma ^{nu}\times A \rightarrow h\times A$ (resp.
$\Sigma ^{nu}\times B \rightarrow h\times B$, $\Sigma
^{nu}\times C \rightarrow h\times
C$) are weak equivalences inducing isomorphisms on objects. As remarked
above, the proof of \ref{pushout} below doesn't use Theorem \ref{ce} in degree
$n$ when concerning pushouts by weak equivalences which are isomorphisms on
objects. On the other hand, the morphism
$\Sigma ^{nu}\times P\rightarrow h\times P$ may be obtained as a successive
coproduct by these previous morphisms; thus it is a weak equivalence.

We may apply this observation to infinite iterations of cofibrant pushouts.
But note that any $n$-precat $B$ may be expressed as an iterated pushout of
representable objects $h(M)$ over the boundaries $\partial h(M) \rightarrow
h(M)$. The boundaries are in turn iterated pushouts over representable objects.
>From the remark of the previous paragraph, it follows that we are reduced to
proving that $\Sigma ^{nu}\times B\rightarrow h\times B$ is a weak equivalence
when $B$ is a representable object. We will write $B=h(u , M)$,
distinguishing the first variable from the rest.

We next define the following operation. Suppose $C$ is an
$n-1$-precat
and suppose $D$ is a $1$-precat; then we define $D\oplus C$ to be the
$n$-precat
with $(D\oplus C)_0:= D_0$ and
for any $p$, $(D\oplus C)_{p/}$ is the union for $f\in D_p$ of
$C$ when $f$ is not totally degenerate and $\ast$ when $f$ is totally degenerate
(we say that $f\in D_p$  is {\em totally degenerate} if it is in the image of
the morphism $D_0 \rightarrow D_p$).
The morphisms of functoriality for $p\rightarrow q$ are obtained by projecting
$C$ to the final $n-1$-precat $\ast$ for elements $f\in D_q$ going to
totally degenerate elements in $D_p$.

This notation is useful because
$h(u , M)= h(u ) \oplus h(M)$. Thus $h(m, 1^k) = h(m)\oplus h(1^k)$;
and finally
$$
\Sigma ^{nu}([m], 1^k)= \Sigma ^{nu} ([m]) \oplus h(1^k).
$$

The operation $\oplus$ is compatible with pushouts:
if $B\leftarrow A \rightarrow C$ is a diagram of $n-1$-precats then for any
$1$-precat $X$,
$$
X\oplus (B\cup ^AC) = (X\oplus B) \cup ^{X\oplus A} (X\oplus C).
$$
On the other hand, if $C\rightarrow C'$
is a weak equivalence then $X\oplus C \rightarrow X\oplus C'$ is a weak
equivalence (by applying \ref{partialCat1}).  The object $h(M)$ is weak
equivalent to a pushout of objects of the form $h(1^k)$, with the pushouts
being over boundaries which are themselves weak equivalent to pushouts of
objects of the same form.  Since, as we have seen above, changing things by
pushouts or by weak equivalences (which are isomorphisms on objects) in the
second variable, preserves the statement in question.
Thus it suffices to treat, in the previous notations, the case $M=1^j$.
We have now boiled down to the basic case which needs to be treated: we must
show that
$$
\Sigma ^{nu}([m], 1^k) \times h(u, 1^j)  \rightarrow h(m, 1^k)\times h(u,1^j)
$$
is a weak equivalence.

To do this, use the standard subdivision of the product
of simplices $h(m)\times h(u)$ into a coproduct of simplices identified over
their boundaries.
\footnote{
This is basically the only place in the paper where we really use the fact that
we have taken the category $\Delta$ and not some other category such as
the semisimplicial category or a
truncation of $\Delta$ using only objects $m$ for $m\leq m_0$.
One can see for example that the statement \ref{ce} for $1$-precats
is no longer true if we try to replace $\Delta$ by the semisimplicial category
throwing out the degeneracy maps---this example comes down to saying that the
product of two free monoids on two sets of generators is not the free monoid on
the product of the sets of generators.}
In this last part of the proof there was an error in version 1: on p. 31
the line claiming that ``$B^{(a,b)} = B^{<(a,b)} \cup ^{B^{\hat{x}}}
B^x$'' was not true. Furthermore the notation of this part of the proof
was relatively difficult to follow. Thus we rewrite things in the
present version 2. 

This error and its correction were found during discussions with R.
Pellissier, so I would like to thank him. 

The basic idea remains the same as what was said in version 1. The
objects of $h(m, 1^k)\times h(u, 1^j)$ may be denoted by $(i,j)$ with
$i=0,\ldots , m$
and $j=0,\ldots , u$. These should be arranged into a rectangle. The
problem is to understand the composition as we go from $(0,0)$ to 
$(m,u)$. There are many different paths (i.e. sequences of points which
are adjacent on the grid and where $i$ and $j$ are nondecreasing) and the goal is to say
that the composition along the various different paths is the same. 
The reason is that when one changes the path by the smallest amount
possible at a single square, that is to say changing ``up then over'' to
``over then up'', the composition doesn't change. This elementary step
was done correctly in the original proof (see in v1 the statements that
the  morphisms $A^{\hat{x}}\rightarrow B^{\hat{x}}$ and $A^x\rightarrow
B^x$ are weak equivalences). Then one has to put these elementary steps
together to conclude the desired statement for the big rectangular grid. This
requires an inductive argument along the lines of what was done in v1
but in v1 the induction wasn't organized correctly and was based on a
mistaken claim as pointed out above. 

So let's rewrite things and hope for the best! Put
$$
A:= \Sigma ^{nu}([m], 1^k) \times h(u, 1^j)
$$
and
$$
B:= h(m, 1^k)\times h(u,1^j).
$$
Note that $A$ is a coproduct of things of the form
$h(1,1^k) \times h(u, 1^j)$. 
We want to show that the morphism $A\rightarrow B$ is a weak
equivalence. The precats $A$ and $B$ share the same set of objects which
we denote by $Ob$, equal to the set of pairs $(i,j)$ with $0\leq i
\leq m$ and $0 \leq j\leq u$. Suppose $S \subset Ob$ is a subset of
objects. We denote by $A\{ S\} $ (resp. $B\{ S\}$) the full sub-precat of
$A$ (resp. $B$) whose object-set is $S$. By ``full sub-precat'' we mean that for any
sequence $x_0, \ldots , x_k$ in $S$, 
$$
A\{ S\} _{k/}(x_0,\ldots , x_k):= A _{k/}(x_0,\ldots , x_k)
$$
and the same for $B\{ S\}$. We will use this for subsets $S$ of the form
``notched sub-rectangle plus a tail''. By a ``sub-recatangle'' we mean a subset
of the form
$$
S' = \{ (i,j) :  \;\;\; 0 \leq i \leq m',\;\;\ 0\leq j \leq u' \}
$$
and by a ``tail'' we mean a subset of the form 
$$
S'' = \{ (i_k, j_k) : \;\;\; 0 \leq k \leq r \}
$$
where $i_k \leq i_{k+1} \leq i_k + 1$  and $j_k \leq j_{k+1} \leq j_k +
1$. A tail $S''$ that goes with a rectangle $S'$ as above, is assumed to
 have $i_0= m'$, $j_0 = u'$, $i_r = m$, $j_r=u$. In other words, the
tail is a path going from the upper corner of $S'$ to the upper corner
of $Ob$, and the path goes by steps of at most one in both the $i$ and
the $j$ directions. 

Finally, a ``notched sub-rectangle'' is a subset of the form 
$$
S' = \{ (i,j) : \;\;\; (0 \leq i \leq m' \, \mbox{and} \, 0 \leq j
\leq u'-1 ) \, \mbox{or} \, (v'\leq i \leq m' \, \mbox{and}
j=u')\} .
$$
We call $(m',u', v')$ the {\em parameters} of $S'$, and if necessary we
denote $S'$ by $S'(m',u',v')$. Note that $u'\geq 1$
here,
and $0\leq v' \leq m'$. A rectangle with $u'=0$ may be considered as
part of a tail; thus, modulo the initial case where all of $S$ is a
tail, which will be treated below, it is safe to assume $u' \geq 1$. 
If $v'=0$ then $S'$ is just the rectangle of size $m'
\times u'$. If $v'=m'$ then $S'$ is a rectangle of size $m'\times
(u'-1)$, plus the first segment of a tail going from $(m', u'-1)$ to
$(m',u')$. Thus if $S''$ is a tail from $(m',u')$ to $(m,u)$,
we get that 
$$
S'(m',u',0) \cup S''
$$ 
is a rectangle of size $m'\times u'$ plus a tail, whereas 
$$
S'(m',u',m') \cup S'' 
$$ 
is a rectangle of size $m'\times (u'-1)$ plus a tail.

We prove by induction on $(m',u')$ that if $R$  is of the form $R=S'\cup S''$ for $S'$
a rectangle of size $m'\times u'$ and $S''$ a tail from $(m',u')$ to
$(m,u)$,
then $A\{ R\} \rightarrow B\{ R\}$ is a weak equivalence. We treat the
case of $m',u'$ and suppose that it is known for all strictly smaller
rectangles (i.e. for  $(m'', u'')\neq (m',u')$ with
$m'' \leq m'$ and $u''\leq u'$) and all tails. 
In the current part of the induction we assume that $u' \geq 1$. The
case $u'=0$ (which is really the case of a tail $S''$ going all the way from $(0,0)$
to $(m,u)$) will be treated below.

In particular, we know
that
$$
A \{ S'(m',u',m') \cup S'' \} \rightarrow B \{ S'(m',u',m') \cup S'' \}
$$
is a weak equivalence (cf the above description of
$S'(m',u',m')$). Now we treat the case where $S$ is a notched
rectangle plus tail of the form 
$$
S = S'(m',u',v') \cup S'' 
$$
for $0\leq v' \leq m'$.
For an $S$ of
this notched form, we claim again that $A\{ S\} \rightarrow B\{ S\}$ is a weak
equivalence. We prove this by descending induction on $v'$,
the initial case $v'=m'$ being obtained above from
the case of rectangles of size $m'\times (u'-1)$.
Thus we 
may fix $v'< m'$ and assume that it is known for 
$$
\overline{S}=  S'(m',u',v'+1) \cup S'',
$$
in other words we may assume that $A\{ \overline{S}\} \rightarrow B\{
\overline{S}\}$
is a weak equivalence.

We analyze how to go from $\overline{S}$ to $S$. Note that $S$ has
exactly one object more than $\overline{S}$, the object 
$$
x := (v',u').
$$
Let $S^x$ denote the subset of objects $(i,j)\in S$ such that either 
$(i,j)\leq x$ or $(i,j) \geq x$. Here we define the order relation by
$$
(i,j) \leq (k,l) \Leftrightarrow i\leq k \;\; \mbox{and}\;\; j \leq l.
$$
With respect to this order relation, note that for a sequence of objects 
$(x_0, \ldots , x_p)$, we have that 
$B_{p/} (x_0,\ldots , x_p)$ is nonempty,
only if $x_0\leq x_1 \leq \ldots \leq x_p$. (The same remark holds {\em
a fortiori} for $A$ because of the map $A\rightarrow B$.) In particular, if
$(x_0,\ldots , x_p)$ is a sequence of objects of $S$ such that 
$B_{p/} (x_0,\ldots , x_p)$ is nonempty and such that some $x_a=x$ then
all of the $x_b$ are in $S^x$. 

Let $\overline{S}^x = \overline{S} \cup S^x$. We claim $(\ast )$ that
$$
A\{ S\} = A\{ \overline{S}\} \cup ^{A\{ \overline{S}^x\} } A\{ S^x \} ,
$$ 
and similarly that 
$$
B\{ S\} = B\{ \overline{S}\} \cup ^{B\{ \overline{S}^x\} } B\{ S^x \} .
$$ 
These are the statements that replace the faulty lines of the proof in
version 1.  To prove the claim, suppose $(x_0, \ldots , x_p)$ is a
sequence of objects of $S$. It suffices to show that
$$
B\{ S\}_{p/}(x_0,\ldots ,x_p) 
$$
is the pushout of 
$$
B\{ \overline{S}\}_{p/}(x_0,\ldots ,x_p)
$$
and
$$
B\{ S^x\}_{p/}(x_0,\ldots ,x_p)
$$
over 
$$
B\{ \overline{S}^x\}_{p/}(x_0,\ldots ,x_p)
$$
(the proof is the same for $A$, we give it for $B$ here).

If none of the objects $x_a$ is equal to $x$,
then either the sequence stays inside $\overline{S}^x$, in which case:
$$
B\{ S\}_{p/}(x_0,\ldots ,x_p) = 
B\{ \overline{S}\}_{p/}(x_0,\ldots ,x_p)=
B\{ \overline{S}^x\}_{p/}(x_0,\ldots ,x_p)=
B\{ S^x\}_{p/}(x_0,\ldots ,x_p);
$$
or else the sequence doesn't
stay inside $\overline{S}^x$, in which case 
$$
B\{ S\}_{p/}(x_0,\ldots ,x_p) = 
B\{ \overline{S}\}_{p/}(x_0,\ldots ,x_p)
$$
but
$$
B\{ \overline{S}^x\}_{p/}(x_0,\ldots ,x_p)=
B\{ S^x\}_{p/}(x_0,\ldots ,x_p)=\emptyset .
$$
In both of these cases one obtains the required pushout formula. On the
other hand, if some $x_a$ is equal to $x$, then either the sequence
doesn't stay inside $S^x$ in which case all terms are empty (cf the
above remark), or else it stays inside $S^x$ in which case 
$$
B\{ S\}_{p/}(x_0,\ldots ,x_p) = B\{ S^x\}_{p/}(x_0,\ldots ,x_p)
$$
but 
$$
B\{ \overline{S}\}_{p/}(x_0,\ldots ,x_p)=
B\{ \overline{S}^x\}_{p/}(x_0,\ldots ,x_p)= \emptyset .
$$
Again one obtain the required pushout formula.
This proves the claim $(\ast )$.

We now note that both $S^x$ and 
$\overline{S}^x$ are of the form, a rectangle of size $v'\times (u'-1)$,
plus a tail. The first step of the tail for $S^x$ goes from 
$(v',u'-1)$ to $x=(v',u')$. The next step goes to $(v'+1,u')$.
The first step of the tail for $\overline{S}^x$ goes from $(v',u'-1)$
directly to $(v'+1, u')$. Both tails continue horizontally from 
$(v'+1,u')$ to $(m',u')$ and then continue as the tail $S''$ from there
to $(m,u)$. By our induction hypothesis, we know that 
$$
A\{ S^x\} \rightarrow B\{ S^x\}
$$
and
$$
A\{ \overline{S}^x\} \rightarrow B\{ \overline{S}^x\}
$$
are weak equivalences. Recall from above that we also know that 
$A\{ \overline{S}\} \rightarrow B\{ \overline{S}\}$ is a weak
equivalence. These morphisms from full sub-precats of $A$ to full
sub-precats of $B$ are all isomorphisms on sets of objects, so by the
remark at the start of the proof of Theorem \ref{ce}, we can use Lemma 
\ref{pushout} for pushouts along these morphisms. 
A standard argument shows that
the morphism
$$
A\{ \overline{S}\} \cup ^{A\{ \overline{S}^x\} }A\{ S^x\}
\rightarrow 
B\{ \overline{S}\} \cup ^{B\{ \overline{S}^x\} }B\{ S^x\} 
$$
is a weak equivalence. For clarity we now give this standard argument.
First,
$$
A\{ \overline{S}\} \rightarrow
A\{ \overline{S}\} \cup ^{A\{ \overline{S}^x\} }B\{ \overline{S}^x\} 
$$
is a trivial cofibration (which again induces an isomorphism on sets of
objects).
One can verify that the morphism
$$
A\{ \overline{S}\} \cup ^{A\{ \overline{S}^x\} }B\{ \overline{S}^x\}
\rightarrow B\{ \overline{S}\}
$$
is a cofibration. It is a weak equivalence by Lemma 3.8. Thus it is a
trivial cofibration inducing an isomorphism on sets of objects. Similarly, 
the morphism
$$
B\{ \overline{S}^x\} \cup ^{A\{ \overline{S}^x\} }A\{ S^x\}
\rightarrow B\{ S^x\}
$$
is a trivial cofibration inducing an isomorphism on sets of objects.
Thus the pushout morphism (pushing out by these two morphisms at once)
$$
(A\{ \overline{S}\} \cup ^{A\{ \overline{S}^x\} }B\{ \overline{S}^x\} )
\cup ^{B\{ \overline{S}^x\} } 
(B\{ \overline{S}^x\} \cup ^{A\{ \overline{S}^x\} }A\{ S^x\} )
$$
$$
\rightarrow 
B\{ \overline{S}\} \cup ^{B\{ \overline{S}^x\} }B\{ S^x\} 
$$
is a weak equivalence. 
Finally, note that 
$$
(A\{ \overline{S}\} \cup ^{A\{ \overline{S}^x\} }B\{ \overline{S}^x\} )
\cup ^{B\{ \overline{S}^x\} } 
(B\{ \overline{S}^x\} \cup ^{A\{ \overline{S}^x\} }A\{ S^x\} )
$$
$$
= (A\{ \overline{S}\} \cup ^{A\{ \overline{S}^x\} }A\{ S^x\} )
\cup  ^{A\{ \overline{S}^x\} }  B \{ \overline{S}^x\} 
$$
so the morphism from 
$A\{ \overline{S}\} \cup ^{A\{ \overline{S}^x\} }A\{ S^x\} $
to this latter, is also a weak equivalence. Putting these all together
we have shown that the morphism
$$
A\{ S\} = A\{ \overline{S}\} \cup ^{A\{ \overline{S}^x\} }A\{ S^x\}
$$
$$
\rightarrow 
B\{ \overline{S}\} \cup ^{B\{ \overline{S}^x\} }B\{ S^x\} 
= B\{ S\}
$$
is a weak equivalence. This completes the proof of the inductive step
for the  descending induction on $v'$, so we obtain the result for $v'=0$, in which case
$S$ is a rectangle of size $m'\times u'$ plus a tail; in turn, this
gives the inductive step for the induction on $(m',u')$.

After all of this induction we are left having to treat the initial case $u'= 0$,
where
all of $S$ is a tail going from $(0,0)$ to $(m,u)$.
This part of the proof is exactly the same as in version 1: what we call
the ``tail'' here corresponds to the sequence which was denoted $x$ in
version 1. If $S$ is a tail, then $A\{ S\} $ and $B\{ S\}$ are $1$-free
ordered $n$-precats, so by Lemma  \ref{freeness} it suffices to check
that for two adjacent objects $x,y$ in the sequence corresponding to
$S$, the morphisms 
$$
A_{1/}(x,y)\rightarrow B_{1/}(x,y)
$$
are weak equivalences. In fact these morphisms are isomorphisms.
If $x=(i,j)$ and $y=(i,j+1)$ then
$$
A_{1/}(x,y) = B_{1/}(x,y) = h(1^j);
$$
if $x=(i,j)$ and $y=(i+1,j)$ then
$$
A_{1/}(x,y) = B_{1/}(x,y) = h(1^k);
$$
and if $x=(i,j)$ and $y=(i+1,j+1)$ then
$$
A_{1/}(x,y) = B_{1/}(x,y) = h(1^k)\times  h(1^j).
$$
Thus the criterion of \ref{freeness} implies that $A\{ S\}\rightarrow
B\{ S\}$ is a weak equivalence. This completes the initial case of the
induction.

Combining this initial step with the inductive step that was carried out above,
we obtain the result in the case
where $S=Ob$ is the whole rectangle of size $m\times u$; in this case $A\{
S\} = A$ and $B\{ S\} = B$, so we have completed the proof that
$$
A=\Sigma ^{nu}([m], 1^k) \times h(u, 1^j)  \rightarrow h(m, 1^k)\times h(u,1^j)=B
$$
is a weak equivalence. This completes the proof of Theorem \ref{ce}.
\eop

The above proof went basically along the same lines as the proof of version
1, but here we met {\em all} possible tails going from $(0,0)$
to $(m,u)$ along the way in our induction, whereas in  version 1 only
some of the tails were met. This should have been taken as an indication
that the proof in version 1 was incorrect.  I would like again to thank R.
Pellissier for an ongoing careful reading which turned up this problem.
His reading has also turned up numerous other problems in the exposition
or organisation of the argument (the reader has no doubt noticed!); 
however, I have chosen in the present version 2 to make only a minimalist
correction of the above problem.


\numero{Proofs of the remaining lemmas and Theorem \ref{cmc}}

We can assume \ref{equiv}, \ref{partialCat1} and the
corollary \ref{partialCat} for degree $n$ also.

{\em Proof of \ref{coherence}:}
We have to prove that the morphism $f:Cat (A) \rightarrow Cat(Cat(A))$ is
an equivalence. Note that this is not the same morphism as the standard
inclusion, rather it is the morphism induced by $A\rightarrow Cat (A)$. In
particular, \ref{coherence} is not just an immediate corollary of \ref{equiv}.
To obtain the proof, note that $Cat(A)$ is marked, so we have a morphism
$$
r:Cat (Cat(A))\rightarrow Cat (A)
$$
of marked easy $n$-categories, inducing the identity on the standard
map $i:Cat (A) \rightarrow Cat (Cat (A))$.  The morphism $r$ is an equivalence
because the standard map $i$ is an equivalence by \ref{equiv}. On the other
hand, the morphism of $n$-precats $A\rightarrow Cat (A)$ induces the morphism
$f$ of marked easy $n$-categories.   We obtain a morphism $r\circ f$ of marked
easy $n$-categories $Cat (A)\rightarrow Cat (A)$ extending the standard map
$A\rightarrow Cat (A)$.  By the universal property of $Cat (A)$, $r\circ f$ is
the identity.  Thus (applying our usual Lemma \ref{remark}) the
morphism $f$ is an equivalence.

{\em Proof of \ref{pushout}:}
We first treat the following special case of \ref{pushout}. Say that a morphism
of $n$-categories $A\rightarrow B$ is an {\em iso-equivalence} if it is an
equivalence and an isomorphism on objects. This is equivalent to the condition
that for all $m$, the morphism $A_{m/}\rightarrow B_{m/}$ is an equivalence of
$n-1$-categories.

\begin{lemma}
\label{isopushout}
Suppose $A\rightarrow B$ is an iso-equivalence of $n$-categories and
$A\rightarrow C$ is a morphism of $n$-categories. Then the morphism
$$
C \rightarrow Cat (B\cup ^AC)
$$
is an equivalence (in fact, even an iso-equivalence).
\end{lemma}
{\em Proof:}
By our inductive hypotheses, the morphisms
$$
C_{m/} \rightarrow Cat (B_{m/} \cup ^{A_{m/}}C_{m/})
$$
are equivalences of $n-1$-categories.  Setting $D= B\cup ^AC$ we have
$$
D_{m/}=B_{m/} \cup ^{A_{m/}}C_{m/}
$$
and note that $D_0=C_0$.  Hence the morphism of $n$-precats
$$
C\rightarrow Cat _{\geq 1}(D)
$$
is an equivalence on the level of  each  $C_{m/}$.  But the condition of
being an
$n$-category depends only on the equivalence type of the $C_{m/}$, in particular
$Cat _{\geq 1}(D)$ is an $n$-category (in this special case only---this is not a
general principle!). Note in passing that the morphism
$$
C \rightarrow Cat _{\geq 1}(D)
$$
is an equivalence of $n$-categories. By Lemma \ref{equiv} in degree $n$,
the morphism
$$
Cat _{\geq 1}(D)\rightarrow Cat (Cat _{\geq 1}(D))
$$
is an equivalence.  On the other hand, by Lemma \ref{partialCat} in degree $n$
the morphism
$$
Cat (D) \rightarrow Cat (Cat _{\geq 1}(D))
$$
is an equivalence.  By Lemma \ref{remark} at the start of the proof of all
of the
lemmas, this shows that $C\rightarrow Cat(D)$ is an equivalence. This proves
Lemma \ref{isopushout}.
\eop

The next lemma is the main special case which has to be treated by hand.
This proof uses Corollary \ref{ProdInterval}, which in turn is where we use the
full simplicial structure of $\Delta$.
It seems likely (from some considerations in topological examples) that the
following lemma would not be true if we looked only at functors
from $\Delta_{\leq k}$ to $n$-precats.

\begin{lemma}
\label{specialcase}
Suppose $\overline{I}$ is the category with two objects and
exactly one isomorphism between them. Let $0$ denote one of the objects. Then
for any $n$-category $A$ and object $c\in A_0$, if $0\rightarrow A$ denotes
the corresponding morphism, the push-out morphism
$$
A \rightarrow Cat (\overline{I} \cup ^0 A)
$$
is an equivalence.
\end{lemma}
{\em Proof:}
Let $\overline{I}^2$ denote $\overline{I}\times \overline{I}$.
There is a morphism $h:\overline{I}^2\rightarrow
\overline{I}\times \overline{I}$ equal to the identity on
$\overline{I}\times \{
0\}$ and on $\{ 0\} \times \overline{I}$, and sending $\overline{I}\times \{
1\}$ to $(0, 1)$.  Let $B:=  \overline{I} \cup ^0A$. Then
$$
B\times \overline{I} = \overline{I}^2 \cup ^{\{ 0\} \times \overline{I}}C\times
\overline{I}.
$$
Using $h$ on the first part of this pushout we obtain a map
$$
f:B\times \overline{I}\rightarrow B\times \overline{I}
$$
such that $f|_{B\times \{ 0\}}$ is the identity and $f|_{B\times \{ 1\}}$
is the projection $B\rightarrow A$ obtained from the projection
$\overline{I}\rightarrow \{ 0\}$. By Corollary \ref{ProdInterval}, the morphisms
$$
i_0:Cat (B)\times \{ 0 \} \rightarrow Cat (B \times \overline{I})
$$
and
$$
i_1:Cat (B)\times \{ 1 \} \rightarrow Cat (B \times \overline{I})
$$
are equivalences of $n$-categories. Next, the morphism $f$ induces a morphism
$$
g:Cat (B\times \overline{I})\rightarrow Cat (B)
$$
such that the composition with $i_0$ is the identity $Cat(B)\rightarrow Cat(B)$
and the composition with $i_1$ is the factorization $Cat(B)\rightarrow
Cat(A)\rightarrow Cat(B)$. Looking at $g\circ i_0$ we conclude that $g$ is
an equivalence of $n$-categories. Therefore (since $i_1$ is an equivalence)
the composition $g\circ i_1$ is an equivalence of $n$-categories.  Now we have
morphisms $Cat(A)\rightarrow Cat(B)$ and $Cat(B)\rightarrow Cat(A)$ such
that the composition in  one direction is the identity, and the composition in
the other direction is an equivalence of $n$-categories. This implies that the
two morphisms are equivalences of $n$-categories by Lemma \ref{remark}.
\eop

In preparation for the next corollary, we discuss a sort of
``versal semi-interval'' $\overline{J}$. Ideally we would like to have an
$n$-precat which is weakly equivalent to $\ast$, containing two objects $0$ and
$1$ such that for any $n$-category (easy, perhaps) $A$ with two equivalent
objects $a$ and $b$, there exists a morphism from our ``interval'' to $A$
taking $0$ to $a$ and $1$ to $b$.  I didn't find an easy way to make this
construction. The problem is somewhat analogous to the problem of finding a
canonical inverse for a homotopy equivalence, solved in a certain topological
context in \cite{flexible} but which seems quite complicated to put into action
here in view of the fact that our $n$-category $A$ might not be fibrant (we
don't yet have the closed model structure!).  Thus we will be  happy with a
cruder version. Let $\overline{J}$ be the universal easy $n$-category
with two objects $0$ and $1$ and a ``marked inner equivalence'' $u :0\rightarrow
1$. The quasi-inverse of $u$ will be denoted by $v$. The marking means a
structure of choice of morphism whenever necessary for the definition of
inner equivalence, as well as a choice of diagram (i.e. a partial marking in
the sense defined at the start of the paper) whenever necessary for things to
make sense.  In practice this means that we start with objects $0$ and $1$,
add the morphisms $u$ and $v$, add the diagrams over $2\in \Delta$ mapping to
$(u,v)$ and $(v,u)$, and (letting $w$ and $y$ denote the compositions resulting
from these diagrams) add (inductively by the same construction for
$n-1$-categories) equivalences between $w$ and $e$ (resp. $y$ and $e$) where $e$
denote the identities.  Let $\overline{L}\subset \overline{J}$ be the
full-sub-$n$-category whose object set is $\{ 0\}$.  The morphism
$\overline{L}\rightarrow \overline{J}$ is automaticallly an equivalence since
it is an isomorphism on morphism $n-1$-categories and is essentially surjective
since by construction $1\in \overline{J}$ is equivalent to $0$. By the
universal property of $\overline{J}$ we obtain a morphism
$\overline{J}\rightarrow \overline{L}$ sending $u$ to $e$ and $v$ to $w$,
sending our $2$-diagrams to degenerate diagrams in $\overline{L}$ and sending
our homotopies to the corresponding homotopies in $\overline{L}$. The
composition
$$
\overline{L}\rightarrow
\overline{J}\rightarrow \overline{L}
$$
is the identity.  On the other hand we have an obvious map
$\overline{J}\rightarrow \overline{I}$, so we obtain a map
$$
\overline{J}\rightarrow \overline{L}\times \overline{I}.
$$
This map is  compatible with the inclusions of
$\overline{L}$ and hence is an equivalence of $n$-categories.
It is also an isomorphism on objects so it is an iso-equivalence.

\begin{corollary}
\label{specialJ}
Let $\overline{L}\subset \overline{J}$ be as above.
Then for any
$n$-category $A$ and morphism $\overline{L}\rightarrow A$, the push-out
morphism
$$
A \rightarrow Cat (\overline{J} \cup ^{\overline{L}} A)
$$
is an equivalence.
\end{corollary}
{\em Proof:}
Let $B= (\overline{L}\times \overline{I}) \cup  ^{\overline{L}} A$
and $C= \overline{J} \cup  ^{\overline{L}}
A$ The morphism $\overline{J}\rightarrow (\overline{L}\times \overline{I})$
is an iso-equivalence so  it satisfies the hypothesis of \ref{partialCat1}.
By \ref{pushout} for $n-1$-precats, the morphism $C\rightarrow B$ also
satisfies the hypothesis of \ref{partialCat1}. Therefore the morphism
$Cat(C)\rightarrow Cat(B)$ is an equivalence. It suffices to show that
$A\rightarrow Cat(B)$ is
an equivalence. For this, note that
$$
\overline{L} \cup ^0\overline{I} \rightarrow \overline{L}\times \overline{I}
$$
is a weak equivalence by Lemma \ref{specialcase}. Similarly
$$
A \rightarrow A \cup ^0\overline{I}
$$
is a weak equivalence, and
$$
B = (A\cup ^0\overline{I}) \cup ^{(\overline{L} \cup ^0\overline{I})}
\overline{L}\times \overline{I}
$$
so composing these two gives that $A\rightarrow B$ is a weak equivalence.
\eop

\begin{corollary}
Suppose $A\rightarrow B$ is a cofibrant equivalence of $n$-categories such that
the objects of $A$ form a subset of the objects of $B$ whose complement has one
object. Then  the push-out morphism
$$
C \rightarrow Cat (B \cup ^A C)
$$
is an equivalence.
\end{corollary}
{\em Proof:}
We may replace
$B$ by $Cat(B)$ since the morphism $B\rightarrow Cat(B)$ is a sequence of
standard pushouts, so the corresponding morphism on pushouts of $C$ is also a
sequence of standard pushouts so the conclusion for $Cat(B)$ implies the
conclusion for $B$ (by Lemma \ref{remark}). Thus we may assume that $B$ is an
easy $n$-category.

Let $A'$ be the full sub-$n$-category of $B$ consisting of the objects of $A$.
The pushout of $C$ from $A$ to $A'$ is a weak equivalence by Lemma
\ref{isopushout}. Thus we may assume that $A=A'$.

Let $b$ denote the single new
object of $B$.  It is equivalent to an object $a\in A$.
By the universal property of $\overline{J}$ there is a morphism
$\overline{J}\rightarrow B$ sending $0$ to $a$ and $1$ to $b$. Since $A$ is now
a full sub-$n$-category of $B$, this morphism sends $\overline{L}$ to $A$. Let
$E$ denote the push-out
$$
E := Cat (A \cup ^{\overline{L}} \overline{J}).
$$
By the prevoius corollary, $A\rightarrow E$ is an equivalence.
Our morphism $\overline{J} \rightarrow B$ gives a morphism
$$
E \rightarrow B
$$
(use the marking of $B$ to go from $Cat(B)$ back to $B$)
and this is an equivalence since
$$
A\rightarrow B \rightarrow Cat(B) \;\;\; \mbox{and}\;\;\; A \rightarrow E
$$
are equivalences. But the morphism $E\rightarrow B$ induces an isomorphism
on objects.  Now we have
$$
C \cup ^A (A \cup ^{\overline{L}} \overline{J})= C \cup
^{\overline{L}}\overline{J}
$$
so
$$
C \rightarrow Cat (C \cup ^A (A \cup ^{\overline{L}} \overline{J}))
$$
is an equivalence by the previous corollary.
It is obvious from the construction
of $Cat$ (resp. $BigCat$) via pushouts, together with the reordering of these
pushouts, that
$$
BigCat (C \cup ^A (A \cup ^{\overline{L}} \overline{J})) =
BigCat (C \cup ^ACat (A \cup ^{\overline{L}}
\overline{J})) =BigCat(C\cup ^AE).
$$
Thus (since taking $BigCat$ is equivalent to taking $Cat$ by Lemma \ref{equiv}
---which we now know---and the reordering principle)
$$
C\rightarrow Cat(C \cup ^AE)
$$
is an equivalence.
Now
$$
Cat (C\cup ^AE) \rightarrow Cat (C \cup ^A B)
$$
is an equivalence because $E_{p/}\rightarrow B_{p/}$ is an equivalence so
by \ref{pushout} in degree $n-1$,
$$
(C\cup ^AE )_{p/} \rightarrow (C \cup ^A B)_{p/}
$$
is an equivalence and by Lemma \ref{partialCat1} we get the desired statement.
Combining, we get that
$$
C\rightarrow Cat (C \cup ^A B)
$$
is an equivalence.
\eop

{\em Proof of Lemma \ref{pushout}:} Suppose $A\rightarrow B$ is a cofibration
of $n$-categories which is an equivalence.  By applying the previous corollary
inductively (adding one object at a time) we conclude that the push-out is an
equivalence.

Finally we treat the case where $A$, $B$ and $C$ are only $n$-precats rather
than $n$-categories.
If $\Sigma \rightarrow h$ is one of our standard pushout diagrams and
if $\Sigma \rightarrow A$ is  a morphism then
$$
(B\cup ^{\Sigma} h) \cup ^{A\cup ^{\Sigma} h}(C\cup ^{\Sigma} h) = (B\cup
^AC)\cup ^{\Sigma}h.
$$
This implies that
$$
BigCat(B)\cup ^{BigCat(A)}BigCat(C)
$$
is obtained by a collection of standard pushouts from $B\cup ^AC$, so in
particular (by reordering)
$$
BigCat(BigCat(B)\cup ^{BigCat(A)}BigCat(C))=BigCat(B\cup ^AC).
$$
Now our hypothesis is that $BigCat(A)\rightarrow BigCat(B)$ is an equivalence
(note also that it is a cofibration since $A\rightarrow B$ is a cofibration).
By our proof of \ref{pushout} for the case of $n$-categories (and the
equivalence between $Cat$ and $BigCat$ which we now know by \ref{equiv}) we
conclude that
$$
BigCat(C)\rightarrow BigCat(BigCat(B)\cup ^{BigCat(A)}BigCat(C))
$$
is an equivalence, which is to say that
$$
BigCat(C)\rightarrow BigCat(B\cup ^AC)
$$
is an equivalence. Thus $C\rightarrow B\cup ^AC$ is a weak equivalence.
\eop

{\em Remark:}  The semi-interval $\overline{J}$ we have constructed above is
not contractible (i.e. equivalent to $\ast$).  However for some purposes we
would  like to have such an object. We have the following fact (which is not
used in the proof of Theorem \ref{cmc} in degree $n$---but which we put here
for expository reasons):
\begin{proposition} \label{intervalK}
There is an
$n$-category $K$ such that $K\rightarrow \ast$ is an equivalence, together with
objects $0,1\in K$
such that
if $A$ is an
$n$-category and if $a,b$ are two equivalent objects of $A$ then there is a
morphism $K\rightarrow A$ sending $0$ to $a$ and $1$ to $b$.
\end{proposition}
{\em Proof:}
Since this proposition is not used in degree $n$  in the proof of Theorem
\ref{cmc}, we can apply Theorem \ref{cmc}. We say that $K$ is {\em
contractible} if the morphism $K\rightarrow \ast$ is an equivalence. In view of
the versal property of $\overline{J}$, it suffices to construct a contractible
$K$ with objects $0,1$ and a morphism $K\rightarrow \overline{J}$ sending $0$ to
$0$ and $1$ to $1$. From the original discussion of $\overline{J}$ we have
an equivalence $\overline{J}\rightarrow \overline{L}\times \overline{I}$.
Using the closed model structure, factor the constant morphism as
$$
\overline{I}\rightarrow M \rightarrow \overline{L}\times \overline{I}
$$
into a composition of a trivial cofibration followed by a fibration. Note that
$M$ is contractible. Set
$$
K:=\overline{J}\times _{\overline{L}\times \overline{I}}M.
$$
The morphism  $\overline{J} \rightarrow \overline{L}\times
\overline{I}$ is an isomorphism on objects, so for each $p$,
$\overline{J}_{p/} \rightarrow (\overline{L}\times \overline{I})_{p/}$ is an
equivalence of $n-1$-categories. Note also that
$M _{p/}\rightarrow (\overline{L}\times \overline{I})_{p/}$ are fibrations.
By the fact that weak equivalences are stable under fibrant pullbacks for
$n-1$-categories (Theorem \ref{properness}), we have that
$$
K_{p/}= \overline{J}_{p/}\times _{(\overline{L}\times \overline{I})_{p/}}M
_{p/} \rightarrow M_{p/}
$$
are weak equivalences, which in turn implies that $K\rightarrow M$ is a weak
equivalence.
In particular, $K$ is contractible. Since the morphism
$$
M \rightarrow \overline{L}\times \overline{I}
$$
is surjective on objects, there are objects $0,1\in K$ mapping to $0,1\in
\overline{J}$.  This completes the construction.
\eop

\begin{corollary}
\label{FibImpliesCat}
If $A$ is a fibrant $n$-precat then $A$ is automatically an
easy $n$-category.
\end{corollary}
{\em Proof:}
To show this it suffices to show that the morphisms
$$
\varphi :\Sigma (M, [m], \langle k,k+1 \rangle ) \rightarrow h(M,[m], 1^{k+1})
$$
are trivial cofibrations. But $\varphi$ is a cofibration which is the first step
in  an addition of arbitrary pushouts of our standard morphisms $\varphi$, so by
reordering of these pushouts the above inclusion extends to an isomorphism
$$
BigCat (\Sigma (M, [m], \langle k,k+1 \rangle ) )\cong
BigCat (h(M,[m], 1^{k+1})).
$$
Since $Cat (B) \rightarrow BigCat(B)$ is an equivalence  by Lemma \ref{equiv}
applied to $Cat(B)$ plus reordering, this implies that the morphism $\varphi$
above is a trivial cofibration.  By the
definition of fibrant, $A$ must then satisfy the extension property to be an
easy $n$-category. \eop

\subnumero{The proof of Theorem \ref{cmc}}

We follow the proof of Jardine-Joyal that simplicial presheaves form a closed
model category, as described in \cite{Jardine}.  The proof is based on the
axioms CM1--CM5 of \cite{QuillenAnnals}.

{\em Proof of CM1:} \,\, The category of $n$-precats is a category of presheaves
so it is closed under finite (and even arbitrary) direct and inverse limits.

{\em Proof of CM2:} \,\, Given composable morphisms
$$
X\stackrel{f}{\rightarrow}Y\stackrel{g}{\rightarrow}Z,
$$
if any two of $f$ or $g$ or $g\circ f$
are weak equivalences  then the same two of $Cat(f)$, $Cat(g)$ or $Cat(g\circ
f)$ are equivalences of $n$-categories in the sense of \cite{Tamsamani} and
by Lemma \ref{remark} the third is
also an equivalence; thus the third of our original morphisms is a weak
equivalence.

{\em Proof of CM3:} \,\,
This axiom says that ``the classes of cofibrations, fibrations and weak
equivalences are closed under retracts''. Jardine \cite{Jardine} doesn't
actually
discuss the retract condition other than to say that it is obvious in his case,
and a look at Quillen  yields only the conclusion that the diagram on p. 5.5 of
\cite{Quillen} for the  definition of retract is wrong (that diagram has
no content related to the word ``retract'', it just  says that one arrow is the
composition of three others). Thus---since I am not sufficiently well acquainted
with other possible references for this---we are reduced to speculation about
what Quillen means by ``retract''. Luckily enough, this speculation comes out to
be non-speculative in the end. We say that  $f: A\rightarrow B$ is a {\em weak
retract} of $g:X\rightarrow Y$  if there is a diagram
$$
\begin{array}{ccccc}
A&\stackrel{i}{\rightarrow}&X&\stackrel{r}{\rightarrow}&A\\ \downarrow
&&\downarrow && \downarrow\\
B&\stackrel{j}{\rightarrow}&Y&\stackrel{s}{\rightarrow}&B \end{array}
$$
such
that $r\circ i=1_A$ and $s\circ j=1_B$.  There is also another notion which we
call {\em strong retract} obtained  by using the same diagram but with the
arrows going in the opposite direction on the bottom.  It turns out that if $f$
is a strong retract of $g$ then $f$ is also a weak retract of $g$: the strong
retract condition can be stated as the condition $j\circ f \circ r = g$
(along with the retract conditions $ri=1$ and $sj=1$).  Applying $s$ on the
left we obtain $fr=sg$ and applying $i$ on the right we obtain $jf=gi$, these
two conditions giving the weak retract condition.  Thus for our purposes, if
we can show that the classes of maps in question are closed under weak retract,
this implies that they are also closed under strong retract, and we don't
actually care which of the two definitions was intended in \cite{Quillen}!

We start out, then, with the condition that $f$ is a weak retract of $g$ using
the notations of the diagram given above.
If $g$ satisfies any lifting property then $f$ satisfies the same lifting
property, using the retractions.  This shows that if $g$ is a fibration then
$f$ is a fibration.  Furthermore, if $g$ is a cofibration then it is injective
over any object $M=(m_1,\ldots , m_k)$ with $k<n$. It follows from the
retractions that $f$ satisfies the same injectivity conditions (one has the same
diagram of retractions on the values of all of the presheaves over the object
$M$). Thus $f$ is a cofibration.

Suppose $g$ is a weak equivalence, we would like to show that $f$ is a weak
equivalence.  Replacing the whole diagram by $Cat$ of the diagram, we may assume
that $A$, $B$, $X$, and $Y$ are $n$-categories and $g$ is an equivalence.
Suppose $x,y$ are objects of $A$. Then denoting by the same letters their
images in $B$, $X$ and $Y$ we obtain morphisms of $n-1$-categories
$f_1(x,y): A_{1/}(x,y)\rightarrow B_{1/}(x,y)$ and
$g_1(x,y): X_{1/}(x,y)\rightarrow Y_{1/}(x,y)$ such that $g$ is a retract
of $f$ in
the category of arrows of $n-1$-precats. Furthermore $g_1(x,y)$ is an
equivalence of $n-1$-categories. It follows by induction (since we may
assume CM3
known for $n-1$-categories) that $f_1(x,y)$ is an equivalence of
$n-1$-categories. In order to prove that $f$ is an equivalence we have to prove
that it is essentially surjective. Suppose $w$ is an object of $B$. Then $i(w)$
is an object of $Y$ so by essential surjectivity of $g$ there is an object $u$
of $X$ with an equivalence $e:g(u)\cong i(w)$ (i.e. a pair of elements
$e'\in Y_1(g(u),i(w))$ and $e'' \in Y_1(g(u),i(w))$ such that their
compositions, which are well defined in the truncation $T^{n-1}Y_1(
g(u),g(u))$ and $T^{n-1}Y_1(
i(w),i(w))$ are the identities in these truncations). Applying the retractions
$r$ and $s$ we obtain an element $r(u)\in A_0$ and an equivalence $s(e)$
between $fr(u)$ and $si(w)=w$. This proves essential surjectivity of $f$,
completing the verification of CM3.

{\em Proof of CM4:}\,\, The first part of CM4 is exactly Lemma
\ref{pushout}. The
second part follows from the first by the same trick as used by Jardine
(\cite{Jardine} pp 64-65) and ascribed by him to Joyal \cite{Joyal}.

{\em Proof of CM5(1):}\,\, For our situation, the cardinal $\alpha$ refered to
in Jardine is the countable infinite one $\omega$. Suppose $A\rightarrow C$ is a
trivial cofibration. We claim that if $B$ is an $\omega$-bounded subobject of
$C$ (by this we mean a sub-presheaf over $\Theta^n$) then there is an
$\omega$-bounded subobject $B_{\omega}\subset C$ as well as an $\omega$-bounded
subobject $A_{\omega}\subset A\times _BB_{\omega}$ such that $B\subset
B_{\omega} \subset C$ and such that $A_{\omega} \rightarrow
B_{\omega}$ is a trivial cofibration.  (Note that in our situation cofibrations
are not necessarily injective morphisms of presheaves, so $A_{\omega}$ is not
necessarily equal to  $A\times _BB_{\omega}$ the latter of which could be
uncountable).

To prove the claim, note that for a given element in $B_M$ for some $M$, the
statement that it is contained in an $n$-precat which is weakly equivalent to
$A$ can in principal be written out explicitly involving only a countable
number of elements of various $A_{M'}$ and $B_{M'}$.  Iterate this operation
starting with all of the elements of $B$ and repeatedly applying it to all of
the new elements that are added.  The iteration takes place a countable number
of times, and each time we add on a countable union of countable objects. At
the end we arrive at $A_{\omega} \subset B_{\omega}$ which is an
$\omega$-bounded trivial cofibration.

Using this claim, the rest of Jardine's arguments of (\cite{Jardine}, Lemmas 2.4
and 2.5) work and we obtain the statement that every morphism $f:X\rightarrow
Y$ of $n$-precats can be factored as $f=p\circ i$ where $i$ is a trivial
cofibration and $p$ is fibrant---\cite{Jardine} Lemma 2.5, which is CM5(1).
Note that the only sentence in Jardine's argument which needs further
verification is the fact that filtered colimits of trivial cofibrations are
again trivial cofibrations; and this holds in our case too.

{\em Proof of CM5(2):}
We have to prove that any morphism $f$ may be factored as $f=q\circ j$
where $q$ is a fibrant weak equivalence and $j$ a cofibration.

It suffices to construct a factorization $f=q\circ j$ with $j$ a cofibration and
$q$ a weak equivalence, for then we can apply CM5(1) to factor $q$ as
a product of a trivial cofibration and a fibration, the latter of which is
automatically also  a weak equivalence by CM2.
Thus we now search for $f=q\circ j$ with $q$ a weak equivalence and $j$ a
cofibration.

The reader may wish to think about this in the case of $1$-categories to get an
idea of what is happening and to see why this part is actually easy modulo
some small details: we multiply the number of objects in each isomorphism class
in the target category to have the morphism injective on the sets of objects.

If $f:A\rightarrow B$ is a morphism of $n$-precats then we define a
canonical factorization $A\rightarrow N(A,B)\rightarrow B$
in the following way. Let $L(A)$ denote the $1$-category (considered as an
$n$-category) whose set of objects is equal to $A_0$ and which has exactly one
morphism between any pair of objects. Note that $L(A)\rightarrow \ast$ is a weak
equivalence. The tautological map $A_0 \rightarrow L(A)_0$ lifts to a
unique map
of $n$-precats $t:A\rightarrow L(A)$.  Set $N(A,B):= L(A) \times B$ with the
diagonal map $(t,f) :A\rightarrow N(A,B)$ and the second projection
$p: N(A,B)\rightarrow B$.  Note that $p$ is a weak equivalence (by an
appropriate generalization of Corollary \ref{ProdInterval}) and $(t,f)$ is
injective on objects.

Now suppose by induction that we have constructed for every morphism
$f':A'\rightarrow B'$ of $n-1$-precats a  factorization
$A'\rightarrow M(A',B')\rightarrow B'$ as a composition of a weak  equivalence
and a cofibration, functorial in $f'$.

(To start the induction for $n=0$ we set
$M(A',B'):= B'$ recalling that all morphisms are cofibrations in this case.)

Suppose $f: A\rightarrow C$ is a morphism of $n$-precats such
that $A_0 \hookrightarrow C_0$ is injective.  Define a presheaf on $\Delta
\times
\Theta ^{n-1}$ denoted  $P(A,C)$, with factorization  $A\rightarrow
P(A,C)\rightarrow C$ as follows.  Put
$$
P(A,C)_{p/} := M(A_{p/}, C_{p/}).
$$
By functoriality this is a functor from $\Delta ^o$ to $n-1$-precats, and we
have a factorization
$$
A_{p/} \rightarrow P(A,C)_{p/} \rightarrow C_{p/}.
$$
The second morphisms in the factorization are equivalences, and the first
morphisms are cofibrations. The only problem is that $P(A,C)_{0/}$ is not a
set: it is an $n$-category which is equivalent to the set $C_0$.

For any $p$ there is a morphism $\psi _p:P(A,C)_{0/}\rightarrow P(A,C)_{p/}$
which, because it is a section of any one of the morphisms back to $0$, is a
cofibration and in fact even injective in the top degree.  If $p\rightarrow
q$ is
a morphism in $\Delta$ then  $P(A,C)_{q/} \rightarrow P(A,C)_{p/}$ composed with
$\psi _q$ is equal to $\psi _p$.  Hence if we set
$$
Q(A,C)_{p/}:= P(A,C)_{p/} \cup ^{P(A,C)_{0/}} C_0
$$
then $Q(A,C)_{p/}$ is functorial in $p\in \Delta$.  Now $Q(A,C)_{0/}= C_0$
is a set rather than an $n-1$-precat
so $Q(A,C)$ descends to a presheaf on $\Theta ^n$. We have a morphism
$A\rightarrow Q(A,C)$ projected from the above morphism into $P(A,C)$.
We also have a morphism $Q(A,C)\rightarrow C$ because the composed morphism
$$
P(A,C)_{0/} \rightarrow P(A,C)_{p/} \rightarrow C_{p/}
$$
factors through the unique morphism $C_0 \rightarrow C_{p/}$.  The composition
of these morphisms is the morphism $A\rightarrow C$.  We next claim that the
second morphism $Q(A,C)\rightarrow C$ is a weak equivalence. It suffices
by Lemma \ref{partialCat1} to prove that
$Q(A,C)_{p/}\rightarrow C_{p/}$ are weak equivalences, but the facts that
$P(A,C)_{0/}\rightarrow C_0$ is a weak equivalence and $P(A,C)_{0/}\rightarrow
P(A,C)_{p/}$ a cofibration imply (inductively using the closed model structure
for $n-1$-precats) that $P(A,C)_{p/}\rightarrow Q(A,C)_{p/}$ is a weak
equivalence. Now $P(A,C)_{p/} \rightarrow C_{p/}$ being a weak equivalence
implies that $Q(A,C)_{p/}\rightarrow C_{p/}$ is a weak equivalence.  This
proves the claim.

Finally we claim that $A\rightarrow Q(A,C)$ is a cofibration. It suffices to
prove that the $A_{p/}\rightarrow Q(A,C)_{p/}$ are cofibrations.  We know by
the inductive hypothesis that $A_{p/} \rightarrow P(A,C)_{p/}$ are cofibrations.
By the pushout definition of $Q(A,C)_{p/}$ and using the fact that
$P(A,C)_{0/}$ is a sub-presheaf of $P(A,C)_{p/}$, it suffices to
prove that the map
$$
A_{p/} \times _{P(A,C)_{p/}} P(A,C)_{0/} \rightarrow C_0
$$
is cofibrant. In fact we show below that for any $M$ of length $<n-1$,
$$
A_{p,M} \times _{P(A,C)_{p,M}} P(A,C)_{0,M} = A_0\subset A_{p,M}
$$
which implies what we want, since we have assumed that $A_0\rightarrow C_0$ is
injective. Note that the notation $P(A,C)_{0,M}$ means $(P(A,C)_{0/})_M$, and we
don't have in this case that this is a constant $n-1$-category so the usual rule
saying that $P(A,C)_{0,M}$ should be equal to $P(A,C)_{0,0}$ doesn't apply.

Fix any one of the maps $e:p\rightarrow 0 \rightarrow p$. This gives a map
$A_{p/} \rightarrow A_{p/}$ whose image is automatically $A_0$. This implies
that the fixed subsheaf of the endomorphism $e$ is equal to $A_0$. The
endomorphism acts compatibly on $P(A,C)_{p/}$ and the fixed point subsheaf there
is $P(A,C)_{0/}$.  For any $M$ of length $<n-1$ we have an inclusion $A_{p,M}
\hookrightarrow P(A,C)_{p,M}$.  This is compatible with the endomorphisms
$e$  on both sides, so the intersection of $A_{p,M}$ with the fixed point
set $P(A,C)_{0,M}\subset P(A,C)_{p,M}$ is the fixed point set $A_{0}\subset
A_{p,M}$. This shows the statement of the previous paragraph.

This completes the proof that $A\rightarrow Q(A,C)\rightarrow C$
is a factorization of the desired type, when $A_0 \rightarrow C_0$ is injective
on objects.  Note also that $Q(A,C)$ is functorial in the morphism
$A\rightarrow C$. Suppose now that $A\rightarrow B$ is any morphism. We put
$$
M(A, B):= Q(A, N(A,B)).
$$
We have the factorization $A\rightarrow N(A,B) \rightarrow B$ with the first
arrow injective on objects and the second arrow a weak equivalence, discussed
at the start of the proof.  The first arrow is then factored as $A\rightarrow
M(A,B) \rightarrow N(A,B)$ with the first arrow a cofibration and the second
arrow a weak equivalence. The factorization $A\rightarrow M(A,B)\rightarrow B$
therefore has the desired properties, and furthermore it is functorial in
$A\rightarrow B$
(this is needed in order to continue with the induction on $n$). This completes
the proof of CM5(2).

\eop

We refer to \cite{Quillen} for all of the consequences of Theorem \ref{cmc}.

Recall
also that a closed model category is said to be {\em proper}
if it satisfies the following two axioms:
\newline
{\bf Pr}(1)\,\, If $A\rightarrow B$ is a weak equivalence and $A\rightarrow
C$ a cofibration then $C\rightarrow B\cup ^AC$ is a weak equivalence;
\newline
{\bf Pr}(2)\,\, If $B\rightarrow A$ is a weak equivalence and $C\rightarrow
A$ a fibration then $B\times _AC\rightarrow C$ is a weak equivalence.

\begin{theorem}
\label{properness}
The closed model category $PC_n$ satisfies axiom {\bf Pr}(1); and it satisfies
axiom {\bf Pr}(2) for equivalences $B\rightarrow A$ between $n$-categories;
however it doesn't satisfy axiom {\bf Pr}(2) in general.
\end{theorem}
{\em Proof:}
We will prove stability of weak equivalences under coproducts.
Suppose $A\rightarrow B$ is a cofibration, and suppose $A\rightarrow C$ is a
weak equivalence. We would like to show that $B\rightarrow B\cup ^AC$
is a weak equivalence.  For this we use a version of the ``mapping cone''.
Recall that $\overline{I}$ is the category with two isomorphic objects $0,1$ and
no other morphisms.  The morphism
$B\times \{ 1\} \rightarrow B \times \overline{I}$ is a trivial cofibration,
so
$$
B\cup ^AC\rightarrow D:= (B\cup ^AC)\cup ^{B\times \{ 1\} } B\times \overline{I}
$$
is a trivial cofibration. It follows that the projection $D\rightarrow
B\cup ^AC$ deduced from $B\times \overline{I}\rightarrow B$ is a weak
equivalence.  Let
$$
E:= (B\times \{ 0\} )\cup ^{A\times \{ 0\} }A \times \overline{I}
$$
and note that $B\times \{ 0\} \rightarrow E$ is a weak equivalence
(since it is pushout by the trivial cofibration $A\times \{ 0\} \rightarrow
A\times \overline{I}$) hence $E\rightarrow B\times \overline{I}$ is a trivial
cofibration.  Thus the morphism
$$
E \cup ^{A\times \{ 1\} } C \rightarrow B\times \overline{I} \cup
^{A\times \{ 1\} } C
$$
is a trivial cofibration. But note that
$B\times \overline{I} \cup
^{A\times \{ 1\} } C =D$ so
$$
E \cup ^{A\times \{ 1\} } C \rightarrow D
$$
is a weak equivalence.  Finally,
$$
E \cup ^{A\times \{ 1\} } C = B\times \{ 0\} \cup ^{A\times \{ 0\}}
(A\times \overline{I} \cup ^{A\times \{ 1\} } C)
$$
and the morphism
$$
A\times \{ 0\}\rightarrow
A\times \overline{I} \cup ^{A\times \{ 1\} } C
$$
is a weak equivalence because it projects to  $A\rightarrow C$ which is by
hypothesis a weak equivalence. Therefore the map
$$
B\times \{ 0\} \rightarrow E \cup ^{A\times \{ 1\} } C
$$
is a weak equivalence, and from above $B\times \{ 0\} \rightarrow D$
is a weak equivalence. Following by the projection $D\rightarrow B\cup ^AC$
which we have seen to be a weak equivalence, gives the standard map
$B\rightarrow B\cup ^AC$ which is therefore a weak equivalence. This proves the
first half of properness.

We now prove the second statement, proceeding as usual by induction on $n$.
Factoring $B\rightarrow A$ into a cofibration followed by a fibration and
treating the fibration, we can assume that the morphism is a cofibration
(note that a fibration which is a weak equivalence, over an $n$-category, is
again an $n$-category so the hypotheses are preserved).
Let $A'\subset A$ be the full sub-$n$-category consisting of the objects
which are in the image of $B_0$.  Let $C':= A'\times _AC$. The morphism
$B\rightarrow A'$ is an iso-equivalence so the $B_{p/}\rightarrow A'_{p/}$ are
equivalences of $n-1$-categories. The morphisms $C'_{p/} \rightarrow
A'_{p/}$ are fibrant, so by the inductive hypothesis
$$
(C'\times _{A'} B)_{p/} = C'_{p/} \times _{A'_{p/}} B_{p/} \rightarrow C'_{p/}
$$
are equivalences. This implies that
$$
C\times _{A} B = C'\times _{A'} B \rightarrow C'
$$
is an equivalence.  Now $C'$ is a full sub-$n$-category of $C$ (meaning that
for any objects $x_0,\ldots , x_m$ of $C'$ the morphism
$C'_{m/}(x_0,\ldots , x_m)\rightarrow C_{m/}(x_0,\ldots , x_m)$ is an
isomorphism), so to prove that $C'\rightarrow C$ is an equivalence it suffices
to prove essential surjectivity.  Suppose $x$ is an object of $C$. It projects
to an object $y$ in $A$ which is equivalent to an object $y'$ coming from $A'$.
By \ref{intervalK} there is a morphism $K\rightarrow A$ sending $0$ to $y$ and
$1$ to $y'$. The object $x$  provides a lifting to $C$ over $\{ 0\}$, so by
the condition that $C\rightarrow A$ is fibrant there is a lifting to
$K\rightarrow C$ sending $0$ to $x$ and $1$ to an object lying over $y'$.
In particular $1$ goes to an object of $C'$. This shows that $x$ is equivalent
to an object of $C'$, the essential surjectivity we needed.

We now sketch an example showing why axiom {\bf Pr}(2) can't be true in
general. Let $A$ be the category $I^{(2)}$ with objects $0,1,2$ and one
morphism $i\rightarrow j$ for $i\leq j$ ($i,j=0,1,2$).  Let $B$ be the
sub-$1$-precat obtained by removing the morphism $0\rightarrow 2$ (it is the
pushout of two copies of $I$ over the object $1$).   The morphism $B\rightarrow
A$ is a weak equivalence. Let $C$ be a $1$-category with three objects $x_0,
x_1, x_2$ and morphisms from $x_i$ to $x_j$ only when $i\leq j$. There is
automatically a unique morphism $C\rightarrow A$ sending $x_i$ to $i$.
One can see that this morphism is fibrant. We can choose $C$ so that the
composition morphism
$$
C_1(x_0,x_1)\times C_1(x_1,x_2)\cong C_2(x_0, x_1, x_2)\rightarrow C_1(x_0,x_2)
$$
is not an isomorphism. Let $D:= Cat(B\times _AC)$. There is a unique morphism
$D\rightarrow C$ extending the second projection morphism, and this morphism
takes $D_1(x_0,x_1)$ (resp. $D_1(x_1,x_2)$) isomorphically to
$C_1(x_0,x_1)$ (resp. $C_1(x_1,x_2)$). However, the composition morphism for
$D$ is an isomorphism
$$
D_1(x_0,x_1)\times D_1(x_1,x_2)\stackrel{cong}{\rightarrow} D_1(x_0,x_2).
$$
Thus the morphism $D_1(x_0,x_2)\rightarrow C_1(x_0,x_2)$ is not an isomorphism;
thus $D\rightarrow C$ is not an equivalence and $B\times _AC\rightarrow C$ is
not a weak equivalence.
\eop

As one last comment in this section we note the following potentially useful
fact.

\begin{lemma}
\label{automatic}
If $f:A\rightarrow A'$ is a fibrant morphism of $m$-precats
then it is again fibrant when considered as a morphism of $n$-precats for $n\geq
m$.
\end{lemma}
{\em Proof:}
Suppose $m< n$.
Define the {\em brutal truncation} denoted $\beta \tau _{\leq m}$ from
$n$-precats to $m$-precats as follows. If $B$ is an $n$-precat then put
$$
\beta\tau _{\leq n-1}(B)_M := B_M
$$
for $M= (m_1, \ldots , m_k)$ with $k <m$ whereas for $M$ of length $m$ put
$$
\beta\tau _{\leq n-1}(B)_M:=B_M /\langle B_{M, 1} \rangle
$$
where $\langle B_{M, 1} \rangle$ denotes the equivalence relation on
$B_M$ generated by the image of $B_{M,1}\rightarrow B_M\times B_M$. This should
not be confused with the ``good'' truncation operation $T^{n-m}$ of
\cite{Tamsamani}, as in general they will not be the same (however they are
equal in the case of $n$-groupoids).
We claim that brutal truncation is compatible with the operation $BigCat$, that
is $$
BigCat (\beta \tau _{\leq m} B) = \beta \tau _{\leq m} (BigCat(B)).
$$
To prove this claim, we note two things:
\newline
(1) \,\, that if $\Sigma \rightarrow h$ is
one of our standard cofibrations for $n$-precats then
$\beta \tau _{\leq m}\Sigma \rightarrow \beta \tau _{\leq m}h$ is a standard
cofibration for $m$-precats; and
\newline
(2)\,\, that any map $\beta \tau _{\leq m}\Sigma \rightarrow
\beta \tau _{\leq m}B$ comes from a map $\Sigma \rightarrow B$ or---in the top
degree case---at least from a map $\Sigma \rightarrow BigCat(B)$.

By reordering
we find that the two sides of the above equation are the same, which gives the
claim.
Next
we claim that brutal truncation preserves weak equivalences. From the
previous claim it suffices to note that it preserves equivalences of
$n$-categories, and this follows from the fact that brutal truncation of
$1$-categories takes equivalences to isomorphisms of sets.

Finally, it is
immediate from the definitions that brutal truncation takes cofibrations of
$n$-precats to cofibrations of $m$-precats (using of course the fact that
there is no injectivity on the top degree morphisms for cofibrations of
$m$-precats).

Suppose $A$ is an  $m$-precat, and let
$Ind^n_m(A)$ denote $A$ considered as an $n$-precat (for this we simply set
$Ind^n_m(A)_{M,M'}:= A_M$ for $M$ of degree $m$ and any $M'$ or for $M$ of
degree $<m$ and empty $M'$).  Then (speaking of absolute $Hom$ here rather than
internal $Hom$ as in the next section)
$$
Hom (\beta \tau _{\leq m} B, A) = Hom (B, Ind^n_m(A)),
$$
in other words $\beta \tau _{\leq m}$ and $Ind^n_m$ are adjoint functors.

We can now prove the lemma. If $f$ is a fibrant morphism of  $m$-precats and
$B\rightarrow C$ is a trivial cofibration of $n$-precats then $\beta\tau _{\leq
m}B\rightarrow \beta\tau _{\leq
m}C$ is a trivial cofibration of $m$-precats, so $f$ has the lifting property
for this latter.  By adjointness $Ind^n_m(f)$ has the lifting property for
$B\rightarrow C$. Therefore $Ind ^n_m(f)$ is fibrant.
\eop

\numero{Internal $\underline{Hom}$ and $nCAT$}

Recall the result of Corollary \ref{forInternalHom}: that direct product with
any $n$-precat  preserves weak equivalences. Direct product also preserves
cofibrations, so it preserves trivial cofibrations. This property is not a
standard property of any closed model category, it is one of the nice things
about our present situation which allows us to obtain the right thing by
looking at internal $Hom$ of $n$-precats.

\begin{theorem}
\label{hom}
Suppose $A$ is an $n$-precat  and $B$ is a fibrant $n$-precat. Then the
internal  $
\underline{Hom}(A,B)$ of
presheaves over $\Theta ^n$ is a fibrant easy $n$-category. Furthermore if
$B'\rightarrow B$ is a fibrant morphism then
$\underline{Hom} (A, B')\rightarrow \underline{Hom} (A, B)$ is fibrant.
Similarly if $A\rightarrow A'$ is a cofibration and if
$B$ is fibrant then $\underline{Hom}(A', B)\rightarrow \underline{Hom}(A,B)$ is
fibrant. \end{theorem}
{\em Proof:}
Note that it suffices to prove that $\underline{Hom} (A,B)$ is fibrant, for
\ref{FibImpliesCat} then shows that it is an easy $n$-category.
A morphism $S\rightarrow \underline{Hom} (A,B)$ is the same thing as a
morphism $S\times
A\rightarrow B$.   Suppose $S\rightarrow T$ is a trivial cofibration. Then
$S\times A \rightarrow T \times A$ is  a trivial cofibration.
It follows immediately from the definition of $B$ being fibrant that any map
$S\times A \rightarrow B$ extends to a map $T \times A \rightarrow B$.
Thus $\underline{Hom} (A, B)$ is fibrant.
Similarly if $B'\rightarrow B$ is fibrant then any map $T\times A\rightarrow B$
with lifting $S\times A\rightarrow B'$ admits a compatible lifting
$T\times A\rightarrow B$.
Thus $\underline{Hom} (A,B')\rightarrow \underline{Hom}(A, B)$ is fibrant.

Suppose $A\rightarrow A'$ is cofibrant, and $B$ fibrant.
We show that  $\underline{Hom} (A',B)\rightarrow \underline{Hom}(A, B)$
satisfies the lifting property to be fibrant.  Say $S\rightarrow S'$ is a
trivial cofibration, and suppose we have maps $S'\rightarrow
\underline{Hom}(A, B)$ and lifting $S\rightarrow \underline{Hom}(A', B)$.
These are by definition maps $S'\times A \rightarrow B$ and $S\times
A'\rightarrow B$ which agree over $S\times A$.  These give a morphism
$$
f:S\times A' \cup ^{S\times A} S ' \times A \rightarrow B.
$$
The
morphism
$$
g:S\times A' \cup ^{S\times A} S ' \times A \rightarrow S' \times A'
$$
is a cofibration.  Lemma \ref{pushout}
applied to the trivial cofibration $S\times A\rightarrow S'\times A$ implies
that the  morphism  $$
S\times A'\rightarrow S\times A' \cup ^{S\times A} S '\times A
$$
is a weak equivalence.  On the other hand the morphism
$S\times A' \rightarrow S'\times A'$ is a weak equivalence by
\ref{forInternalHom}, so by Lemma \ref{remark} the morphism  $g$
is a weak equivalence. Thus the fact that $B$ is fibrant means that our
morphism $f$ extends to a morphism $S'\times A' \rightarrow B$, and this gives
exactly the desired lifting property for the last statement of the theorem.
\eop

\begin{lemma}
\label{stability}
Suppose $A\rightarrow A'$ is a weak equivalence,  and  $B$
fibrant. Then $\underline{Hom} (A', B)\rightarrow \underline{Hom} (A, B)$ is an
equivalence of $n$-categories.

If $B\rightarrow B'$ is an equivalence
of fibrant $n$-precats then $\underline{Hom}(A,B)\rightarrow
\underline{Hom}(A,B')$ is an equivalence.

Suppose $A\rightarrow B$ and $A\rightarrow C$ are cofibrations. Then
$$
\underline{Hom} (B\cup ^AC, D) = \underline{Hom} (B, D)
\times _{\underline{Hom}(A,D)}\underline{Hom}(C,D).
$$
\end{lemma}
{\em Proof:}
The last statement is of course immediate, because for any $S$ we  have
$(B\cup ^AC)\times S = (B\times S ) \cup ^{(A\times S)} (C\times S)$.
We treat the other statements.

Suppose first that $A\rightarrow A'$ is a trivial
cofibration.  Suppose that $S\rightarrow T$ is any
cofibration.  Suppose we have maps
$T\rightarrow \underline{Hom}(A, B)$ lifting over $S$  to $S\rightarrow
\underline{Hom}(A', B)$. We claim that the lifting extends to $T$;
then the
characterization of weak equivalences in (\cite{Quillen} \S 5, Definition 1,
Property M6, part (c)) will imply that $\underline{Hom}(A', B)\rightarrow
\underline{Hom}(A,B)$ is a weak equivalence. To prove the claim,
note that our data correspond to a morphism
$$
T\times A \cup ^{S\times A} S\times A' \rightarrow B.
$$
The morphism
$S\times A \rightarrow S\times A'$ is a trivial cofibration, so the morphism
$$
T\times A \rightarrow T\times A \cup ^{S\times A} S\times A'
$$
is a trivial cofibration, and since $T\times A \rightarrow T\times A'$ is
a weak equivalence we get that the morphism
$$
T\times A \cup ^{S\times A} S\times A'
\rightarrow T\times A'
$$
is a trivial cofibration. The fibrant property of $B$ implies that our map
extends to a map $T\times A' \rightarrow B$, so we get the required lifting
to $T\rightarrow \underline{Hom}(A', B)$. This implies that
$\underline{Hom}(A',B)\rightarrow \underline{Hom}(A,B)$ is a fibrant
weak equivalence.

Next we treat the case of any weak equivalence $A\rightarrow A'$.
Let $C$ be the $n$-precat pushout of $A\times 0\rightarrow A\times \overline{I}$
and $A\times 0 \rightarrow A' \times 0$.  Since $0\rightarrow \overline{I}$ is
a trivial cofibration, the various morphisms
$$
A\hookrightarrow C, \;\; A' \hookrightarrow C, \;\; C \rightarrow A'
$$
(the first sending $A$ to $A\times 1$, the second sending $A'$ to $A' \times 0$
and the third coming from the projection $A\times \overline{I} \rightarrow A'$)
are all weak equivalences.
We have a composable pair of morphisms
$$
\underline{Hom} (A', B) \rightarrow \underline{Hom} (C, B) \rightarrow
\underline{Hom} (A', B)
$$
composing to the identity, and where the second arrow is an equivalence
by the previous paragraph
since $A'\rightarrow C$ is a trivial cofibration. Therefore the first arrow is
an equivalence.  Next, the morphism $\underline{Hom} (C, B)\rightarrow
\underline{Hom} (A,B)$ obtained from the trivial cofibration $A\rightarrow C$
(going to $A\times 1$) is an equivalence, so the composed map $\underline{Hom}
(A', B)\rightarrow \underline{Hom} (A,B)$ is an equivalence. This is the map
induced by our original $A\rightarrow A'$. This completes the proof of the first
part of the lemma.

We now turn to the second part and treat first a fibrant weak equivalence
$B\rightarrow B'$.  Note first that such a morphism satisfies the lifting
property for any cofibrations (this is the other half of CM4 which comes from
Joyal's trick).
We prove that $\underline{Hom}(A,B)\rightarrow \underline{Hom}(A,B')$
satisfies lifting for any cofibration (which as above implies that it is a
fibrant weak equivalence).  Suppose $S\rightarrow T$ is a cofibration.
and $T\rightarrow \underline{Hom} (A,B')$ is
a map with lifting over $S$ to a map $S\rightarrow \underline{Hom}
(A,B)$. These correspond to maps $T\times A \rightarrow B'$ and lifting to
$S\times A \rightarrow B$. The morphism $S\times A \rightarrow
T\times A $ is a cofibration so by the lifting property of $B\rightarrow B'$
for
any cofibration, there is a lifting to $T\times A\rightarrow B$  compatible with
the given map on $S$.  This establishes the necessary lifting property to
conclude that $\underline{Hom}(A,B)\rightarrow \underline{Hom}(A,B')$ is a
fibrant equivalence.

Next suppose that $i:B\rightarrow B'$ is a trivial cofibration of fibrant
$n$-precats. The lifting property for $B$ lets us choose a retraction $r:
B'\rightarrow B$ such that $ri= 1_B$.
Let $p:B' \cup ^BB' \rightarrow B'\rightarrow
B'$ be the  projection which induces the identity on both of the
components $B'$. Note that $B'\rightarrow B' \cup _BB'$ is a trivial
cofibration by \ref{pushout} so the projection $p$ is a weak equivalence
(using \ref{remark}).  Choose a
factorization
$$
B' \cup ^BB' \rightarrow P \rightarrow B'
$$
with the first morphism cofibrant and the  second morphism fibrant (whence $P$
fibrant itself); and both morphisms weak equivalences. Let $q: B' \cup
^BB'\rightarrow B'$ be the morphism inducing the retraction $r$ on the first
copy of $B'$ and the identity on the second copy. Since $B'$ is fibrant
this extends to a morphism we again denote $q: P\rightarrow B'$.
The result of the previous paragraph implies that the morphism
$$
p:\underline{Hom}(A, P)\rightarrow \underline{Hom}(A, B')
$$
is an equivalence, which implies that either of the two morphisms
$$
j_0, j_1:\underline{Hom}(A, B')\rightarrow \underline{Hom}(A, P)
$$
(refering to the two inclusions $j_0,j_1: B' \rightarrow P$)
are equivalences. Now we have that the composition
$$
\underline{Hom}(A, B')\stackrel{j_1}{\rightarrow}  \underline{Hom}(A, P)
\stackrel{q}{\rightarrow} \underline{Hom}(A, B')
$$
is the morphism induced by $qj_1= 1_{B'}$ thus it is the identity. The fact
that the morphism induced by $j_1$ (the first of the above pair) is an
equivalence implies that the morphism induced by $q$ (the second in the above
sequence) is an equivalence. But since the morphism induced by $j_0$ is
an equivalence, we get that the morphism induced by $qj_0=ir$
is an auto-equivalence of $\underline{Hom}(A,B')$.  The morphism induced by
$ri=1_B$ is of course the identity. The last part of Lemma \ref{remark}
now implies that the morphism
$$
i: \underline{Hom}(A,B)\rightarrow \underline{Hom}(A,B')
$$
is an equivalence.
This proves the statement in case of a trivial
cofibration.

Finally note that any equivalence of $n$-categories  $B\rightarrow B'$
decomposes as a composition $B\rightarrow C \rightarrow B'$ where the first
arrow is a trivial cofibration and the second a fibration and weak equivalence.
Note that $C$ is fibrant since by hypothesis $B'$ is fibrant. Thus our two
previous discussions apply to give that the two morphisms
$$
\underline{Hom}(A,B)\rightarrow \underline{Hom}(A,C)
\rightarrow \underline{Hom}(A,B')
$$
are equivalences, their composition is therefore an equivalence. This completes
the proof of the first paragraph of  Lemma \ref{stability}.
\eop

For any
fibrant $n$-categories $A$, $B$  and $C$  we have
composition morphisms
$$
\underline{Hom} (A, B) \times \underline{Hom}(B,C) \rightarrow
\underline{Hom}(A,C),
$$
which are associative.

Define a simplicial $n$-category $nCAT$ by setting
$nCAT_0 $ equal to a set of representatives for {\em isomorphism} classes of
fibrant $n$-categories, and by setting
$$
nCAT_m(A_0, \ldots , A_m):= \underline{Hom} (A_0, A_1) \times \ldots \times
\underline{Hom} (A_{m-1},
A_m),
$$
with simplicial structure given by the above compositions.  Since
$nCAT_0$ is a set considered as $n$-precat, this simplicial $n$-precat
(presheaf on $\Delta \times \Theta ^{n}$) descends
to a presheaf on $\Theta ^{n+1}$, in other words it is an $n+1$-precat.
The composition gives the necessary conditions in the first degree, and in
higher degrees the fact that $\underline{Hom} (A,B)$ are $n$-categories
completes what we
need to know to conclude that $nCAT$ is an  $n+1$-category.

If $A$ and $B$ are $n$-categories but not necessarily fibrant then let $Fib(A)$
and $Fib(B)$ be their fibrant replacements (given by the above construction for
Theorem \ref{cmc} for example). We call these the {\em fibrant envelopes}. The
``right'' $n$-category of morphisms from $A$ to $B$ is $\underline{Hom}(Fib(A),
Fib(B))$. We will sometimes use the notation
$$
HOM (A, B):= \underline{Hom}(Fib(A), Fib(B)).
$$
We obtain an $n+1$-category equivalent to $nCAT$ by taking all $n$-categories
as objects and taking the $HOM (A,B)$ as morphism $n$-categories.

{\em Question:} Describe the fibrant envelope of the $n+1$-category $nCAT$.
This would be important if one wants to consider weak morphisms $A\rightarrow
nCAT$ as families of $n$-categories indexed by $A$ in a meaningul way.

We have almost proved Conjecture (\cite{Tamsamani} between 1.3.6 and 1.3.7)
on the existence of $nCAT$. We
just have to check that the truncation of $nCAT$ down to a $1$-category is
equivalent to the localization of the category of $n$-categories by the
subcategory of morphisms which are  equivalences.

In Corollary \ref{htytype} above, we have seen that the localization in
question is equal to the localization of $PC_n$ by the weak equivalences.
We now know that $PC_n$ is a closed model category, and Quillen shows in this
case that the $Hom$ in the localized category is equal to the set of {\em
homotopy classes} of morphisms between fibrant and cofibrant objects. In our
(second) definition above, we took $nCAT$ to be the category of fibrant (and
automatically cofibrant) $n$-precats.  The $Hom$ in the localized category is
thus the set of homotopy classes of maps. On the other hand, the truncation of
$nCAT$ down to a $1$-category is obtained by replacing the
$\underline{Hom}(A,B)$
$n$-categories by their sets of equivalence classes of objects. Thus, to prove
the conjecture we simply need to show that for $A$ and $B$ fibrant, the set of
equivalence classes of objects in the $n-1$-category $\underline{Hom} (A,B)$
is equal to the set of homotopy classes of maps from $A$ to $B$.  Note that the
objects of $\underline{Hom} (A,B)$ are again just the maps from $A$ to $B$
so we are
reduced to showing the following lemma.

\begin{lemma}
If $A$ and $B$ are fibrant $n$-precats then two morphisms $f,g: A\rightarrow B$
are homotopic in the sense of \cite{Quillen} if and only if the corresponding
elements of the $n$-category $\underline{Hom}(A,B)$ are equivalent.
\end{lemma}
{\em Proof:}
Suppose $f$ and $g$ are homotopic (\cite{Quillen} p. 0.2).  Then there is an
object $A'$ with  morphisms $i,j: A\rightarrow A'$ each inducing a weak
equivalence, with a projection $p:A'\rightarrow A$ such that the compositions
$pi$ and $pj$ are the identity, and a morphism $h: A'\rightarrow B$ such that
$hi=f$ and $hj=g$. We may assume that $A'$ is fibrant.  Then we obtain pullback
morphisms on the $\underline{Hom}$ $n$-categories and in particular, two
morphisms
$$
i^{\ast}, j^{\ast} :\underline{Hom} (A', B) \rightarrow \underline{Hom} (A, B)
$$
and a morphism
$$
p^{\ast}: \underline{Hom} (A, B)\rightarrow \underline{Hom} (A', B)
$$
which are weak equivalences by Lemma \ref{stability}.  These induce isomorphisms
on the sets of equivalence classes which we denote $T^n\underline{Hom}
(A,B)$ etc., so we
have $$
T^ni^{\ast}, T^nj^{\ast} :T^n\underline{Hom} (A', B) \cong T^n
\underline{Hom} (A, B)
$$
and
$$
T^np^{\ast}:T^n \underline{Hom} (A, B)\cong T^n  \underline{Hom} (A', B).
$$
Here, as before, we have that $T^ni^{\ast} \circ T^np^{\ast}$ and
$T^n j^{\ast}\circ T^np^{\ast}$ are equal to the identity. This implies that
$T^ni^{\ast} = T^nj^{\ast}$ and hence that their applications to the class of
$h$ give the same equivalence class. The results are respectively the classes
of $f$ and of $g$, hence $f$ is equivalent to $g$.

Conversely suppose $f$ and $g$ are equivalent as objects in
$\underline{Hom}(A,B)$. Then
by Proposition \ref{intervalK}
there is a contractible $K$ with $0,1\in K$ and a morphism
$K\rightarrow \underline{Hom}(A,B)$ taking $0$ to $f$ and $1$ to $g$. This
yields
(by the universal property of the internal $\underline{Hom}$) a morphism
$h:A\times
K\rightarrow B$. This morphism together with the various others gives
a homotopy from $f$ to $g$. \eop

\begin{corollary}
{\rm (\cite{Tamsamani} Conjecture 1.3.6-7)}
The $n+1$-category $nCAT$ yields when truncated down to a $1$-category the
localization $Ho-n-Cat$ of \cite{Tamsamani}.
\end{corollary}
\eop

\numero{$n$-stacks}

We can give a preliminary discussion of the notion of $n$-stack, following the
lines that are already well known for simplicial presheaves and even $n$-stacks
of $n$-groupoids (approached via topological spaces in \cite{flexible},
discussed for $n=2,3$ in \cite{Breen}). Our present discussion will be
incomplete, basically for the following reason: if $\Xx$ is a $1$-category,
there are several natural types of objects which represent the idea of a family
of $n$-categories indexed (contravariantly) by $\Xx$, and we would like to know
that all of these notions are equivalent.  The main possible versions are:
\newline 1. \,\, A functor $\Xx ^o\rightarrow nCAT$, which if we take the second
point of view on $nCAT$ presented above, is the same thing as a presheaf of
fibrant $n$-categories over $\Xx$;
\newline
2.\,\, A weak functor from $\Xx ^o$ to $nCAT$, in other words a functor from
$\Xx ^o$ to $Fib(nCAT)$ or (what is basically the same thing) an element of
$HOM(\Xx ^o, nCAT)$ i.e. a morphism in $(n+1)CAT$;
\newline
3.\,\, A ``fibered $n$-category over $\Xx$'', which would be a morphism of
$n$-categories $\Ff \rightarrow \Xx$ (note that a $1$-category considered as an
$n$-category is automatically fibrant by \ref{automatic}) satisfying some
condition analogous to the definition of fibered $1$-category---I haven't
written down this condition (note however that it is distinct from the
condition that the morphism be fibrant in the sense we use in this paper).

Here is what I currently know about the relationship between these points of
view. From (1) one automatically gets (2) just by composing with the morphism
$nCAT \rightarrow Fib(nCAT)$ to the fibrant envelope.  From (2) one should be
able to get (3)  by pulling back a universal fibered $n+1$-category over
$Fib(nCAT)$. To construct this universal object, first construct a universal
$n+1$-category $\Uu \rightarrow nCAT$ (with fibers the $n$-categories being
parametrized---in particular this morphism is relatively $n$-truncated) then
replace the composed morphism $\Uu \rightarrow Fib(nCAT)$ by a fibrant
morphism.  Finally, from (3) one should be able to get (1) by applying the
``sections functor'': if $\Ff \rightarrow \Xx$ is a fibered $n$-category then
define $\Gamma (\Xx , \Ff )$ to be the $n$-category fiber (calculated in the
correct homotopic sense) over $1_{\Xx} \in HOM (\Xx , \Xx)$ of
$$
HOM (\Xx , \Ff )\rightarrow HOM (\Xx , \Xx ).
$$
Require now that $\Ff \rightarrow \Xx$ be a fibrant morphism (if this doesn't
come into the condition of being fibered already).  Then
$$
X\in \Xx\mapsto \Gamma (\Xx /X, \Ff \times _{\Xx} (\Xx /X ))
$$
should be strictly functorial in the variable $X$ yielding a presheaf $\Xx ^o
\rightarrow nCAT$ which is notion (1). The condition of being a fibered
category should imply (as it does in the case $n=1$) that the morphism
$$
\Gamma (\Xx /X, \Ff \times _{\Xx} (\Xx /X ))\rightarrow \Ff _X:= \Ff \times
_{\Xx} \{ X\}
$$
be an equivalence of $n$-categories (one might even try to take this condition
as the definition of being fibered but I'm not sure if that would work).
Finally we would like to show that doing these three constructions in a circle
results in an essentially equivalent object.

The previous paragraph is for the moment speculative, the main questions left
open  being the definition of ``fibered $n$-category'' and the construction of
the universal family.  However, for the rest of this section we will discuss
the theory of $n$-stacks supposing that the above equivalences are known.
Denote by $\int$ the operation going from (1) to (3).

Suppose $\Xx$ is  a site.
There are a couple of different ways of approaching the notion of $n$-stack
over $\Xx$. Our first definition will be modelled on what was done in
\cite{flexible}.
A fibered $n$-category $\Ff \rightarrow \Xx$ is an {\em $n$-stack} if for any
$X\in \Xx /X$ and any sieve $\Bb \subset \Xx /X$ the morphism
$$
\Gamma (\Xx /X, \Ff \times _{\Xx} (\Xx /X ))\rightarrow
\Gamma (\Bb , \Ff \times _{\Xx} \Bb )
$$
is an equivalence of $n$-categories.
If $\Ff \rightarrow \Xx$ is a fibered $n$-category then we (should be able to)
construct the {\em associated $n$-stack} by iterating $n+2$ times the operation
$$
L(\Ff ):= \int \left( X\mapsto \lim _{\rightarrow , \Bb \subset \Xx /X}
\Gamma (\Bb , \Ff \times _{\Xx} \Bb ) \right) .
$$
This conjecture is based on the corresponding result for flexible sheaves in
\cite{flexible}.

The second main type of approach is to combine the
theory of simplicial presheaves of Jardine-Joyal-Brown-Gersten (cf
\cite{Jardine}) with the discussion in the present paper to obtain a closed
model structure for the category of presheaves of $n$-precats over $\Xx$. In
this case the fibrant condition would imply the condition of being an
$n$-stack. To give the definitions (without proving that we get a closed model
category) it suffices to define weak equivalence---the cofibrations being just
the maps which over each object of $\Xx$ are cofibrations of $n$-precats, and
the fibrations then being defined by the lifting property for trivial
cofibrations. (As usual the main problem would then be to prove that pushout by
a trivial cofibration is again a trivial cofibration---and for this we could
probably just combine the proofs of Jardine/Joyal \cite{Jardine} and the present
paper.) If $A\rightarrow B$ is a morphism of presheaves of $n$-precats over
$\Xx$ then we obtain a morphism $Cat(A)\rightarrow Cat(B)$ of presheaves of easy
$n$-categories (where $Cat(A)(X):= Cat(A(X))$). We will say that $A\rightarrow
B$ is a weak equivalence if $Cat(A)\rightarrow Cat(B)$ is a weak equivalence of
presheaves of $n$-categories, notion which we now define. Let $T$ denote
Tamsamani's truncation operation \cite{Tamsamani} which is functorial so it
extends to presheaves of $n$-categories.  A morphism of presheaves of
$n$-categories $A\rightarrow B$ is a {\em weak equivalence in top degree} if for
every $n$-morphism of $B$ and lifting of its source and target to
$n-1$-morphisms in $A$, there exists a unique lifting to an $n$-morphism in $A$.
Now we say that a morphism $A\rightarrow B$ of presheaves of $n$-categories
is a {\em weak equivalence} if for every $k$ the morphism
$$
Sh(T^kA)\rightarrow Sh(T^kB)
$$
is a weak equivalence in top degree, where $Sh$ denotes the stupid
sheafification operation (i.e. sheafify each of the presheaves $A_M$).

In this point of view, if $A$ is a presheaf of $n$-categories over $\Xx$ then
we define the {\em associated $n$-stack} to be the fibrant object
equivalent to $A$ in the
previous presumed closed model category.

\subnumero{$n$-categories as $n-1$-stacks}
Heuristically we can define a structure of {\em
site} on $\Delta$ where the coverings of an object $m$ are the collections of
morphisms $\lambda _i \hookrightarrow m$ where
$\lambda _i = \{ a_i , a_i + 1, \ldots , b_i\}$ such that $a_{i+1}=b_i$.
If $A$ is an $n$-precat then the collection $\{ A_{p/}\}$ may be thought of as
a presheaf of $n-1$-precats over $\Delta $. The condition to be an $n$-category
is that this should be a presheaf of $n-1$-categories which should satisfy
descent for coverings, i.e. it should be an $n-1$-stack of $n-1$-categories over
this site.   The construction $Cat$ is essentially just finding the  $n-1$-stack
associated to an $n-1$-prestack by an operation similar to that described in
\cite{flexible}. The main problems above are caused by the fact
that this site doesn't admit fiber products. It might be a good idea to replace
this site by  its associated topos, the category of categories, which would lead
to the yoga: {\em that an $n$-category is an $n-1$-stack over the topos of
categories.}

It might be possible, by treating
$n$-stacks at the same time as $n$-categories, to simplify the arguments of the
present paper by recursively defining $n$-categories as $n-1$-stacks. I haven't
thought about this any further.

\numero{The generalized Seifert-Van Kampen theorem}

Our closed model category structure allows us (with a tiny bit of extra work)
to obtain the analogue of the Siefert-Van Kampen theorem for the Poincar\'e
$n$-groupoid of a topological space $\Pi _n(X)$ defined by Tamsamani
(\cite{Tamsamani}, \S 2.3).

\begin{theorem}
\label{svk}
If $X$ is a space covered by open subsets $X = U\cup V$ then (setting $W:=
U\cap V$) $\Pi _n (X)$ is equivalent to the category-theoretic pushout of the
diagram
$$
\Pi _n (U) \leftarrow \Pi _n (W) \rightarrow \Pi _n (V).
$$
\end{theorem}

In order to prove this theorem we recall  Tamsamani's
realization functor from $n$-precats to topological spaces (\cite{Tamsamani}
\S 2.5).  There is a covariant functor $R: \Delta ^n \rightarrow Top$ which
associates to $M= (m_1, \ldots , m_n)$ the product $R^{m_1}\times
\ldots \times R^{m_n}$ where $R^m$ denotes the usual topological
$m$-simplex.   If $A: (\Delta ^n)^o\rightarrow Sets$ is a
presheaf of sets then Tamsamani defines the {\em realization of $A$} in the
standard way combining $R$ and $A$. We denote this $\langle R, A\rangle$ because
it is a sort of pairing of functors.  If $A$ is an $n$-precat in our notations
then pull it back to a presheaf on $\Delta ^n$ and apply the realization (we
still denote this as  $\langle R, A\rangle$). The functor $\langle R, \cdot
\rangle$ obviously  preserves pushouts.

{\em Caution:} The realization functor does not preserve cofibrations.
It takes injective morphisms of presheaves over $\Theta^n$ to cofibrations of
spaces, but the cofibrations which are not injective in the top degree
are taken to non-injective morphisms.

Recall Tamsamani's Proposition 3.4.2(ii):
\begin{proposition}
If $A\rightarrow B$ is an equivalence of $n$-categories then
$$
\langle R, A \rangle \rightarrow \langle R, B \rangle
$$
is an equivalence of spaces.
\end{proposition}
{\em Proof:}
The proof of \cite{Tamsamani} for $n$-groupoids using induction on $n$ also
works
for $n$-categories.
\eop

Say that a morphism $X\rightarrow Y$ of topological spaces is an {\em $n$-weak
equivalence}  if it is an isomorphism on $\pi _0$ and for any choice of
basepoint in $X$ it is an isomorphism on $\pi _i$ for $i\leq n$. This is
equivalent to saying that it induces a weak equivalence on the Postnikov tower
up to stage $n$.

\begin{corollary}
\label{equivtoequiv}
The realization functor $\langle R, \cdot \rangle$ from  $n$-precats to $Top$
takes weak equivalences to $n$-weak equivalences of topological spaces.
\end{corollary}
{\em Proof:}
Realization takes our standard trivial cofibrations $\Sigma \rightarrow h$
to homotopy equivalences of topological spaces.
This is essentially the content of the constructions of retractions in the
proof of Theorem 2.3.5 (that $\Pi _n(X)$ is an $n$-nerve) of \cite{Tamsamani}.

For all except the upper boundary cases, the standard trivial cofibrations are
taken to cofibrations of topological spaces. Pushout by the injective
standard trivial cofibrations becomes pushout by a trivial cofibration of
spaces, whence a homotopy equivalence.
In order to deal with the upper boundary cases we introduce the following
notation:
$$
\rangle R, A \langle _{q}
$$
denotes the $q$-skeleton of the realization of $A$, that is the realization
taken over all $M$ with $\sum m_i \leq q$.  Then, if $q\leq n-1$
the functor $\rangle R, \cdot \langle _{q}$ takes cofibrations to cofibrations
of spaces.

Suppose $\varphi : \Sigma \rightarrow h$ is a standard trivial cofibration
in one
of the boundary cases. Using the notation of \S 2 we can write $\Sigma$ as the
coequalizer of
$$
h' \sqcup  h' \sqcup  \Upsilon ' \rightarrow h^a \sqcup  h^a
$$
(the component $\Upsilon$ which appears on the right in the general case
disappears in the upper boundary case).  The map $\Sigma \rightarrow h$
is given by the map $h^a \sqcup  h^b\rightarrow h$
which in this case is two times the identity (because $h(M, m, 1^{k+1})$ which
doesn't exist is replaced by $h:= h(M, m, 1^k)= h^a=h^b$). The cells in
$\langle R,\Sigma \rangle $ of dimension $n$ are automatically of the form
$h(1^l, 1, 1^k)$ for maps $1^l\rightarrow M$ and $1\rightarrow m$. There are two
such which are identified whenever $1\rightarrow m$ is not one of the principal
morphisms. The cells coming from principal $1\rightarrow m$ occur only once in
the realization of $\Sigma$ already. It follows (since any non-principal
$1\rightarrow m$ is a path which is homotopic to a concatenation of principal
$1\rightarrow m$) that the $n$-cells which are identified are homotopic.

Note that on the level of cells of dimension $<n$ the morphism $\Sigma
\rightarrow h$ is an isomorphism. In particular, pushout via $\varphi$ over any
$\Sigma \rightarrow A$  preserves the $n-1$-skeleton of the realization, and the
$n$-cells which are identified are homotopic.  In this boundary case the pushout
by $\varphi$ is surjective, in particular it is surjective on $n+1$-cells.
A surjective morphism of cell complexes which is an isomorphism on
$n-1$-skeleta and which only identifies $n$-cells which are homotopic (relative
the $n-1$-skeleton) is an $n$-weak equivalence. This completes the proof that
pushout by any of our standard trivial cofibrations $\varphi$ induces an
$n$-weak equivalence.

It follows by construction of the operation $Cat$ that for any $n$-precat
$A$ the
morphism
$$
\langle R, A \rangle \rightarrow \langle R, Cat(A) \rangle
$$
is an $n$-weak equivalence of spaces.
Now we can complete the proof: if $A\rightarrow B$ is a weak equivalence then
by definition $Cat(A)\rightarrow Cat(B)$ is an equivalence of $n$-categories so
in the diagram
$$
\begin{array}{ccc}
\langle R, A \rangle  & \rightarrow & \langle R, B \rangle  \\
\downarrow && \downarrow \\
\langle R, Cat(A) \rangle  & \rightarrow & \langle R, Cat(B) \rangle
\end{array}
$$
the vertical arrows are $n$-weak equivalences from the previous argument and the
bottom arrow is a weak equivalence by the proposition, so the top arrow is  an
$n$-weak equivalence of spaces.
\eop

\begin{lemma}
\label{groupoid}
If $A\rightarrow B$ and $A\rightarrow C$ are morphisms of $n$-groupoids with
one being a cofibration,
then the category-theoretic pushout $Cat(B\cup ^AC)$ is
an $n$-groupoid.
\end{lemma}
{\em Proof:}
We say that an $n$-precat is a {\em pre-groupoid} if its associated
$n$-category is a groupoid.
We prove that the pushout of pre-groupoids is again a pre-groupoid, and we
proceed by induction on $n$ so we may assume this is known for
$n-1$-pre-groupoids.

Suppose now that $A$, $B$ and $C$ are $n$-groupoids with morphisms as in the
statement of the lemma.   Then the $A_{p/}$, $B_{p/}$ and $C_{p/}$ are
$n-1$-groupoids, and
$$
(B\cup ^AC)_{p/} = B_{p/}\cup ^{A_{p/}}C_{p/}.
$$
In particular, $(B\cup ^AC)_{p/}$ are $n-1$-pre-groupoids. The process of going
from this collection of $n-1$-precats to the collection corresponding to $Cat
(B\cup ^AC)$ as described in \S 4, uses only iterated pushouts by the various
$(B\cup ^AC)_{p/}$ in various combinations. Since we know by induction that
pushouts of $n-1$-pre-groupoids are again $n-1$-pre-groupoids, it follows that
$Cat(B\cup ^AC)_{p/}$ are $n-1$-groupoids.  It suffices now to show that the
truncation of $Cat(B\cup ^AC)$ down to a $1$-category is a groupoid. But this
truncation is the same as the brutal truncation since we know that the
$Cat(B\cup ^AC)_{p/}$ are $n-1$-groupoids. On the other hand, brutal
truncation commutes with the operations $Cat$ and pushout, therefore
the truncation of  $Cat(B\cup ^AC)$ is the pushout of the truncations of $A$,
$B$ and $C$ which are groupoids. Finally, the $1$-category pushout of groupoids
is again a groupoid, so $Cat(B\cup ^AC)$ is a groupoid.

To complete the proof it remains to be seen that the pushout of
$n$-pregroupoids is a pregroupoid. Suppose $A$, $B$ and $C$ are
$n$-pre-groupoids. Then by reordering
$$
BigCat(B\cup ^AC) = BigCat(Cat(B) \cup ^{Cat(A)}Cat(C)).
$$
Thus $BigCat(B\cup ^AC)$ is the category-theoretic pushout of
$n$-groupoids so by the previous argument it is an $n$-groupoid. This shows
that $B\cup ^AC$ is an $n$-pre-groupoid, completing the proof of the
induction step.
\eop

{\em Proof of Theorem \ref{svk}:}
Note first of all that Tamsamani's proof that  $\Pi_n(X)$ is an $n$-category
(\cite{Tamsamani} Theorem 2.3.6) actually shows that it is an easy
$n$-category.
With the notations of the theorem, we have a diagram of easy $n$-categories
$$
\begin{array}{ccc}
\Pi _n(W) & \rightarrow & \Pi _n(U) \\
\downarrow && \downarrow \\
\Pi _n(V) & \rightarrow & \Pi _n(X).
\end{array}
$$
Let $A$ be the pushout $n$-precat of the upper and left arrows. We have a
morphism $A\rightarrow \Pi _n(X)$, and hence (non-uniquely) $Cat(A)\rightarrow
\Pi _n(X)$ since the latter is an easy $n$-category. The realization of $A$ is
the pushout of the realizations of    $\Pi _n(U)$ and $\Pi _n(V)$ over $\Pi
_n(W)$. These last realizations are $n$-weak equivalent to $U$, $V$ and $W$
respectively (\cite{Tamsamani} 3.3.4), so the realization of $A$ is $n$-weak
equivalent to the pushout of $U$ and $V$ over $W$, in other words to $X$. Thus
the morphism $A\rightarrow \Pi _n(X)$ induces an $n$-weak equivalence on
realizations.  On the other hand we have seen above that $A\rightarrow Cat(A)$
induces an $n$-weak equivalence on realizations. Thus the morphism
$Cat(A)\rightarrow \Pi _n(X)$ induces an $n$-weak equivalence on realizations.
Lemma \ref{groupoid} implies that $Cat(A)$ is an $n$-groupoid.
Applying the functor $\Pi _n$ again and using Proposition 3.4.4 of
\cite{Tamsamani} we conclude that $Cat(A)\rightarrow \Pi _n(X)$ is an
equivalence of $n$-groupoids. This proves the theorem.
\eop

\numero{Nonabelian cohomology}

 If $A$ is a fibrant $n$-category and
$X$ a topological space then define the {\em nonabelian cohomology of $X$
with coefficients in $A$} to be $H(X, A): = \underline{Hom} (\Pi _n(X),
A)$. It is an
$n$-category.  This satisfies Mayer-Vietoris: if $U,V\subset X$ and $W=U\cap V$
then
$$
m:H(X, A) \rightarrow H(U,A)\times _{H(W,A)}H(V,A)
$$
is an equivalence of $n$-categories (where the fiber product is understood to
be the homotopic fiber product obtained by replacing one of the morphisms with
a fibrant morphism). To see this, note that if $A$ is fibrant then for any
cofibration $B\rightarrow C$ the morphism $\underline{Hom}(C,A)\rightarrow
\underline{Hom}(B,A)$ is
fibrant. To prove this claim it suffices to remark that if $S\rightarrow T$ is a
trivial cofibration then  $$ S\times C \cup ^{S\times B} T\times B \rightarrow T
\times C $$
is a trivial cofibration, now apply the universal property of the internal
$\underline{Hom}$ to obtain the lifting property in question.

In particular, note that $H(U,A)\rightarrow H(W,A)$ and $H(V,A)\rightarrow
H(W,A)$ are fibrations (since open inclusions induce cofibrations of $\Pi _n$).
Thus the above fiber product is the homotopic fiber product.

We now prove that the Mayer-Vietoris map $m$ is an equivalence of
$n$-categories. It is the same as the map
$$
\underline{Hom} (\Pi _n(X), A)\rightarrow \underline{Hom} (\Pi _n(U)\cup
^{\Pi _n(W)}\Pi _n(V), A).
$$
But we have seen in \ref{svk} that the morphism
$$
\Pi _n(U)\cup ^{\Pi _n(W)}\Pi _n(V)\rightarrow \Pi _n(X)
$$
is  a trivial cofibration. Thus, to complete the proof it suffices to note
(from \ref{stability}) that for any trivial cofibration $B\rightarrow C$ the
morphism  $\underline{Hom} (C,A)\rightarrow \underline{Hom} (B,A)$ is an equivalence of $n$-categories.

If we take cohomology of  a CW complex $X$ with coefficients in a fibrant
groupoid $A$ then $H(X,A)$ is equivalent to $\Pi _n(\underline{Hom} _{Top}(X,
\langle R,A\rangle )$.  To prove this note that $\Pi _n$ is adjoint to the
realization $\langle R, A \rangle$, which implies on the level of internal
$\underline{Hom}$ that for any $n$-precat $B$ and space $U$,
$$
\underline{Hom}(B, \Pi _n(U)) = \Pi _n (\underline{Hom}_{Top}(\langle R,
B\rangle , U)
$$
where $\underline{Hom}_{Top}$ denotes the compact-open mapping space.
Corolary \ref{equivtoequiv} and the adjointness imply that for any space $U$,
$\Pi _n(U)$ is fibrant.  On the other hand, Tamsamani proves in
\cite{Tamsamani} \S 3 that for any $n$-groupoid $A$ the morphism
$$
A\rightarrow \Pi _n (\langle R, A \rangle )
$$
is an equivalence of $n$-groupoids.  Thus if $A$ is a fibrant $n$-groupoid we
have an equivalence
$$
\underline{Hom}(\Pi _n(X), A) \rightarrow
\underline{Hom}(\Pi _n(X), \Pi _n (\langle R, A \rangle ))
= \Pi _n \underline{Hom}_{Top}(\langle R, \Pi _n(X) \rangle ,
\langle R, A \rangle  ).
$$
On the other hand, again from \cite{Tamsamani} \S 3 we know that
$\langle R, A \rangle $ is $n$-truncated, and there is a space $W$ and diagram
$$
X \leftarrow W \rightarrow \langle R, \Pi _n(X) \rangle
$$
where the left morphism is a weak homotopy equivalence and the right morphism
induces an isomorphism on homotopy groups in degrees $\leq n$.
Thus, under the assumption that $X$ is a CW complex (which allows us to obtain
weak equivalences when we apply $\underline{Hom}_{Top}$) we obtain a diagram of
weak homotopy equivalences
$$
\underline{Hom}_{Top}(X ,
\langle R, A \rangle  )
\rightarrow
\underline{Hom}_{Top}(W ,
\langle R, A \rangle  )
\leftarrow
\underline{Hom}_{Top}(\langle R, \Pi _n(X) \rangle ,
\langle R, A \rangle  ).
$$
Combining with the above we get a diagram of equivalences
$$
\underline{Hom}(\Pi _n(X), A)\rightarrow
\Pi _n \underline{Hom}_{Top}(W ,
\langle R, A \rangle  )
\leftarrow
\Pi _n \underline{Hom}_{Top}(X ,
\langle R, A \rangle  ).
$$
Thus the nonabelian cohomology with coefficients in an $n$-groupoid coincides
with the approach using topological spaces.

Of course even the  nonabelian cohomology with
coefficients in an $n$-category $A$ which isn't a groupoid doesn't really give a
new homotopy invariant since all of the information is contained in the
Poincar\'e $n$-groupoid $\Pi _n(X)$. However, it might give some interesting
special cases to study.

Once the theory of $n$-stacks gets off the ground, we should be able to
interpret $H(X,A)$ as the $n$-category of global sections of the $n$-stack
associated to the constant presheaf $U\mapsto A$ over the site $Site(X)$ of
disjoint unions of open subsets of $X$.

More generally if $\Xx$ is any site and $\underline{A}$ is a presheaf of
$n$-categories over $\Xx$ (or a fibered $n$-category over $\Xx$) then the
$n$-category of global sections of the $n$-stack associated to $\underline{A}$
is the {\em nonabelian cohomology} $H(\Xx, \underline{A})$.

We now  treat the example mentionned in the footnote of the introduction.
Suppose $G$ is a group and $U$ an abelian group, and let $A$ (resp. $B$) be the
strict $1$-category with one object and group of automorphisms $G$ (resp. the
strict $n$-category with one arrow in each degree $<n$ and group $U$ of
automorphisms in degree $n$).  Let $A'$ (resp. $B'$) be fibrant replacements for
$A$ and $B$. We would like to show that the set $T^n\underline{Hom}(A',B')$ is
equal to the group cohomology $H^n(G, U)$.  For the moment, the only way I see
to do this is to pass by topology using the Seifert-Van Kampen theorem.  Let
$X=K(G,1)$ and $Y= K(U,n)$. Then $\Pi _n(X)$ is equivalent to $A$ and $\Pi
_n(Y)$ is equivalent to $B$.  Similarly in the other direction $\langle R, B'
\rangle $ is equivalent to $\langle R , B \rangle$ which in turn is equivalent
to $Y$.  By the above discussion $H(X, B')$ is equivalent to $\Pi
_n(\underline{Hom} (X, Y))$. The truncation $T^nH(X,B')$ is thus equal to $\pi
_0(\underline{Hom} (X,Y))$ which (as is well-known) is $H^n(X, U)= H^n(G, U)$.
But by definition $H(X,B')=\underline{Hom} (\Pi _n(X), B')$ which is equivalent
to $\underline{Hom} (A', B')$. This gives the desired statement.

The above argument is clearly not ideal, since we are looking for a purely
algebraic approach to these types of problems. It seems likely that the
algebraic
techniques of \cite{Quillen} with appropriate small additional lemmas would
permit us to give a purely algebraic proof of the result of the previous
paragraph.

\numero{Comparison}

As pointed out in the introduction, there are many different theories of weak
$n$-categories in the process of becoming reality, and this will pose the
problem of comparison.  As an initial step we give a construction of functors
modeled on Tamsamani's Poincar\'e $n$-groupoid construction. We denote our
``Poincar\'e $n$-category'' functors by $\Upsilon_n$ to avoid confusion with
the $\Pi _n$ (specially on the fact that the $\Upsilon_n$ will not take images
in $n$-groupoids).

This section is only a sketch, with many details of proofs missing. In
particular the following proposed set of axioms for internal model categories
is a preliminary attempt only.

Suppose $\Cc$ is a closed model category with the following additional
properties:
\newline
(IM1)---$\Cc$ admits an internal $\underline{Hom}$;
\newline
(IM2)---If $A$ and $B$ are fibrant and cofibrant objects then $\underline{Hom}
(A,B)$ is fibrant;
\newline
(IM3)---If $A\rightarrow A'$ is a cofibration (resp. trivial cofibration) of
fibrant and cofibrant objects, and if $B'\rightarrow B$ is a fibration
(resp. trivial fibration) of fibrant and cofibrant objects, then
$\underline{Hom}(A', B')\rightarrow \underline{Hom}(A,B)$ is a fibration (resp.
a trivial fibration);
\newline
(IM4)---Internal $\underline{Hom}$ takes cofibrant pushout in the first variable
(resp. fibrant fiber product in the second variable) to fiber product.

We call $\Cc$ an {\em internal closed model category}.

Suppose now that $\Cc$ is an internal closed model category with an inclusion
$i:Cat\subset \Cc$ having the following properties:
\newline
(a)---$i(\emptyset )$ is the initial object and $i(\ast )$
is the final object of
$\Cc$;
(b)---$i$ is compatible with disjoint union;
\newline
(c)---$i$ takes values in the fibrant and cofibrant objects of $\Cc$;
\newline
(d)---$i$ takes the internal $\underline{Hom}$ in $Cat$ to the internal
$\underline{Hom}$ of $\Cc$.
\newline
(e)---Let $I$ denote the category with two objects $0,1$ and one non-identity
morphism from $0$ to $1$.  Let $I^{(m)}$ denote the symmetric product of $m$
copies of $I$, which is the category with objects $0,\ldots , m$ and
one morphism from $i$ to $j$ when $i\leq j$.  Then we require that the morphism
from the $\Cc$-pushout of the diagram
$$
i(I) \leftarrow \ast  \rightarrow i(I) \leftarrow \ast \ldots \ast \rightarrow
i(i)
$$
to  $i(I^{(m)})$ be a cofibrant weak equivalence in $\Cc$.

{\em Remark:}  Our closed model categories $PC_n$ are internal, and (for
$n\geq 1$) have functors $i: Cat\hookrightarrow PC_n$ satisfying properties
(a)--(e) above.

These seem to be reasonable properties to ask of any closed model category
representing a theory of $n$-categories (or $\infty$-categories). However, some
of the properties are of a rather technical nature so it is possible that some
technically slightly different approach to comparison would be needed---the
present section is just a first attempt.

Suppose $(\Cc , i)$ is an internal closed model category with inclusion $i$
having the above properties.  Let $\Cc _f$ denote the subcategory of fibrant
objects. Then for any $n$ we define a functor $\Upsilon _n: \Cc \rightarrow
n-Cat \subset PC_n$, which we call the ``Poincar\'e $n$-category'' functor.
These
functors will have the property that they take weak equivalences in $\Cc _f$ to
equivalences of $n$-categories, and will be compatible with direct products
(hence  with fiber products over sets).

The definition is by induction. First of all,
$\Upsilon _0(X)$ is defined bo be equal to the set of homotopy classes of maps
$\ast\rightarrow X$.  Then, supposing that we have defined $\Upsilon _{n-1}$ we
define for any $X\in \Cc _f$ the simplicial object $U$ of $\Cc$
by: $U_0$ is the set (note that sets are categories so $i$ gives an inclusion
of sets into $\Cc$) of morphisms $\ast\rightarrow X$; and for $x_0,\ldots ,
x_m\in U_0$,
$U_m(x_0,\ldots , x_m)$ is the fiber of
$$
\underline{Hom}(I^{(m)}, X)\rightarrow \underline{Hom}(\{ 0, \ldots , m\} , X)
$$
over the point $(x_0, \ldots , x_m)$. Then $U_m$ is the disjoint union of the
$U_m(x_0,\ldots ,  x_m)$ over all sequences of $x_i \in U_0$.

Axiom IM4 and condition (e) imply that the usual morphism
$$
U_m \rightarrow U_1\times _{U_0} \ldots \times _{U_0}U_1
$$
is a weak equivalence.

With the above notation set
$$
\Upsilon _n(X)_{m/}:= \Upsilon _{n-1}(U_m).
$$
This simplicial $n-1$-category is an $n$-category since $\Upsilon _{n-1}$ is
compatible with direct products and preserves weak equivalences. Note that
$\Upsilon _n$ is obviously compatible with direct products.

One has to prove
that $\Upsilon _n$ preserves weak equivalences (we leave this out for now).

\subnumero{Examples}

Tamsamani's functor $\Pi _n$ is essentially an example of the functor
functor $\Upsilon _n$ for $\Cc = Top$ and $i: Cat \rightarrow Top$ the
realization functor.  Our definition of $\Upsilon _n$ is a
generalization of the definition of (\cite{Tamsamani} \S 3).

For $\Cc = PC_{n'}$ we obtain the functor $\Upsilon _n$.
If $n=n'$ it is essentially the identity; for $n< n'$ is is the truncation
$T^{n'-n}$; and for $n>n'$ it is the induction $Ind ^{n'}_n$. Note however that
the
induction doesn't preserve pushouts, so $\Upsilon _n$ will not necessarily
preserve pushouts in general (where by pushouts here we mean the replacement of
pushouts by weak-equivalent objects of $\Cc _f$).

For any given theory $\Cc$ of $n$-categories satisfying the above properties,
one would like to check that  the functor $\Upsilon_n$ is an equivalence
of homotopy theories in the sense of \cite{Quillen} (or at least that it induces
an isomorphism of localized categories). If $n$ is correctly
chosen to correspond to the level of $\Cc$ then one would try to show that
$\Upsilon _n$ preserves pushouts.

There are examples of $\Cc$ which are not equivalent to $PC_n$, such as
$Top$ or, for example,
the category of ``Segal categories'', i.e. simplicial spaces whose first object
is a set and which satisfy Segal's condition (cf \cite{Tamsamani} \S 3).
Even if we look at Segal categories whose elements are $n-1$-truncated, the
functor $\Upsilon _n$ will go into the $n$-precats $A$ whose $A_{p/}$ are
$n-1$-groupoids, in particular $\Upsilon _n$ will not be essentially surjective.

Similarly, one can imagine looking at a theory $\Cc$ of $n$-categories with
extra structure. For example $Top$ is basically the theory of $n$-categories
where the $i$-morphisms have essential  inverses. Baez and Dolan propose
another type of extra structure of ``adjoints'' rather than inverses, in
relation to topological quantum field theory \cite{BaezDolan}.  It is possible
in this case that $\Cc$ would again be an internal closed model category and
that
we would have a functor $i$. The resulting functor $\Upsilon_n$ would then be
essentially the functor of ``forgetting the extra structure'' and taking the
underlying $n$-category.

A more fundamental example of the above phenomenon will be the closed model
category of $n$-stacks.  This retracts onto that of $n$-categories: the
inclusion being the constant stack functor and the morphism $\Upsilon _n$ being
the global section functor.  Of course in this situation we don't expect
$\Upsilon _n$ to be an equivalence of theories. This example shows that more is
needed than just the above axioms for $\Cc$ in order to prove that the
composition in the other direction is the identity.

\end{document}